\documentclass[a4paper,12pt,leqno]{article}
\usepackage[utf8]{inputenc}
\usepackage{amssymb,amsthm}
\usepackage{amsmath}
\usepackage{mathrsfs}
\usepackage{csquotes}
\usepackage{scalerel}
\usepackage[english]{babel}
\usepackage{arydshln}
\usepackage{shuffle}
\usepackage{hyperref,color}
\usepackage{url}
\usepackage{float} 
\usepackage{filecontents}
\hypersetup{
	colorlinks=true,
	linkcolor=blue,
	filecolor=magenta,
	urlcolor=cyan,}
\usepackage[left=2.3cm,top=3cm,right=2.3cm,bottom=3cm,bindingoffset=0.5cm]{geometry}
 \newcommand\blfootnote[1]{%
 	\begingroup
 	\renewcommand\thefootnote{}\footnote{#1}%
 	\addtocounter{footnote}{-1}%
 	\endgroup
 }
 
 \makeatletter
 \newcommand*{\rom}[1]{\expandafter\@slowromancap\romannumeral #1@}
 \makeatother

\usepackage{titlesec}
\makeatletter
\patchcmd{\ttlh@hang}{\parindent\z@}{\parindent\z@\leavevmode}{}{}
\patchcmd{\ttlh@hang}{\noindent}{}{}{}
\makeatother
\titlelabel{\thetitle.\quad}
\titleformat{\section}[runin]
{\large\bf}{\thesection.\quad\hspace{0.1em}}{0.2em}{}[{\\[0.1em]}]
\titlespacing*{\section}{0cm}{0em}{0em}
\titleformat{\subsection}[runin]
{\bf}{\thesubsection.\quad}{0.1em}{}[{\\[0.1em]}]
\titlespacing*{\subsection}{0cm}{0em}{0em}
\makeatletter
\def\thm@space@setup{%
	\thm@preskip=\parskip \thm@postskip=0pt
}
\makeatother
\usepackage{lineno,hyperref}
\makeatletter
\providecommand*{\boxast}{%
	\mathbin{
		\mathpalette\@boxit{*}%
	}%
}
\newcommand*{\@boxit}[2]{%
	\sbox0{$\m@th#1\Box$}%
	\ifx#1\displaystyle \ht0=\dimexpr\ht0+.05ex\relax \fi
	\ifx#1\textstyle \ht0=\dimexpr\ht0+.05ex\relax \fi
	\ifx#1\scriptstyle \ht0=\dimexpr\ht0+.04ex\relax \fi
	\ifx#1\scriptscriptstyle \ht0=\dimexpr\ht0+.065ex\relax \fi
	\sbox2{$#1\vcenter{}$}
	\rlap{%
		\hbox to \wd0{%
			\hfill
			\raisebox{%
				\dimexpr.5\dimexpr\ht0+\dp0\relax-\ht2\relax
			}{$\m@th#1#2$}%
			\hfill
		}%
	}%
	\Box
}
\makeatother

\date{}
\font\pc=cmcsc10
\newtheorem{thm}{\pc Theorem }

\newtheorem{prop}{\pc Proposition}
\newtheorem{cor}{\pc Corollary}
\newtheorem{lem}{\pc Lemma }
\newtheorem{coj}{\pc Conjecture }

\def\cvirg{\,\raise2pt\hbox{,}}
\author{}
\newlength{\oldparindent}
\newcommand{\myindent}{\hspace{\oldparindent}}
\newtheorem{innercustomgeneric}{\customgenericname}
\providecommand{\customgenericname}{}
\newcommand{\newcustomtheorem}[2]{%
	\newenvironment{#1}[1]
	{%
		\renewcommand\customgenericname{#2}%
		\renewcommand\theinnercustomgeneric{##1}%
		\innercustomgeneric
	}
	{\endinnercustomgeneric}
}

\newcustomtheorem{customthm}{\pc Theorem}
\newcustomtheorem{customlem}{\pc Lemma}
\newcustomtheorem{customcor}{\pc Corollary}
\makeatletter
\@addtoreset{cor}{subsection}
\makeatother

\title{ Multiple zeta values and  multiple Ap{\' e}ry-like sums}
\author{P. Akhilesh}
\newcommand{\Addresses}{{
		\bigskip
		\footnotesize
		P.~Akhilesh, \textsc{Kerala School of Mathematics, kozhikode, Kerala, India.\\Chennai Mathematical Institute, Chennai,  India \\ The first address is present address. }\par\nopagebreak 
		\textit{E-mail }: \texttt{akhi@cmi.ac.in},  \texttt{akhi@ksom.res.in}
	}}

\begin{document}
\maketitle
\begin{abstract}
	In this paper, we formally introduce the notion of Ap{\'e}ry-like sums and we show that every multiple zeta values can be expressed as a $\bf Z$-linear combination of them. We even describe a canonical way to do so. This allows us to put in a new theoretical context several identities scattered in the literature, as well as to discover many new interesting ones. We give in this paper new integral formulas for multiple zeta values and Ap{\'e}ry-like sums. They enable us to give a short direct proof of Zagier's formulas for $\zeta(2,\ldots,2,3,2,\ldots,2)$ as well as of similar ones in the context of Ap{\'e}ry-like sums. The relations between Ap{\'e}ry-like sums themselves still remain rather mysterious, but we get significant results and state some conjectures about their pattern.
\end{abstract}

\vspace{2em}
\section*{Introduction}\label{s0}
\myindent Multiple zeta values are the real numbers
\begin{equation}
\label{EF1}
\zeta({\bf a})= \sum_{n_1>\ldots>n_r>0}n_1^{-a_1}\ldots n_r^{-a_r},
\end{equation}
where ${\bf a}=(a_1, \ldots ,a_r) $ is an {\it admissible composition}, {\it i.e.} a finite sequence of positive integers, with $a_1 \geqslant 2$ when $r\neq 0$.

In this paper, we introduce the {\it multiple Ap{\'e}ry-like sums} defined by
\begin{equation}\label{EF2}
\sigma({\bf a})=\sum_{n_1>\ldots>n_r>0}\left({2 n_1 \atop n_1}\right)^{-1}n_1^{-a_1}\ldots n_r^{-a_r}
\end{equation}
when ${\bf a}\neq\varnothing$ and by $\sigma(\varnothing)=1$. We show that for any admissible composition ${\bf a}$, there exists a finite formal $\bf Z$-linear combination $\sum \lambda_{\bf b} {\bf b}$ of admissible compositions such that
 \begin{equation}\label{EU3}
\zeta({\bf a})=\sum \lambda_{\bf b}\, \sigma({\bf b}).
 \end{equation} 
 The simplest instance of this fact is the identity
\begin{equation}\label{EF3}
\sum_{n=1}^{\infty}\frac{1}{n^2}=3\sum_{n=1}^{\infty}\frac{1}{\left({2n \atop n}\right)n^2}
\end{equation} 
due to Euler, which expresses that $\zeta(2)=3\,\sigma(2)$. Note that multiple Ap{\' e}ry-like sums have the advantage on multiple zeta values to be exponentially quickly convergent.

The  formal $\bf Z$-linear combination $\sum \lambda_{\bf b}\,{\bf b}$ satisfying (\ref{EU3}) is not uniquely determined by ${\bf a}$ in general. But we show in Section \ref{s3s2} that there exists a unique one satisfying the  stronger requirement that
\begin{equation}\label{EF4}
\forall n \in {\bf N}, \quad \zeta({\bf a})_{n,n}=\sum \lambda_{\bf b} \;\sigma({\bf b})_n,
\end{equation}
 where $\bf N$ denotes the set of non-negative integers, $\sigma({\bf b})_n$ is the $n$-tail of $\sigma({\bf b}) $ (defined in Section \ref{s2s1}) and  $\zeta({\bf a})_{n,n}$ is the symmetric $(n,n)$-double tail of $\zeta({\bf a})$ (defined in Section~\ref{s1s8}). We deduce the  unicity  from a strong  linear independence theorem   about tails of multiple Ap{\'e}ry-like sums, stated and proved in  Section \ref{s2s7}.
 
 The unique formal $\bf Z$-linear combination $\sum \lambda_{\bf b}{\bf b}$ associated to $\bf a$ in this way is the same as the one associated to the dual composition $\overline{\bf a}$ (defined in \ref{s1s7}), since $\zeta({\bf a})_{n,n}=\zeta(\overline{\bf a})_{n,n}$ for all $ n \in {\bf N}$ (see \ref{s1s7}, formula (\ref{EF13})). Therefore it depends only on the equivalence class $[{\bf a}]=\{{\bf a},\overline{\bf a}\}$ of ${\bf a}$ up to duality. We write it $\delta([{\bf a}])$ and extend $\delta$ to a $\bf Z$-linear  map from $ {\bf Z}^{(\mathcal{B})}$ to  ${\bf Z}^{(\mathcal{A})}$,
 where $\mathcal{A}$ is the set of admissible compositions and $\mathcal{B}$  the set of admissible compositions up to duality.
 
 We show (see Section  \ref{s3s2}, Remark) that $\delta$ is the unique $\bf Z$-linear map from ${\bf Z}^{(\mathcal{B})}$ to ${\bf Z}^{(\mathcal{A})}$, graded of degree $0$ for the weight, such that $\delta([\varnothing])=\varnothing$ and satisfying the functional equation
 \begin{equation}\label{E7}
 \forall {\bf a}\in \mathcal{A}-\{\varnothing\},\qquad \quad \delta([{\bf a}])^{\rm init}=\delta([{\bf a}^{\rm init}]+[{\bf a}^{\rm mid}]+[{\bf a}^{\rm fin}]),
 \end{equation}
 where the initial, middle and final parts of  admissible compositions are defined in \ref{s1s9}, and   extended to ${\bf Z}^{(\mathcal{A})}$ by $\bf Z$-linearity.

Using these results, we not only provide elegant proofs of various identities between multiple zeta values and Ap{\'e}ry-like sums (due to by Bailey, Borwein, Bradley, Leshchiner, Zagier etc.), but we also discover many new ones, such as for example
 	\begin{equation}
 	\zeta(a)=2\;\sum_{b=3}^{a}\sigma(b,\underbrace{1,\ldots ,1}_{a-b})+3\;\sigma(2,\underbrace{1,\ldots ,1}_{a-2}),
 	\end{equation}
which extends Euler identity (\ref{EF3}) to all integers  $a\geqslant 2$ (see Section \ref{S3S3}, Corollary \ref{CT2} of Theorem~\ref{Th8}).

In Section \ref{s3as2}, we give an explicit combinatorial expression of $\delta([{\bf a}])$ for any admissible composition $\bf a$, and we use it in Section \ref{s4as4} to compute of $\delta([{\bf a}])$ in cases that were not easy to handle by using only the functional equation (\ref{E7}).

In our previous paper \cite{akhi}, we had discovered non zero elements $\sum \lambda_{\bf a}{\bf a}$ in ${\bf Z}^{(\mathcal{B})}$ for which
\begin{equation}\label{EF5}
\forall n \in {\bf N}, \;\; \sum \lambda_{\bf a}\zeta({\bf a})_{n,n}=0,
\end{equation} 
  holds. The first  one was found in weight $6$ (see Sections \ref{s1s12} and \ref{s4s6}).

In Section \ref{s4s1}, we study the $\bf Z$-module ${\rm M}$  of finite formal ${\bf Z}$-linear combinations \mbox{ $\sum \lambda_{\bf a}[{\bf a}]\in {\bf Z}^{(\mathcal{B})}$}  for which $(\ref{EF5})$ holds. We prove that it is the kernel of the  map $\delta$, hence that it is graded by the weight. Note that  this is a non trivial result, since it is not known at present whether the $\bf Z$-module of finite formal $\bf Z$-linear combinations\mbox{ $\sum \lambda_{\bf a}[{\bf a}]\in {\bf Z}^{(\mathcal{B})}$}  such that $\sum \lambda_{\bf a}\zeta({\bf a})=0$ is graded.

We give in Section \ref{s4s4} a second independent description of the  module $\rm M$. We use it to prove that the relation $ \flqq{} \, \forall n \in {\bf N}, \,\zeta({\bf{a}})_{n,n}=\zeta({\bf{b}})_{n,n}\,\frqq{} $ between two admissible compositions $\bf a$ and $\bf b$ implies that $\bf b$ is equal to $\bf a$ or $\overline{\bf a}$ (see Section \ref{s4s5}). Here again, it is not known whether the single hypothesis $\zeta({\bf a})=\zeta({\bf b})$ suffices to imply the same conclusion.

In the course of this work, we have established new integral formulas for multiple zeta values and Ap{\'e}ry-like sums (see Section \ref{s2s6}). They seem well suited to study F.~Brown's motivic depth filtration (see the remark in Section \ref{s4s3}). They enable us to give a very short and direct proof of D. Zagier's  evaluation of $\zeta(2,\ldots,2,3,2\ldots,2)$ in terms of  values of the Riemman zeta function at odd integers $\geqslant3$ (see Section \ref{s2s8a}). We also get analogous expressions for the Ap{\'e}ry-like sums $\sigma(2,\ldots,2,3,2\ldots,2)$ and $\zeta(2,\ldots,2,1,2\ldots,2)$ (see Section \ref{s2s8}). The $\bf Q$-subspace of $\bf R$ generated by the Ap{\'e}ry-like sums of this form and of weight $k$ is then shown to be equal to the $\bf Q$-subspace of $\bf R$ generated by the real numbers $\pi^{k-2r-1}\zeta(2r+1)$ and $\pi^{k-2r}L(2r,\chi)$, where $1\leqslant r\leqslant \frac{k-1}{2}$ and $\chi$ is the principal Dirichlet character of conductor $3$ (see Section \ref{s2s9}, Theorem \ref{th7}).

The $\bf Q$-linear relations between Ap{\'e}ry-like sums themselves and  the combinatorial behaviour of the map $\delta$ remain rather mysterious. We start their study in the last Section \ref{s5} of the paper, where we both obtain significant partial results and state some  conjectures based on numerical experimentation. We intend to pursue this study in a forthcoming paper.

\vspace{2em}
\section{Multiple zeta values and their double tails}\label{s1}
\subsection{Compositions}\label{s1s1}
\myindent A finite sequence ${\bf a}=(a_1,\ldots,a_r)$ of positive integers  is called {\it a composition}. The integer $r$ is called {\it the depth} of $\bf a$. The integer $k=a_1+\ldots+a_r$ is called {\it the weight} of $\bf a$ and is denoted by $|{\bf a}|$. We denote by $\mathscr{C}$ the set  of all compositions and by $\mathscr{C}^{\ast}$ the set $\mathscr{C}-\{\varnothing\}$ of non-empty ones.

The composition~$\bf a$ is said to be {\it admissible} if either it is empty or if $a_1\geqslant 2$.
We denote by  $\mathcal{A}$ the set of admissible compositions and by  $\mathcal{A}^*$ the set  $\mathcal{A} - \{ \varnothing\}$. For each integer $k \geqslant 0$, we denote by $\mathcal{A}_k$ the set of admissible compositions of weight $k$.
\vspace{2em}
\subsection{Multiple zeta values}\label{s1s2}
\myindent To each admissible composition ${\bf a}=(a_1,\ldots,a_r)$, one associates a real number~$\zeta(\bf a)$. It is defined by the convergent series
\begin{equation}
\label{E1}
\zeta({\bf a})= \sum_{n_1>\ldots>n_r>0}n_1^{-a_1}\ldots n_r^{-a_r}.
\end{equation}
We have in particular $\zeta(\varnothing)=1$ when $r=0$. These numbers are called {\it multiple zeta values} or
{\it Euler-Zagier numbers}.
\vspace{2em}
\subsection{Binary word associated to a composition }\label{s1s3}
\myindent A {\it binary word}  is by definition a word $w$ constructed on the alphabet $\{0,1\}$. Its letters are  called {\it bits}. The number of bits of $w$ is called {\it the weight} of $w$ and denoted by~$|w|$. The number of bits of $w$ equal to $1$ is called {\it the depth} of $w$. To any composition ${\bf a}=(a_1,\ldots,a_r)$, one associates the binary word
\begin{equation}
\label{E2}
{\bf w({\bf{a}})}= \{0\}_{a_1-1}1\ldots \{0\}_{a_r-1}1,
\end{equation}
where for each integer $u\geqslant 0$, $\{0\}_u$ denotes the binary word consisting of $u$ bits equal to~$0$. When $\bf a$ is the empty composition, ${\bf w}({\bf a})$ is the empty word.
The weight of $\bf w(\bf a)$ is equal to the weight of $\bf a$ and its depth to the depth of $\bf a$.

The map $\bf w$ is a bijection from the set of compositions onto the set of binary words non ending in $0$. Words corresponding to admissible compositions will be called admissible. Non empty compositions corresponds to binary words ending in $1$, and non empty admissible compositions to  binary words starting with $0$ and ending in $1$.
\vspace{2em}
\subsection{Kontsevich integral formula for multiple zeta values}\label{s1s4}
\myindent M. Kontsevich noticed that for each  admissible composition $\bf a$, the multiple zeta value $\zeta(\bf a)$ can be written as an iterated integral. More precisely, if $w=\varepsilon_1\ldots\varepsilon_k$ denotes  the associated binary word $\bf w(\bf a)$, we have
\begin{equation}
\label{E3}
\zeta({\bf{a}})={\rm{It}}\int_0^1(\omega_{\varepsilon_1},\ldots,\omega_{\varepsilon_k})= \int_{1>t_1>\ldots>t_k>0} f_{\varepsilon _1}(t_1)\ldots f_{\varepsilon _k}(t_k)dt_1\ldots dt_k,
\end{equation}
where $\omega_{\varepsilon}=f_{\varepsilon}(t)dt$, with  $f_0(t)={1\over t}$ and $f_1(t)={1\over 1-t}\cdot$
\vspace{2em}
\subsection{Tails and double tails of multiple zeta values}\label{s1s5}
\myindent Let ${\bf a}=(a_1,\ldots,a_r)$ be a non-empty admissible composition and $n$ be a non-negative integer. The $n$-tail of
$\zeta(\bf a)$ is by definition
 the sum of the series
$$\sum_{n_1>\ldots>n_r>n}n_1^{-a_1}\ldots n_r^{-a_r}.$$

This $n$-tail can be written as the iterated integral
\begin{equation}
\label{E6}
{\rm It}\int_0^1(\omega_{\varepsilon_1},\ldots,t^n\omega_{\varepsilon_k})=\int_{1>t_1>\ldots>t_k>0}f_{\varepsilon _1}(t_1)\ldots
f_{\varepsilon _k}(t_k)t_k^ndt_1\ldots dt_k,
\end{equation}
where $\varepsilon_1\ldots\varepsilon_k$ is the binary word $\bf w(\bf a)$ (\cite{akhi}, Theorem 1).

In our paper \cite{akhi}, we have defined  for any non empty admissible composition $\bf{a}$ and any integers $m,n \in {\bf N}$,  the double tails $\zeta({\bf{a}})_{m,n}$ by the iterated integrals
\begin{eqnarray}
\label{E3a}
\zeta({\bf{a}})_{m,n}&=&{\rm{It}}\int_0^1((1-t)^m\omega_{\varepsilon_1},\ldots,t^n\omega_{\varepsilon_k})\\ \nonumber&=& \int_{1>t_1>\ldots>t_k>0} (1-t_1)^m f_{\varepsilon _1}(t_1)\ldots f_{\varepsilon _k}(t_k)t_k^ndt_1\ldots dt_k,
\end{eqnarray}
and we have given in In Theorem 1 of ${ loc.\;cit.}$  the following series expressions of these double tails:
\begin{equation}
\label{E4}
\zeta({\bf a})_{m,n} =\sum_{n_1>\ldots>n_r>n}\left({n_1+m\atop m}\right)^{-1}n_1^{-a_1}\ldots n_r^{-a_r}.
\end{equation} 
It is convenient to extend this definition to the  empty composition by the convention $\sigma(\varnothing)_{m,n}=\left({m+n \atop m}\right)^{-1}$.
\vspace{2em}
\subsection{Bounds for double tails of multiple zeta values}\label{s1s6}
\myindent Let
${\bf a}=(a_1,\ldots,a_r)$ be a non-empty admissible composition and let $m,n\in\bf N$. From the expressions
(\ref{E3}) and (\ref{E3a}) of $\zeta(\bf a)$ and $\zeta({\bf a})_{m,n}$ as iterated integrals, one deduces in \cite{akhi}, Section 4, that 
\begin{equation}\label{EF12}
\zeta({\bf a})_{m,n}\leqslant \frac{m^m n^n}{(m+n)^{m+n}}\zeta({\bf a}).
\end{equation}
Note that we have $\zeta({\bf a})\leqslant \frac{\pi^2}{6}$ (\cite{akhi}, Theorem 3).
\vspace{2em}
\subsection{Duality}\label{s1s7}
\myindent Let $w=\varepsilon_1\ldots\varepsilon_k$ be a binary word. Its {\it dual word} is defined to be $\overline{w}=\overline{\varepsilon}_k\ldots\overline{\varepsilon}_1$, where $\overline{0}=1$ and $\overline{1}=0$. When $w$ is admissible, so is $\overline{w}$. We can therefore define
the {\it dual composition} of an admissible composition $\bf a$ to be the admissible composition $\overline{\bf a}$ such that $\bf w(\overline{\bf a})$ is dual to $\bf w(\bf a)$. When $\bf a$ has weight $k$ and depth $r$, $\overline{\bf a}$ has weight $k$ and depth~$k-r$.

In \cite{akhi}, Theorem 2, we noticed that, by the change of variables $t_i\mapsto 1-t_{k+1-i}$ in the integral (\ref{E3a}), one gets the following {\it duality relation}
\begin{equation}
\label{EF13}
\zeta({\bf a})_{m,n}=\zeta(\overline{\bf a})_{n,m}.
\end{equation}

We denote by    $\mathcal{B}$ the set of admissible compositions up to duality, {\it i.e} the quotient of  $\mathcal{A}$ by the equivalence relation~$\sim$ for which ${\bf{a}}\sim {\bf{b}}$ if and only if ${\bf{b}}$ is equal to ${\bf{a}}$ or $\overline{{\bf{a}}}$. The equivalence class $\{ {\bf a},\;\overline{\bf a}\} $ of an admissible composition ${\bf a}=(a_1, a_2, \ldots , a_r)$ will be denoted by $[{\bf a}]$  or $[a_1, a_2, \ldots , a_r]$. We denote by  $\mathcal{B}^*$ the set  $\mathcal{B} - \{ [\varnothing]\}$. For each integer $k \geqslant 0$, we denote by $\mathcal{B}_k$ the set of equivalence classes of admissible compositions of weight $k$. 

We denote by $\ell\rightarrow [\ell]$ the $\bf Z$-linear map from ${\bf Z}^{(\mathcal{A})}$ to ${\bf Z}^{(\mathcal{B})}$  extending the map ${\bf a}\rightarrow [{\bf a}]$ from $\mathcal{A}$ to $\mathcal{B}$.
\vspace{2em}
\subsection{Symmetric double tails of multiple zeta values}\label{s1s8}
  \myindent Let $\bf a$ be an admissible composition. We call {\it symmetric double tails of} $\zeta({\bf a})$  the double tails of the form $\zeta({\bf a})_{n,n}$ for $n \geqslant 0$. For them, the duality relations $(\ref{EF13})$ become  
  \begin{equation}
  \label{EF14}
  \zeta({\bf a})_{n,n}=\zeta(\overline{\bf a})_{n,n}.
  \end{equation}
  The map ${\bf a} \rightarrow \zeta({\bf a})_{n,n}$ therefore defines by passing to the quotient a map from $\mathcal{B}$ to $\bf R$. We still denote by $\ell \rightarrow \zeta(\ell)_{n,n}$ its $\bf Z$-linear extension from ${\bf Z}^{(\mathcal{B})}$ to $\bf R$

  When $\bf a$ is non-empty, the upper bounds $(\ref{EF12})$ become, for symmetric double tails,
  \begin{equation}
  \label{E12}
  \zeta({\bf a})_{n,n}\leqslant 2^{-2n}\zeta({\bf a})\leqslant 2^{-2n}\frac{\pi^2}{6}.
  \end{equation}
  \vspace{2em}
\subsection{Initial, middle and final part of an admissible composition}\label{s1s9}
\medskip
\myindent Let ${\bf a}=(a_1,a_2, \ldots, a_r)$ be a  composition. We define  the initial part  ${\bf a}^{\rm init}$ of $\bf a$ to be the  composition
$ (a_1,a_2, \ldots,a_{r-1})$ when $r \geqslant 1$ and $\varnothing$ when $r=0$. 

{\it In the remaining part of this subsection, we assume $\bf a$ to be admissible}. Then ${\bf a}^{\rm init}$ is admissible. We define the final part ${\bf a}^{\rm fin}$ of $\bf a$ to be the admissible composition dual to $\overline{\bf a}^{\rm init}$. Note that it is generally not equal to $(a_2,\ldots,a_r)$ when $r\geqslant 1$. In fact, we have
\begin{equation}
{\bf a}^{\rm fin}=\begin{cases}
\varnothing & {\rm if} \; {\bf a}=\varnothing \; {\rm or }\; {\bf a}=\smash{(2,\underbrace{1,\ldots,1}_{r-1})} {\text{ with }} r\geqslant1,\\
(a_1-1,\ldots,a_r) &\; {\rm if}\; r\geqslant 1\; {\rm and} \;a_1 \geqslant 3,\\
(a_i,\ldots,a_r)&\; {\rm if } \;{\bf a}=(2,\underbrace{1,\ldots,1}_{i-2},a_i,\ldots,a_r)\; {\rm with }\; 2\leqslant i \leqslant r\; {\rm and}\; a_i \geqslant 2.
\end{cases}
\end{equation}

We define the middle part $\bf a^{\rm mid}$ of $\bf a$ to be the initial part of $\bf a^{\rm fin}$. We can deduce from the formulas given above that ${\bf a}^{\rm mid}$ is also the final part of ${\bf a}^{\rm init}$, or equivalently, that $ \overline{\bf a}^{\rm mid}$ is the dual composition of ${\bf a}^{\rm mid}$. This also follows from the next remark.
\medskip

{\it Remark.-} Let $w$ denote the binary word ${\bf w}({\bf a})$ associated to the admissible composition~$\bf a$. When $\bf a$ is non-empty, $w$ can be written in a unique way in the form $0\{1\}_{b-1}v\{0\}_{a-1}1$, where $a$ and $b$ are positive integers and $v$ is an admissible word. We defined in \cite{akhi} the initial, middle and final parts of the word $w$ to be the binary words 
\begin{equation}
w^{\rm init}=0\{1\}_{b-1}v, \quad w^{\rm mid}=v, \quad w^{\rm fin}=v\{0\}_{a-1}1.
\end{equation}

When $r\geqslant 2$, we have $w^{\rm init}={\bf w}({\bf a}^{\rm init})$ and $a=|{\bf a}|-|{\bf a}^{\rm init}|$, whereas when $r=1$, we have $w^{\rm init}=0$, ${\bf a}^{\rm init }= \varnothing$ and $|{\bf a}|-|{\bf a}^{\rm init}|=a+1$. Similarly when the depth of $\overline{\bf a}$ is $\geqslant 2$, we have $w^{\rm fin}={\bf w}({\bf a}^{\rm fin})$ and $b=|{\bf a}|-|{\bf a}^{\rm fin}|$, whereas when it is $1$, we have
 $w^{\rm fin}=1$, ${\bf a}^{\rm fin}= \varnothing$ and $|{\bf a}|-|{\bf a}^{\rm fin}|= b+1$. In all cases,
  $w^{\rm mid}$ is equal to ${\bf w}({\bf a}^{\rm mid})$ and $a+b$ to $|{\bf a}|-|{\bf a}^{\rm mid}|$.
   Since $\overline{w}^{\rm mid}$ is dual to $w^{\rm mid}$, $\overline{\bf a}^{\rm mid}$ is dual to ${\bf a}^{\rm mid}$.

\bigskip

We extend the maps ${\bf a}\rightarrow {\bf a}^{\rm init}$, ${\bf a}\rightarrow {\bf a}^{\rm fin}$ and ${\bf a}\rightarrow {\bf a}^{\rm mid}$ from $\mathcal{A}$ to $\mathcal{A}$ to $\bf Z$-linear maps from ${\bf Z}^{(\mathcal{A})}$ to ${\bf Z}^{(\mathcal{A})}$,  denoted by $\ell \rightarrow \ell^{\rm init}$, $\ell \rightarrow \ell^{\rm fin}$ and $\ell \rightarrow \ell^{\rm mid}$.

\vspace{2em}

\subsection{Recurrence relations satisfied by symmetric double tails }\label{s1s10}
\begin{prop}\label{P1}
Let ${\bf a}$ be a non-empty admissible composition. For each integer $n \geqslant 1$, $\zeta({\bf a})_{n-1,n-1}-\zeta({\bf a})_{n,n}$ is equal to
\begin{equation}
n^{|{\bf a}^{\rm init}|-|{\bf a}|}\zeta({\bf a}^{\rm init})_{n,n} +n^{|{\bf a}^{\rm mid}|-|{\bf a}|}\zeta({\bf a}^{\rm mid})_{n,n}+n^{|{\bf a}^{\rm fin}|-|{\bf a}|}\zeta({\bf a}^{\rm fin})_{n,n}.
\end{equation}
\end{prop}

We will check that Proposition \ref{P1} is a reformulation of Theorem 5 in \cite{akhi}, that was stated as follows:
\begin{equation}
\hspace{0.5cm}
\zeta(w)_{n-1,n-1}-\zeta(w)_{n,n}=n^{-a}\zeta(w^{\rm init})_{n,n}+n^{-a-b}\zeta(w^{\rm mid})_{n,n}+n^{-b}\zeta(w^{\rm fin})_{n,n}\,,
\end{equation}
the notations being those of the remark in \ref{s1s9}, and  the multiple zeta values involved being those defined in \cite{akhi}, definition 2. Indeed, by using the \mbox{remark} of \ref{s1s9}, one checks that $n^{-a}\zeta(w^{\rm init})_{n,n}$ is always equal to $n^{|{\bf a}^{\rm init}|-|{\bf a}|}\zeta({\bf a}^{\rm init})_{n,n},\;$~: when the depth of $\bf a$ is at least $2$, we have $|{\bf a}|-|{\bf a}^{\rm init}|=a$ and $w^{\rm init}={\bf w}({\bf a}^{\rm init}),$ hence $\zeta(w^{\rm init})_{n,n}=\zeta({\bf a}^{\rm init})_{n,n},$ whereas when the depth is $1$, we have $|{\bf a}|-|{\bf a}^{\rm init}|=a+1,w^{\rm init }=0 $ and $ {\bf a}^{\rm init}=\varnothing,$ hence $\zeta(w^{\rm init})=\frac{1}{n}\zeta({\bf a}^{\rm init})$. One checks similarly that $n^{-a-b}\zeta(w^{\rm mid})_{n,n}$ is equal to $n^{|{\bf a}^{\rm mid}|-|{\bf a}|}\zeta({\bf a}^{\rm mid})_{n,n}$ and $n^{-b}\zeta(w^{\rm fin})_{n,n}$ to $n^{|{\bf a}^{\rm fin}|-|{\bf a}|}\zeta({\bf a}^{\rm fin})_{n,n}$.  
\vspace{2em}
\subsection{The linear map $\alpha : {\bf Z}^{(\mathcal{B}^*)}\rightarrow {\bf Z}^{(\mathcal{B})}$}\label{s1s11}

\myindent Let $\bf a$ be a non-empty admissible composition. The element $[{\bf a}^{\rm init}]+[{\bf a}^{\rm mid}]+[{\bf a}^{\rm fin}]$ of~${\bf Z}^{(\mathcal{B})} $ remains unchanged when we replace $\bf a$ by its dual composition $\overline{\bf a}$. Indeed, we have $ \;\overline{\bf a}^{\rm init}=\overline{\bf a^{\rm fin}},\;\; \overline{\bf a}^{\rm mid}=\overline{\bf a^{\rm mid}}\;$ and  $\;\overline{\bf a}^{\rm fin}=\overline{\bf a^{\rm init}}$ by \ref{s1s9}, hence $ [\overline{\bf a}^{\rm init}]=[{\bf a}^{\rm fin}],\; [\overline{\bf a}^{\rm mid}]=[{\bf a}^{\rm mid}]$ and $[\overline{\bf a}^{\rm fin}]=[{\bf a}^{\rm init}]$.

Therefore there is a unique $\bf Z$-linear map
\begin{equation}
\alpha : {\bf Z}^{(\mathcal{B}^*)}\rightarrow {\bf Z}^{(\mathcal{B})}
\end{equation}
that maps $[{\bf a}]$ to $[{\bf a}^{\rm init}]+[{\bf a}^{\rm mid}]+[{\bf a}^{\rm fin}]$ for any non-empty admissible composition $\bf a$.

For each integer $k \geqslant 1$, $\alpha$ induces a map
\begin{equation}
\alpha_k : {\bf Z}^{\mathcal{B}_k}\rightarrow \bigoplus_{0\leqslant k'<k}{\bf Z}^{\mathcal{B}_{k'}},
\end{equation}
and we denote by $\alpha_{k',k}:{\bf Z}^{\mathcal{B}_k}\rightarrow {\bf Z}^{\mathcal{B}_{k'}}$ its $k'$-component.
\vspace{2em}

\subsection{The kernel of $\alpha_k$}\label{s1s12}
\myindent Let $k \geqslant 1$ be an integer.  For each $\ell \in {\bf Z}^{\mathcal{B}_k}$ and each integer $n \geqslant 1$, we have, with the notations of Section \ref{s1s8},

\begin{equation}\label{EF20}
\zeta(\ell)_{n-1,n-1}-\zeta(\ell)_{n,n}=\sum_{0\leqslant k'<k}n^{k'-k}\zeta(\alpha_{k',k}(\ell))_{n,n}\,.
\end{equation}
 This follows indeed by $\bf Z$-linearity from  the recurrence relations stated in Proposition \ref{P1} of Section \ref{s1s10}.
 
 \medskip
 
\begin{prop}
For each integer $k\geqslant 1$ and each $\ell \in {\rm Ker}(\alpha_k)$, we have $\zeta(\ell)_{n,n}=0$ for all $n\geqslant 0$.
\end{prop}

Since $\ell$ belongs to the kernel of $\alpha_k$, and hence of $\alpha_{k',k}$ for $0 \leqslant k' < k$, we have $\zeta(\ell)_{n-1,n-1}=\zeta(\ell)_{n,n}$ for each integer $n\geqslant 1$ by formula (\ref{EF20}).
Since we know by Section~\ref{s1s8} that $\zeta(\ell)_{n,n}$ tends to $0$ when $n$ goes to $+\infty$, it follows that $\zeta(\ell)_{n,n}=0$ for all $n \geqslant 0$.
\medskip

We have computed the kernel of $\alpha_k$ for $1 \leqslant k \leqslant 16$. Its rank is given in the following table.
\begin{table}[h]
\centering  
\begin{tabular}{|c|c|c|c|c|c|c|c|c|c|c|c|c|c|c|c|c|}
\hline
$k$  & 1 & 2 & 3 & 4 & 5 & 6 & 7 & 8 & 9 & 10 & 11 & 12 & 13 & 14 & 15 & 16 \\
\hline
 ${\rm rank(Ker}(\alpha_k))$ & 0 & 0 & 0 & 0 & 0 & 1 & 0 & 3 & 2 & 9 & 10 & 31 & 42 & 105 & 165& 364\\
\hline
\end{tabular} 
\end{table}

We see that $k=6$ is the first weight for which $\alpha_k$ is not injective. The kernel of $\alpha_6$ has rank $1$ over $\bf Z$.
It is generated by
\begin{eqnarray}\nonumber
2[6]-2[5,1]+4[4,2]+[4,1,1]+[3,3]-2[3,2,1] -[3,1,2]-2[2,4]+[2,2,2]-2[2,1,3].
\end{eqnarray}
This implies that we have, for all $n\geqslant 0$,
\begin{eqnarray}\nonumber
2\zeta(6)_{n,n}-2\zeta(5,1)_{n,n}+4\zeta(4,2)_{n,n}+\zeta(4,1,1)_{n,n}+\zeta(3,3)_{n,n}-2\zeta(3,2,1)_{n,n}\\ \nonumber-\zeta(3,1,2)_{n,n}-2\zeta(2,4)_{n,n}+\zeta(2,2,2)_{n,n}-2\zeta(2,1,3)_{n,n}=0,
\end{eqnarray}
and moreover this relation does not follow from the duality relations stated in Section~\ref{s1s7}.
\medskip

{\it Remarks.-}  1) One could ask whether conversely any $ \ell \in {\bf Z}^{\mathcal{B}_k}$ such that  $$ \forall n \in {\bf N},\zeta(\ell)_{n,n}=0 $$  belongs to the kernel of $\alpha$. We shall see in Section \ref{s4s6} that this is not true, already when $k=8$.

2) It is neither true that the kernel of the map $\alpha$ is graded, {\it i.e.} is the sum of the kernels of the maps $\alpha_k$ for $k \geqslant 1$, nor that this kernel is annihilated by $\zeta$. For example $\ell=(2,1,3)-2(2,1,2)+(2,2)$ belongs to ker($\alpha$) but not to $\bigoplus_{k\geqslant 1}{\rm Ker}(\alpha_k)$, and we have $\zeta(\ell)\neq 0$.

3) We do not have at present any closed formula for the rank of  ${\rm Ker}(\alpha_k)$ as a function of $k$.
\vspace{2em}
\section{Multiple Ap{\'e}ry-like sums }\label{s2}
\subsection{Definition of multiple Ap{\'e}ry-like functions and their tails }\label{s2s1}
\myindent Let $r \geqslant 0$ be an integer. We define a multiple Ap{\'e}ry-like function $\sigma: {\bf C}^{r}\rightarrow {\bf C}$ by the convergent series
\begin{equation}\label{EF21}
\sigma(s_1,\ldots,s_r)=\sum_{n_1>\ldots>n_r>0}\left({2n_1 \atop n_1}\right)^{-1}n_1^{-s_1}\ldots n_r^{-s_r}.
\end{equation}

Since $|n_1^{-s_1}\ldots n_r^{-s_r}| = n_1^{-{\rm Re}(s_1)}\ldots n_r^{-{\rm Re}(s_r)}$, and since $(t_1,\ldots ,t_r) \rightarrow n_1^{-t_1}\ldots n_r^{-t_r}$ is a decreasing function of each coordinate in ${\bf R}^{r}$, the infinite sum in the left-hand side of (\ref{EF21}) is normally convergent on any  subset of ${\bf C}^{r}$, in which the real parts of the complex numbers $s_i$ are bounded below. Therefore, $\sigma$ is an entire function on ${\bf C}^r$.

 For each integer $n \geqslant 0$, we call {\it $n$-tail of the multiple Ap{\'e}ry-like function} $\sigma: {\bf C}^{r}\rightarrow {\bf C}$ the function ${\bf s} \rightarrow \sigma({\bf s})_{n}$ on ${\bf C}^r$ defined by
\begin{equation}\label{EM23}
\sigma({\bf s})_n=\sum_{n_1>\ldots >n_r>n}\left({2n_1 \atop n_1}\right)^{-1}n_1^{-s_1}\ldots n_r^{-s_r}
\end{equation}
when $r \geqslant 1$, and by the  convention $\sigma(\varnothing)_n=\left({2n \atop n}\right)^{-1}$ when $r=0$.

\vspace{2em}

\subsection{Multiple Ap{\'e}ry-like   sums and their integral representations}\label{s2s6}
\myindent
We call multiple Ap{\'e}ry-like sum associated to a composition ${\bf a}=(a_1, \ldots ,a_r)$ the real number
\begin{equation}
\sigma({\bf a})=\sum_{n_1>\ldots >n_r>0}\left({2n_1 \atop n_1}\right)^{-1}n_1^{-a_1}\ldots n_r^{-a_r}
\end{equation}
when $r\geqslant 1$ and $\sigma(\varnothing)=1$ when $r=0.$ 
\medskip

\begin{thm}
Let ${\bf a}=(a_1, \ldots ,a_r)$ be a non-empty composition and 
 $w=\varepsilon_1\ldots\varepsilon_k$ be its associated binary word $\bf w(\bf a)$. Then we have
\begin{equation}\label{EU28}
\begin{aligned}
\sigma({\bf{a}})&={\rm{It}} \int_{0}^{\frac{1}{4}} (\frac{\omega_{\varepsilon_1}}{\sqrt{1-4t}}, \omega_{\varepsilon_2}\ldots,\omega_{\varepsilon_k})\\  &= \int_{\frac{1}{4}>t_1>\ldots>t_k>0} \frac{f_{\varepsilon _1}(t_1)}{\sqrt{1-4t_1}}f_{\varepsilon _2}(t_2)\ldots f_{\varepsilon _k}(t_k)dt_1\ldots dt_k,
\end{aligned}
\end{equation}
where $\omega_{\varepsilon}=f_{\varepsilon}(t)dt$, with  $f_0(t)={1\over t}$ and $f_1(t)={1\over 1-t}\cdot$
\end{thm}
By definition of the iterated integral, we have
\begin{eqnarray}\label{EF30}
{\rm{It}} \int_{0}^{\frac{1}{4}} (\frac{\omega_{\varepsilon_1}}{\sqrt{1-4t}}, \omega_{\varepsilon_2},\ldots,\omega_{\varepsilon_k})  &=& \int_{0}^{\frac{1}{4}} \frac{\omega_{\varepsilon _1}}{\sqrt{1-4t}}\left({\rm It}\int_{0}^{t}(\omega_{\varepsilon_2},\ldots,\omega_{\varepsilon_k})\right)\\ \nonumber &=&\int_{0}^{\frac{1}{4}}\frac{d{\rm Li}_{\bf a}(t)}{\sqrt{1-4t}},
\end{eqnarray}
where ${\rm Li}_{\bf a}$ denotes the one-variable multiple polylogarithm associated to the composition $\bf a$, defined on the open unit disc $\bf D$ of $\bf C$ by the infinite sum
\begin{equation}
{\rm Li}_{\bf a}(z)=\sum_{n_1>\ldots >n_r>0}n_1^{-a_1}\ldots n_r^{-a_r}z^{n_1},
\end{equation}
which is normally convergent on any compact subset of $\bf D$. We can therefore replace $d{\rm Li}_{\bf a}(t)$ in (\ref{EF30}) by the infinite sum
\begin{equation}
\sum_{n_1>n_2>\ldots >n_r>0}n_1^{1-a_1}n_2^{-a_2}\ldots n_r^{-a_r}t^{n_1-1}dt
\end{equation}
which converges  normally  on $[0,\frac{1}{4}]$, and then exchange the sum and integral.

We conclude by noting that the change of variable $t=x(1-x)$ yields for $n\geqslant 1$
\begin{equation}
\begin{aligned}
\hspace{3em} \int_{0}^{\frac{1}{4}}\frac{t^{n-1}}{\sqrt{1-4t}}dt=&\int_{0}^{\frac{1}{2}}\frac{x^{n-1}(1-x)^{n-1}}{1-2x}(1-2x)dx=\frac{1}{2}\int_{0}^{1}x^{n-1}(1-x)^{n-1}dx\\ =&\frac{\Gamma(n)\;\Gamma(n)}{2\;\Gamma(2n)}=\left({2n \atop n}\right)^{-1}\frac{1}{n}\,.
\end{aligned}
\end{equation}

\medskip
\begin{thm}
Let ${\bf a}$ be a non-empty composition and 
 $w=\varepsilon_1\ldots\varepsilon_k$ be its  associated binary word $\bf w(\bf a)$. Then we have, for each integer $n \geqslant 0$,
\begin{eqnarray}\label{EU33}
\sigma({\bf{a}})_{n}&=&{\rm{It}} \int_{0}^{\frac{1}{4}} (\frac{\omega_{\varepsilon_1}}{\sqrt{1-4t}}, \omega_{\varepsilon_2},\ldots,t^n\omega_{\varepsilon_k}).
\end{eqnarray}

\end{thm}

The proof is analogous to the proof of the previous theorem. We just have to note that $\omega_{\varepsilon_1}{\rm It}\int_{0}^{t}(\omega_{\varepsilon_2}\ldots,t^n\omega_{\varepsilon_k})$ is the differential of the $n$-tail of the multiple polylogarithm~$Li_{\bf a}$.
\vspace{2em}
\subsection{Bounds for  multiple Ap{\'e}ry-like sums and their tails}\label{s2s2}
\begin{prop}
Let ${\bf a}=(a_1,\ldots,a_r)$ be a non-empty composition. We have
	\begin{equation}
	\sigma({\bf a})\leqslant \frac{({\rm log}(4/3))^{r-1}}{(r-1)!} \frac{\pi}{3 \sqrt{3}}
	\end{equation}
	and $\sigma({\bf a})_n \leqslant4^{-n}\sigma({\bf a})$ for each integer $n \geqslant 0$.
\end{prop}

We have indeed
\begin{equation}
\begin{aligned}
\sigma({\bf a})&\leqslant \sigma(\underbrace{1,\ldots,1}_{r}) = {\rm It}\int_{0}^{\frac{1}{4}}\left(\frac{dt}{ \sqrt{1-4t}(1-t)},\underbrace{\frac{dt}{1-t}, \ldots , \frac{dt}{1-t}}_{r-1}\right)
\\  &= \int_{0}^{\frac{1}{4}}\frac{1}{\sqrt{1-4t}(1-t)}\frac{(-{\rm log}(1-t))^{r-1}}{(r-1)!}dt
 \\&\leqslant \frac{({\rm log}(4/3))^{r-1}}{(r-1)!}\int_{0}^{\frac{1}{4}}\frac{1}{\sqrt{1-4t}(1-t)}=\frac{({\rm log}(4/3))^{r-1}}{(r-1)!}\frac{\pi}{3 \sqrt{3}}\cdot
 \end{aligned}
 \end{equation}

\medskip
The last assertion follows immediately from the integral formulas (\ref{EU28}) and (\ref{EU33}).
\vspace{2em}
\subsection{Recurrence relations for tails of  multiple Ap{\'e}ry-like sums }\label{s2s5}
\myindent

\medskip
\begin{prop}\label{P4} Let ${\bf a}=(a_1,\ldots,a_r)$ be a non-empty composition. For each integer $n \geqslant 0$, we have
	\begin{equation}\label{EF28}
	\sigma({\bf a})_{n-1}-\sigma({\bf a})_{n}={n^{-a_r}}\sigma({\bf a}^{\rm init})_{n}
	\end{equation}
\end{prop}
Indeed the left-hand side is equal to
\begin{equation}
\begin{aligned}
&
\sum_{n_1>\ldots >n_r>n-1}\left({2n_1 \atop n_1}\right)^{-1}n_1^{-a_1}\ldots n_r^{-a_r}-
\sum_{n_1>\ldots >n_r>n}\left({2n_1 \atop n_1}\right)^{-1}n_1^{-a_1}\ldots n_r^{-a_r}
\\  
&\qquad=\left\{ {{n^{-a_r}}\sum_{n_1>\ldots >n_{r-1}>n}\left({2n_1 \atop n_1}\right)^{-1}n_1^{-a_1}\ldots n_{r-1}^{-a_{r-1}}\quad {\rm if}\quad r \geqslant 2,\atop n^{-a_1}\left(2n \atop n\right)^{-1}\qquad \quad\qquad\;\qquad\qquad\qquad\qquad\quad{\rm if }\quad r=1.}\right.
\end{aligned}
\end{equation}

\medskip

{\it Remark.-} The recurrence relation (\ref{EF28}), together with the bounds obtained in Section~\ref{s2s2}, yield an efficient algorithm to compute a multiple Apery-like sum $\sigma(a_1,\ldots,a_r)$ by computing simultaneously the values of $\sigma(a_1,\ldots,a_i)$ for $1 \leqslant i \leqslant r$.

\vspace{2em}

 \subsection{Values of multiple Ap{\'e}ry-like functions at all integer entries}
 \medskip
 
 \myindent In this section we show that the values of Ap{\'e}ry-like functions at all integer entries can be expressed in terms of Ap{\' e}ry-like sums (associated to compositions). More precisely:
 
 \medskip
 \begin{thm}
 	For any ${\bf a} \in {\bf Z}^{r}$, there exists a unique finite formal linear combination $\sum_{{\bf b}\in \mathscr{C}} f_{\bf b}{\bf b}$, where the $f_b$ are polynomials in ${\bf Q}[T]$, such that, for each $n \in {\bf N}$, we have
 	\begin{equation}
 	\sigma({\bf a})_{n}=\sum_{{\bf b}\in \mathscr{C}} f_{\bf b}(n)\sigma({\bf b})_n.
 	\end{equation}
 \end{thm}
The unicity immediately follows from the linear independence result 
in Theorem \ref{T2} of Section \ref{s2s7}.

We shall prove the existence by induction on the depth of $\bf a$. When $\bf a$ itself is a composition, we can just take $f_{\bf b}$  equal to $1$ when ${\bf b}={\bf a}$ and to $0$ otherwise. So we assume that ${\bf a}=(a_1,\ldots ,a_r)$ has at least one non-positive entry. We distinguish three cases :

{\it Case} 1 : {\it We have $r=1$ and ${\bf a}=(-c)$ with $c\geqslant 0$}.

 We work in this case by induction on $c$. We can write for each integer $n_1 > 0$
\begin{equation}\label{EM47}
\begin{array}{lll}
\left( {2n_1-2 \atop n_1-1}\right)^{-1}(n_1-1)^{c}-\left({2n_1 \atop n_1}\right)^{-1}n_1^{c} &=& \left({2n_1 \atop n_1}\right)^{-1} \left(\left(4-\frac{2}{n_1}\right)(n_1-1)^c-n_1^c\right)\\ 
&=& \left({2n_1 \atop n_1}\right)^{-1}\left(3 n_1^c + \sum_{d=0}^{c-1} \lambda_{c,d}\,n_1^d+\frac{2 (-1)^{c+1}}{n_1}\right),
\end{array}
\end{equation}
where
$$\lambda_{c,d}=4 \left( {c \atop d} \right) (-1)^{c-d}-2 \left( {c \atop d+1} \right) (-1)^{c-d-1}= 2 (-1)^{c-d}\frac{c!}{(d+1)!(c-d)!}(c+d+2).$$
Summing up these expressions for $n_1 >n$ yields
\begin{equation}
n^c \sigma(\varnothing)_{n}=\left( {2n \atop n}\right)^{-1}n^c=3\sigma(-c)_{n}+\sum_{d=0}^{c-1}\lambda_{c,d}\sigma(-d)_{n}+ 2 (-1)^{c+1}\sigma(1)_n
\end{equation}
and we therefore conclude by the induction hypothesis on c.

{\it Case }2 : {\it We have $r \geqslant 2$ and $a_1 \geqslant 1$}. 

Since $\bf a$ is not a composition, there exists at least one index $i$ such that $ 1 < i \leqslant r$ and $a_i \leqslant 0$. Write $a_i=-c$. In the series expression (\ref{EM23}) for $\sigma({\bf a})_n$,  we can first sum over the variable $n_i$. We note that
\begin{eqnarray}
\nonumber \sum_{n_i=n_{i+1}+1}^{n_{i-1}-1}n_i^c&=&\frac{B_{c+1}(n_{i-1})-B_{c+1}(n_{i+1}+1)}{c+1}\quad \text{when } i \neq r,
\\ \nonumber
\sum_{n_i=n+1}^{n_{i-1}-1}n_i^c&=&\frac{B_{c+1}(n_{i-1})-B_{c+1}(n+1)}{c+1}\quad \text{when } i=r,
\end{eqnarray}
were  $B_{c+1}$ denotes the Bernoulli polynomial  of index $c+1$. Hence
$\sigma({\bf a})_n$ can be expressed in the form $\sum_{{\bf b} \in {\bf Z}^{r-1}}P_{\bf b}(n)\sigma({\bf b})_n$, where $(P_{\bf b})_{{\bf b}\in {\bf Z}^r}$ is a family of polynomials in ${\bf Q}[T]$ with finite support. We conclude by the induction hypothesis on the depth.

{\it Case} 3 : {\it We have $r \geqslant 2$ and $a_1=-c$ with $c \geqslant 0$}. 

We again work by induction on $c$. Multiplying  formula (\ref{EM47}) by $n_2^{a_2}\ldots n_r^{a_r}$ and summing over the tuples $(n_1,\ldots ,n_r)$ with $n_1>\ldots >n_r>n$, we get 
\begin{equation}
\sigma(a_{2}-c,a_3,\ldots ,a_r)_{n}=3\sigma({\bf a})_{n}+\sum_{d=0}^{c-1}\lambda_{c,d}\,\sigma(-d,a_2,\ldots ,a_r)_{n}+ 2 (-1)^{c+1}\sigma(1,a_2,\ldots ,a_r)_n
\end{equation}
We can apply to the left hand side the induction hypothesis concerning the depth, to the terms in the sum the induction hypothesis on $c$ and to the last term the case 2 already treated. The result then follows.

\medskip
{\it Example.}- Let ${\bf a}=(a_1,\ldots,a_r)$ be a composition. We have 
\begin{equation}
\sigma(0,a_1,\ldots,a_r)_n=\frac{2}{3}\sigma(1,a_1,\ldots,a_r)_n+\frac{1}{3}\sigma(a_1,\ldots,a_r)_n \quad \text{ for all $n \geqslant0$.}
\end{equation}
\vspace{2em}

\subsection{Other iterated integral expressions of multiple Ap{\'e}ry-like sums and multiple zeta values}\label{s2s7a}
\myindent For every  non empty composition ${\bf a}=(a_1,\ldots,a_r)$ (admissible or not), we define  a power series $\sigma_{\bf a}$ by
\begin{equation}\label{EU44}
\sigma_{\bf a}(z)=\sum_{n_1>\ldots>n_r>0}\left({2n_1 \atop n_1}\right)^{-1}n_1^{-a_1}\ldots n_r^{-a_r}z^{n_1}.
\end{equation}
 It converges for  $|z|<4$ (and even  normally   for $|z|\leqslant 4$ when $\bf a$ is admissible).  We shall also consider its $n$-tail, defined for $n\geqslant 0$ by
 \begin{equation}
 \sigma_{\bf a}(z)_{n}=\sum_{n_1>\ldots>n_r>n}\left({2n_1 \atop n_1}\right)^{-1}n_1^{-a_1}\ldots n_r^{-a_r}z^{n_1}.
 \end{equation}
 We extend these definitions to the empty composition by setting $\sigma_{\varnothing}(z)=1$, \mbox{$\sigma_{\varnothing}(z)_n=\left({2n \atop n}\right)^{-1}z^n$}. We have $\sigma({\bf a})=\sigma_{\bf a}(1)$ and \mbox{$\sigma({\bf a})_n=\sigma_{\bf a}(1)_n$} for all ${\bf a}\in \mathcal{C}$.
\medskip
\begin{thm}\label{TH6}
Let ${\bf a}=(a_1,\ldots,a_r)$ be a non empty composition and $\varepsilon_1\ldots\varepsilon_k$ denote its associated binary word ${\bf w}({\bf a})$. We  have
	\begin{equation}\label{EU60}
	\sigma_{\bf a}(4\,{\rm sin}^2y)=2^k({\rm tan}\,y)^{\varepsilon_1}{\rm It}\int_{0}^{y}(({\rm tan}\,t)^{\varepsilon_1+\varepsilon_2-1}dt, \ldots,({\rm tan}\,t)^{\varepsilon_{k-1}+\varepsilon_k-1}dt,dt),
	\end{equation}
	and more generally
	\begin{equation}\label{EU47}
	\sigma_{\bf a}(4\,{\rm sin}^2y)_{n}=2^k\left({2n \atop n}\right)^{-1}({\rm tan }\, y)^{\varepsilon_1}{\rm It}\int_{0}^{y}(({\rm tan}\, t)^{\varepsilon_1+\varepsilon_2-1}dt,\ldots, ({\rm tan}\,t)^{\varepsilon_{k-1}+\varepsilon_k-1}dt,(4\,{\rm sin}^2t)^ndt)
	\end{equation}
for $n\geqslant 0$ and $y \in ]-\frac{\pi}{2}, \frac{\pi}{2}[$. (These formulas remain valid  for $y \in [-\frac{\pi}{2}, \frac{\pi}{2}]$  when $\bf a$ is admissible.)
\end{thm}

The proof is by induction on the weight $k$ of $\bf a$. When $k=1$, we have ${\bf a}=(1)$ and we have to prove that the two expressions
\begin{equation}
\begin{aligned}
u_n(y)&=\sum_{n_1>n}\left({2n_1 \atop n_1}\right)^{-1}n_1^{-1}(4\,{\rm sin}^2y)^{n_1}\\
v_n(y)&=2\left({2n \atop n}\right)^{-1} {\rm tan }\; y\int_{0}^{y}(4\,{\rm sin}^2t)^n dt
\end{aligned}
\end{equation}
are equal for $n\geqslant 0$ and $y \in ]-\frac{\pi}{2}, \frac{\pi}{2}[$. Both of them tend to $0$ when $n$ tends to $+\infty$. Hence our assertion follows from the fact that
\begin{equation}
\begin{aligned}
v_{n-1}(y)-v_{n}(y)&=2\left({2n \atop n}\right)^{-1} {\rm tan }\; y\int_{0}^{y}\left(\big(4-\frac{2}{n}\big)(4\,{\rm sin}^2t)^{n-1}-(4\,{\rm sin}^2t)^n\right)dt\\
&=\left({2n \atop n}\right)^{-1} {\rm tan }\; y\int_{0}^{y}\frac{1}{n}\,\frac{d}{dt}\left(\frac{(4\,{\rm sin}^2t)^n}{{\rm tan}\,t}\right)dt\\
&=\left({2n \atop n}\right)^{-1}n^{-1}(4\,{\rm sin}^2y)^n\,=\,u_{n-1}(y)-u_n(y).
\end{aligned}
\end{equation}

We now consider the case where $k \geqslant 2$ and $a_1=1$. By multiplying the equality \mbox{$u_{n_2}(y)=v_{n_2}(y)$} by $n_2^{-a_2}\ldots n_r^{-a_r}$ and summing over the tuples $(n_2,\ldots,n_r)$ such that $n_2>\ldots>n_r>n$, we get 
\begin{equation}
\sigma_{(1,a_2,\ldots,a_r)}(4\,{\rm sin}^2y)_n=2\,{\rm tan}\,y\int_{0}^{y}\sigma_{(a_2,\ldots,a_r)}(4\,{\rm sin}^2t)_{n}dt.
\end{equation}
Formula (\ref{EU47}) then  follows for ${\bf a}=(1,a_2,\ldots,a_r)$ from the induction hypothesis.

We finally consider the case where $a_1 \geqslant 2$. We note that
\begin{equation}
\frac{(4\,{\rm sin}^2y)^{n_1}}{n_1}=2\int_{0}^{y}\frac{(4\,{\rm sin}^2t)^{n_1}}{{\rm tan}\,t}dt
\end{equation}
for $n_1\geqslant 1$ and  $y \in ]-\frac{\pi}{2}, \frac{\pi}{2}[$. By multiplying this identity by $\left({2n_1 \atop n_1}\right)^{-1}n_1^{1-a_1}n_2^{-a_2}\ldots n_r^{-a_r}$  and summing over the tuples $(n_1,\ldots,n_r)$ such that $n_1>\ldots>n_r>n$, we get
\begin{equation}
\sigma_{(a_1,\ldots,a_r)}(4\,{\rm sin}^2y)_{n}=2\int_{0}^{y}\frac{1}{{\rm tan}\, y}\sigma_{a_1-1,a_2,\ldots,a_r}(4\,{\rm sin}^2t)_n dt
\end{equation}
for $y \in ]-\frac{\pi}{2}, \frac{\pi}{2}[$. Formula (\ref{EU47}) then
 follows from the induction hypothesis. It remains true for $y \in [-\frac{\pi}{2},\frac{\pi}{2}]$ by continuity.
\medskip

{\it Examples.}  1) When ${\bf a}=(\underbrace{2,\ldots,2}_r)$ with $r\geqslant 1$, we get from  Theorem \ref{TH6}
\begin{equation}\label{EU69}
\sigma_{(\underbrace{\scriptstyle 2,\ldots,2}_r)}(4\,{\rm sin}^2y)=2^{2r}{\rm It}\int_{0}^{y}(\underbrace{dt,\ldots,dt}_{2r})=\frac{2^{2r}y^{2r}}{2r!}
\end{equation}
for $y \in [-\frac{\pi}{2},\frac{\pi}{2}]$. In particular, by taking $y=\frac{\pi}{6}$, we get
\begin{equation}\label{EU71a}
\sigma(\underbrace{ 2,\ldots,2}_r)=\frac{\pi^{2r}}{3^{2r}(2r)!}\cdot
\end{equation}

Formula (\ref{EU69}) allows us to recover the Taylor expansion at $0$ of even powers of the arcsine function, stated by J. Borwein and M. Chamberland in \cite{BC}, but in fact already found in 1896 in a slightly different, but equivalent form, by F. G.  Teixeira in \cite{Tei}: for every integer $r\geqslant 1$ and every $x \in [-1,1],$ we have
\begin{equation}\label{EU166a}
({\rm arcsin}\;x)^{2r}=\frac{(2r)!}{2^{2r}}\sigma_{(\underbrace{\scriptstyle 2,\ldots,2}_r)}(4x^2)=\frac{(2r)!}{2^{2r}}\sum_{n_1>\ldots >n_r>0}\frac{(2x)^{2n_1}}{\left({2n_1\atop n_1}\right)n_1^2\ldots n_r^2}\cdot
\end{equation}

2) When ${\bf a}=(1,\underbrace{2,\ldots,2}_{r-1})$ with $r\geqslant 1$, we get from Theorem \ref{TH6}
\begin{equation}
\sigma_{\scriptsize(1,\underbrace{2,\ldots,2}_{r-1})}(4\,{\rm sin}^2y)=2^{2r-1}{\rm tan }\,y\; {\rm It}\int_{0}^{y}(\underbrace{dt,\ldots,dt}_{2r-1})=\frac{2^{2r-1}y^{2r-1}}{(2r-1)!}\,{\rm tan }\,y
\end{equation}
for  $y \in ]-\frac{\pi}{2},\frac{\pi}{2}[$. In particular, by taking $y=\frac{\pi}{6}$, we get
\begin{equation}
\sigma(1,\underbrace{2,\ldots,2}_{r-1})=\frac{\pi^{2r-1}}{3^{2r-1}(2r-1)!\sqrt{3}}\cdot
\end{equation}
\medskip

\begin{customcor}{1}
	Let $\bf a$ be a non empty composition and $\varepsilon_1\ldots\varepsilon_k$ denote its associated binary word. We have
	\begin{equation}
	\sigma({\bf a})=2^k(\sqrt{3})^{-\varepsilon_1}{\rm It}\int_{0}^{\pi/6}(({\rm tan}\,t)^{\varepsilon_1+\varepsilon_2-1}dt, \ldots,({\rm tan}\,t)^{\varepsilon_{k-1}+\varepsilon_k-1}dt,dt)
	\end{equation}
	and more generally
	\begin{equation}
	\sigma({\bf a})_n=2^k\left({2n \atop n}\right)^{-1}(\sqrt{3})^{-\varepsilon_1}{\rm It}\int_{0}^{\pi/6}(({\rm tan}\,t)^{\varepsilon_1+\varepsilon_2-1}dt, \ldots,({\rm tan}\,t)^{\varepsilon_{k-1}+\varepsilon_k-1}dt,(4\,{\rm sin}^2t)^n dt).
	\end{equation}
	for $n\geqslant 0$.
\end{customcor}
\medskip

\begin{customcor}{2}\label{cor2}
Let ${\bf a}=(a_1\ldots a_r)$ be a non empty admissible composition and $\varepsilon_1\ldots\varepsilon_k$ denote its associated binary word. We have
\begin{equation}\label{EU60a}
\zeta({\bf a})=\frac{2^{k+1}}{\pi} {\rm It}\int_{0}^{\pi/2} (dt,({\rm tan}\, t)^{\varepsilon_1+\varepsilon_2-1}dt,\ldots,({\rm tan} \, t)^{\varepsilon_{k-1}+\varepsilon_k-1}, dt).
\end{equation}
More generally, the double tails of $\zeta({\bf a})$ {\rm (as defined in Section \ref{s1s5})} are given by
\begin{equation}\label{EU61a}
\begin{aligned}
\zeta({\bf a})_{m,n}& =\frac{2^{k+1}}{\pi}\left({2m \atop m}\right)^{-1}\left({2n \atop n}\right)^{-1}\\&\times  {\rm It}\int_{0}^{\pi/2} ((4\,{\rm cos}^2t)^m dt,({\rm tan}\, t)^{\varepsilon_1+\varepsilon_2-1}dt,\ldots,({\rm tan} \, t)^{\varepsilon_{k-1}+\varepsilon_k-1},({4\,{\rm sin}^2t})^n dt).
\end{aligned}
\end{equation}
\end{customcor}

For all integers $m\geqslant 0$ and $n \geqslant 0$, we have, by the change of variables $x={\rm sin}^2t$,
\begin{equation}
\begin{aligned}
&\frac{2}{\pi}\int_{0}^{\pi}(4\, {\rm cos}^2t)^m(4\,{\rm sin}^2t)^ndt=\frac{1}{\pi}\int_{0}^{1}4^{m+n}(1-x)^mx^n\frac{dx}{\sqrt{x(1-x)}}\\&
4^{m+n}\frac{\Gamma(m+\frac{1}{2})\Gamma(n+\frac{1}{2})}{\pi \Gamma(m+n+1)}=\frac{(2m)!(2n)!}{m!n!(m+n)!}=\left({2m \atop m}\right)\left({2n \atop n}\right)\left({m+n \atop m}\right)^{-1}\cdot
\end{aligned}
\end{equation}
By applying this equality to a pair of integers $(m,n_1)$, then multiplying the result by $\left({2n_1 \atop n_1}\right)^{-1}n_1^{-a_1}\ldots n_r^{-a_r}$ and summing over the tuples $(n_1,\ldots,n_r)$ such that \mbox{$n_1 >\ldots>n_r>n$}, we get
\begin{equation}
\frac{2}{\pi}\int_{0}^{\pi/2}(4{\rm cos}^2t)^m\sigma_{\bf a}(4\,{\rm sin}^2t)_ndt=\left({2m \atop m}\right)\zeta({\bf a})_{m,n}.
\end{equation}
Formula (\ref{EU61a}) follows, by replacing $\sigma_{\bf a}(4{\rm sin}^2t)_n $ by its expression (\ref{EU47}).
\vspace{2em}
\subsection{On some definite integrals involving the tangent and cotangent \mbox{functions} }
 \myindent  For each integer $\ell \geqslant 1$, the polylogarithm ${\rm Li}_{\ell}$ of index $\ell$ is defined on the open unit disc $\bf D$ of $\bf C$ by the convergent power series
 \begin{equation}
 {\rm Li}_{\ell}(z)=\sum_{n=1}^{\infty}\frac{z^n}{n^\ell}.
 \end{equation} 
 When $\ell \geqslant 2$, this power series  converges normally on $\overline{\bf D}$. When $\ell=1$, we have ${\rm Li}_1(z)=-{\rm log}(1-z),$ where $\rm log$ is the principal determination of  the logarithm on $\bf D$,  ${\rm Li}_1$ extends by continuity to $\overline{\bf D}-\{ 1\}$, and $\displaystyle{\lim_{\substack{y \to 0\\y\neq 0}}\,y{\rm Li}_1(e^{2iy})=0}$.
 
 \medskip
 \begin{lem}\label{LU1}
 	For $p\geqslant0$ and $q\geqslant 0$, we have
 	\begin{equation}\label{EU75}
 	\begin{aligned}
 	\int_0^{y}\frac{(y-t)^p}{p!}{\rm tan}\, t \, \frac{t^q}{q!}dt&=(-1)^p\sum_{\ell=p}^{p+q}(-1)^{\ell}\left({\ell \atop p}\right)\frac{y^{p+q-\ell}}{(p+q-\ell)!}\frac{{\rm Li}_{\ell+1}(-e^{2iy})}{(2i)^\ell}\\
 	&-(-1)^q\sum_{\ell=q}^{p+q}\left({\ell \atop q}\right)\frac{y^{p+q-\ell}}{(p+q-\ell)!}\frac{{\rm Li}_{\ell+1}(-1)}{(2i)^\ell}
 	+i\frac{y^{p+q+1}}{(p+q+1)!}
 	\end{aligned}
 	\end{equation}
 	for $y \in ]\frac{\pi}{2},\frac{\pi}{2}[$. For $p\geqslant 0$ and $q\geqslant 1$, we have
 	\begin{equation}\label{EU76}
 	\begin{aligned}
 	\int_0^{y}\frac{(y-t)^p}{p!} \frac{1}{{\rm tan}\, t \,} \frac{t^q}{q!}dt&=-(-1)^p\sum_{\ell=p}^{p+q}(-1)^{\ell}\left({\ell \atop p}\right)\frac{y^{p+q-\ell}}{(p+q-\ell)!}\frac{{\rm Li}_{\ell+1}(e^{2iy})}{(2i)^\ell}\\
 	&+(-1)^q\sum_{\ell=q}^{p+q}\left({\ell \atop q}\right)\frac{y^{p+q-\ell}}{(p+q-\ell)!}\frac{{\rm Li}_{\ell+1}(1)}{(2i)^\ell}
 	-i\frac{y^{p+q+1}}{(p+q+1)!}\cdot
 	\end{aligned}
 	\end{equation}
 	for $y \in\; ]-{\pi},{\pi}[\,-\{0\}$.
 \end{lem}
 
 We shall prove the two equalities (\ref{EU75}) and (\ref{EU76}) at the same time. To do so, we introduce a sign, denoted by $\pm$,  and we denote by $\mp$ the opposite sign. We then have to show that
 \begin{equation}
 A_{p,q}^{\pm}(y)=\int_0^{y}\frac{(y-t)^p}{p!} ({{\rm tan}\, t})^{\pm 1}\,\frac{t^q}{q!}dt
 \end{equation}
 is equal to
 \begin{equation}
 \begin{aligned}
 B_{p,q}^{\pm}(y)=&\pm(-1)^p\sum_{\ell=p}^{p+q}(-1)^{\ell}\left({\ell \atop p}\right)\frac{y^{p+q-\ell}}{(p+q-\ell)!}\frac{{\rm Li}_{\ell+1}(\mp e^{2iy})}{(2i)^\ell}\\
 &\mp(-1)^q\sum_{\ell=q}^{p+q}\left({\ell \atop q}\right)\frac{y^{p+q-\ell}}{(p+q-\ell)!}\frac{{\rm Li}_{\ell+1}(\mp 1)}{(2i)^\ell}
 \pm i\frac{y^{p+q+1}}{(p+q+1)!}
 \end{aligned}
 \end{equation}
 for $y \in ]-\frac{\pi}{2}, \frac{\pi}{2}[$  when $\pm =+$, and for $y \in ]-\pi,\pi[-\{0\}$ when $\pm=-$.
 
 We shall argue by induction on $p$. We first note that both  functions $A_{p,q}^{\pm}$ and $B_{p,q}^{\pm}$ are differentiable on  the previous sets and tend to $0$ when $y$ tends to $0$. It is therefore sufficient to prove that their derivatives are equal at any point of these sets. We now compute these derivatives at such a point $y$. We have
 \begin{equation}
 \frac{d}{dy}A_{p,q}^{\pm}(y)=\begin{cases}
 A_{p-1,q}^{\pm}(y) & \text{ when } p \geqslant 1,\\
 ({\rm tan}\,y)^{\pm 1}\,\frac{y^q}{q!}& \text{ when } p =0.
 \end{cases}
 \end{equation}
 When $p=0$ we have
 \begin{equation}
 \begin{aligned}
 \frac{d}{dy}B_{0,q}^{\pm}(y)=&\pm \sum_{\ell=0}^{q-1}(-1)^\ell \frac{y^{q-1-\ell}}{(q-1-\ell)!}\frac{{\rm Li}_{\ell+1}(\mp e^{2iy})}{(2i)^\ell}\\&
 \pm \sum_{\ell=1}^{q}(-1)^\ell \frac{y^{q-\ell}}{(q-\ell)!}\frac{{\rm Li}_{\ell}(\mp e^{2iy})}{(2i)^{\ell-1}}- \frac{y^q}{q!}\frac{ 2ie^{2iy}}{1\pm e^{2iy}}\pm \frac{iy^q}{q!}.
 \end{aligned}
 \end{equation}
 The two first sums are opposite of each other, and we therefore have
 \begin{equation}
 \begin{aligned}
 \frac{d}{dy}B_{0,q}^{\pm}(y)&=\frac{-iy^q}{q!}\left(\frac{2e^{2iy}}{1\pm e^{2iy}}\mp 1\right)=\frac{-iy^q}{q!}\left(\frac{\mp 1+e^{2iy}}{1\pm e^{2iy}}\right)=({\rm tan}\,y)^{\pm 1}\,\frac{y^q}{q!}\cdot
 \end{aligned}
 \end{equation}
 When $p \geqslant 1$, we have
 \begin{equation}
 \begin{aligned}
 \frac{d}{dy}B_{p,q}^{\pm}(y)=&\pm(-1)^p\sum_{\ell=p}^{p+q-1}(-1)^{\ell}\left({\ell \atop p}\right)\frac{y^{p+q-1-\ell}}{(p+q-1-\ell)!}\frac{{\rm Li}_{\ell+1}(\mp e^{2iy})}{(2i)^\ell}\\
 &\pm(-1)^p\sum_{\ell=p}^{p+q}(-1)^{\ell}\left({\ell \atop p}\right)\frac{y^{p+q-\ell}}{(p+q-\ell)!}\frac{{\rm Li}_{\ell}(\mp e^{2iy})}{(2i)^{\ell-1}}\\
 &\mp(-1)^q\sum_{\ell=q}^{p+q-1}\left({\ell \atop q}\right)\frac{y^{p+q-1-\ell}}{(p+q-1-\ell)!}\frac{{\rm Li}_{\ell+1}(\mp 1)}{(2i)^\ell}\pm i\frac{y^{p+q}}{(p+q)!}\cdot
 \end{aligned}
 \end{equation}
 
 The sum of the first two terms of the right-hand side can be written as
 \begin{equation}
 \begin{aligned}
 &\pm(-1)^p\sum_{\ell=p-1}^{p+q-1}(-1)^{\ell}\frac{\ell!}{p!\,(\ell+1-p)!}((\ell+1-p)-(\ell+1))\frac{y^{p+q-1-\ell}}{(p+q-1-\ell)!}\frac{{\rm Li}_{\ell+1}(\mp e^{2iy})}{(2i)^\ell}\\
 =&\pm(-1)^{p-1}\sum_{\ell=p-1}^{p+q-1}(-1)^{\ell}\left({\ell \atop p-1}\right)\frac{y^{p+q-1-\ell}}{(p+q-1-\ell)!}\frac{{\rm Li}_{\ell+1}(\mp e^{2iy})}{(2i)^\ell},\\
 \end{aligned}
 \end{equation}
 hence we have
 \begin{equation}
 \frac{d}{dy}B_{p,q}^{\pm}(y)=B_{p-1,q}^{\pm}(y).
 \end{equation}
 By the induction hypothesis, we have $\frac{d}{dy}A_{p,q}^{\pm}(y)=\frac{d}{dy}B_{p,q}^{\pm}(y)$. This finishes the proof of Lemma \ref{LU1}.
 \medskip
 
 {\it Remark.}- Since the left-hand sides of (\ref{EU75}) and (\ref{EU76}) are real, we can replace all terms in the right-hand sides by their real parts.

\vspace{2em}
\subsection{Evaluation of $\sigma(2,\ldots,2,1,2,\ldots,2)$ and $\sigma(2,\ldots,2,3,2,\ldots,2)$}\label{s2s8}
\myindent 
In this section, we give explicit expressions of the Ap{\'e}ry-like sums $\sigma({\bf a})$, when $\bf a$ is an admissible composition of the form $(2,\ldots,2,1,2,\ldots,2)$ or $(2,\ldots,2,3,2,\ldots,2)$, as  \mbox{$\bf Q$-linear} combinations of the real numbers of the form $\pi^{k-2\ell-1}\zeta(2\ell+1)$ and $\pi^{k-2\ell}\sqrt{3}\,{\rm L}(2\ell,\chi)$, where $k$ is the weight of $\bf a$, $1 \leqslant \ell \leqslant \frac{k-1}{2}$ and $\chi$ is the unique non principal Dirichlet character modulo $3$.
\medskip

\begin{thm}\label{TH7}
	For all integers $a\geqslant 1$, $b\geqslant 0$, we have
	\begin{equation}
	\begin{aligned}
	\sigma(\underbrace{2,\ldots,2}_a,1,\underbrace{2,\ldots,2}_b)&=\sum_{r=a}^{a+b}\left({2r \atop 2a-1}\right)\frac{\left(\pi/3\right)^{k-1-2r}}{(k-1-2r)!}(-1)^{r}(1-3^{-2r})\zeta(2r+1)
	\\ & +\sum_{r=a}^{a+b}\left({2r-1 \atop 2a-1}\right)\frac{\left(\pi/3\right)^{k-2r}}{(k-2r)!}(-1)^{r}\sqrt{3}\,{\rm L}(2r,\chi)
	\\ &-\sum_{r=b+1}^{a+b}\left({2r \atop 2b+1}\right)\frac{\left(\pi/3\right)^{k-1-2r}}{(k-1-2r)!}(-1)^{r}2(1-2^{-2r})\zeta(2r+1),
	\end{aligned}
	\end{equation}
	where $k=2a+2b+1$.
\end{thm}

\medskip

\begin{thm}\label{TH8}
	For all integers $a\geqslant 0$, $b\geqslant 0$, we have
	\begin{eqnarray}
	\nonumber\sigma(\underbrace{2,\ldots,2}_a,3,\underbrace{2,\ldots,2}_{b})&=&-\sum_{r={\rm max}(a,1)}^{a+b+1}\left({2r \atop 2a}\right)\frac{\left(\pi/3\right)^{k-1-2r}}{(k-1-2r)!}(-1)^{r}(1-2^{-2r})(1-3^{-2r})\zeta(2r+1)
	\\ &-& \sum_{r=a+1}^{a+b+1}\left({2r-1 \atop 2a}\right)\frac{\left(\pi/3\right)^{k-2r}}{(k-2r)!}(-1)^{r}\sqrt{3}(1+2^{1-2r}){\rm L}(2r,\chi)
	\\\nonumber &+&\sum_{r=b+1}^{a+b+1}\left({2r \atop 2b+2}\right)\frac{\left(\pi/3\right)^{k-1-2r}}{(k-1-2r)!}2(-1)^{r}\zeta(2r+1),
	\end{eqnarray}
	where $k=2a+2b+3$.
\end{thm}

 Theorem \ref{TH6} yields
\begin{equation}\label{EU77}
\begin{aligned}\sigma_{(\underbrace{\scriptstyle 2,\ldots,2}_a,1 ,\underbrace{ \scriptstyle 2,\ldots,2}_b)}(4\, {\rm sin}^2y)&=2^{2a+2b+1}\;{\rm It}\int_{0}^{y}(\underbrace{dt,\ldots,dt}_{2a-1},{\rm tan}\, t\, dt,\underbrace{dt,\ldots,dt}_{2b+1})\\
&=2^{2a+2b+1}\int_{0}^{y}\frac{(y-t)^{2a-1}}{(2a-1)!}{\rm tan}\, t\, \frac{t^{2b+1}}{(2b+1)!}dt
\end{aligned}
\end{equation}
for $a\geqslant 1$, $b \geqslant 0$ and $y \in [-\frac{\pi}{2},\frac{\pi}{2}]$, and similarly that
\begin{equation}\label{EU78}
\begin{aligned}\sigma_{(\underbrace{\scriptstyle 2,\ldots,2}_a,3 ,\underbrace{ \scriptstyle 2,\ldots,2}_{b})}(4\, {\rm sin}^2y)&=2^{2a+2b+3}\;{\rm It}\int_{0}^{y}(\underbrace{dt,\ldots,dt}_{2a},({\rm tan}\, t)^{-1}\, dt,\underbrace{dt,\ldots,dt}_{2b+2})\\
&=2^{2a+2b+3}\int_{0}^{y}\frac{(y-t)^{2a}}{(2a)!}({\rm tan}\, t)^{-1}\, \frac{t^{2b+2}}{(2b+2)!}dt
\end{aligned}
\end{equation}
for $a\geqslant 0$, $b \geqslant 0$ and $y \in [-\frac{\pi}{2},\frac{\pi}{2}]$.
\medskip

We now take $y=\frac{\pi}{6}$ in  formulas (\ref{EU77}) and (\ref{EU78}) and apply Lemma \ref{LU1} to compute the integrals on the right-hand sides. Since the left-hand sides are real numbers, we have, for $a\geqslant 1$, $b \geqslant 0$ and  $k=2a+2b+1$,
\begin{equation}
\begin{aligned}
\sigma(\underbrace{2,\ldots,2}_a,1,\underbrace{2,\ldots,2}_b)&=-2\sum_{r=a}^{a+b}(-1)^{r}\left({2r \atop 2a-1}\right)\frac{(\frac{\pi}{3})^{k-1-2r}}{(k-1-2r)!}{\rm Re}({\rm Li}_{2r+1}(-e^{i\pi/3}))\\
&-2\sum_{r=a}^{a+b}(-1)^{r}\left({2r-1 \atop 2a-1}\right)\frac{(\frac{\pi}{3})^{k-2r}}{(k-2r)!}{\rm Im}({\rm Li}_{2r}(-e^{i\pi/3}))\\
&+2\sum_{r=b+1}^{a+b}(-1)^{r}\left({2r \atop 2b+1}\right)\frac{(\frac{\pi}{3})^{k-1-2r}}{(k-1-2r)!}{\rm Li}_{2r+1}(-1),
\end{aligned}
\end{equation}
and for $a\geqslant 0$, $b \geqslant 0$ and $k=2a+2b+3$,
\begin{equation}
\begin{aligned}
\sigma(\underbrace{2,\ldots,2}_a,3,\underbrace{2,\ldots,2}_{b})&=-2\sum_{r=a}^{a+b+1}(-1)^{r}\left({2r \atop 2a}\right)\frac{(\frac{\pi}{3})^{k-1-2r}}{(k-1-2r)!}{\rm Re}({\rm Li}_{2r+1}(e^{i\pi/3}))\\
&-2\sum_{r=a+1}^{a+b+1}(-1)^{r}\left({2r-1 \atop 2a}\right)\frac{(\frac{\pi}{3})^{k-2r}}{(k-2r)!}{\rm Im}({\rm Li}_{2r}(e^{i\pi/3}))\\
&+2\sum_{r=b+1}^{a+b+1}(-1)^{r}\left({2r \atop 2b+2}\right)\frac{(\frac{\pi}{3})^{k-1-2r}}{(k-1-2r)!}{\rm Li}_{2r+1}(1).
\end{aligned}
\end{equation}

Theorems \ref{TH7} and \ref{TH8} follow from these formulas and from the following lemma.
\medskip

\begin{lem}\label{LU2}
	We have ${\rm Re}({\rm Li}_1(e^{i\pi/3}))=0$ and, for $r \geqslant 1$,
	\begin{equation}
	\begin{aligned}
		2\,{\rm Re}({\rm Li}_{2r+1}(-e^{i\pi/3}))&=-(1-3^{-2r})\,\zeta(2r+1),\\
		2\,{\rm Im}({\rm Li}_{2r}(-e^{i\pi/3}))&=-\sqrt{3}\,{\rm L}(2r, \chi),\\
		2\,{\rm Re}({\rm Li}_{2r+1}(e^{i\pi/3}))&=(1-2^{-2r})(1-3^{-2r})\,\zeta(2r+1),\\
		2\,{\rm Im}({\rm Li}_{2r}(e^{i\pi/3}))&=(1+2^{1-2r})\sqrt{3}\,{\rm L}(2r, \chi),\\
		{\rm Li}_{2r+1}(-1)&=-(1-2^{-2r})\,\zeta(2r+1),\\
		{\rm Li}_{2r+1}(1)&=\;\zeta(2r+1).	
	\end{aligned}
	\end{equation}
\end{lem}

We have ${\rm Li}_1(z)=-{\rm log}(1-z)$ for $z \in \overline{\bf D}-\{1\}$, hence ${\rm Re}({\rm Li}_1(z))=-{\rm log}\, |1-z|$ and ${\rm Re}({\rm Li}_1(e^{i\pi/3}))=-{\rm log}\, |1-e^{i\pi/3}|=0$. From now on, we assume $r \geqslant 1$. We have ${\rm Li}_{2r+1}(1)=\zeta(2r+1)$ and
\begin{equation}
\begin{aligned}
2 \,{\rm Re}\; {\rm Li}_{2r+1}(-e^{i\pi/3})&=\sum_{n=1}^{\infty}\frac{1}{n^{2r+1}}\times\begin{cases}
-1 & {\rm if}\quad n \not\equiv 0\,(3),\\
\;\;2 & {\rm if}\quad n \equiv 0\,(3)
\end{cases}\\
&=-\sum_{n=1}^{\infty}\frac{1}{n^{2r+1}}+3\sum_{n=1}^{\infty}\frac{1}{(3n)^{2r+1}}=-(1-3^{-2r})\zeta(2r+1),
\end{aligned}
\end{equation}
\begin{equation}
\begin{aligned}
2\;{\rm Im}\; {\rm Li}_{2r}(-e^{i\pi/3})&=\sum_{n=1}^{\infty}\frac{1}{n^{2r}}\times \begin{cases}
-\sqrt{3} & {\rm if}\quad n \equiv 1\,(3),\\
\,\sqrt{3} & {\rm if}\quad n \equiv 2\,(3),\\
\;0 & {\rm if}\quad n \equiv 3\,(3)
\end{cases}\\
&=-\sqrt{3}\,{\rm L}(2r,\chi).
\end{aligned}
\end{equation}
Since ${\rm Li}_k(z)+{\rm Li}_k(-z)=2^{1-k}{\rm Li}_k(z^2)$ and $\overline{{\rm Li}_k(z)}={\rm Li}_k(\overline{z})$ for $k \geqslant 2$ and $|z|=1$, we have
\begin{equation}
\begin{aligned}
{\rm Re}\; {\rm Li}_{2r+1}(e^{i\pi/3})&=-(1-2^{-2r}) \,{\rm Re}\; {\rm Li}_{2r+1}(-e^{i\pi/3}),\\
{\rm Im}\; {\rm Li}_{2r}(e^{i\pi/3})&=-(1+2^{1-2r}){\rm Im}\; {\rm Li}_{2r}(-e^{i\pi/3}),\\
{\rm Li}_{2r+1}(-1)&=-(1-2^{-2r}) \,{\rm Li}_{2r+1}(1).
\end{aligned}
\end{equation}
This completes the proof.
\vspace{2em}
\subsection{The $\bf Q$-subspace  generated by  the previous Ap{\'e}ry-like sums  }\label{s2s9}
\medskip
\begin{thm}\label{th7}
	Let $k$ be an odd positive integer. Let $ V$ be the $\bf Q$-subspace of $\bf R$ generated by the  $k-1$ Ap{\'e}ry-like sums $\sigma({\bf a})$, where $\bf a$ is an admissible composition of weight $k$ either of the form $({2,\ldots,2},3,2,\ldots,2)$ or of the form $({2,\ldots,2},1,2,\ldots,2)$. It is also generated by the $k-1$ real numbers  $\pi^{k-2r}\sqrt{3}\,L(2r,\chi) $ and $\pi^{k-2r-1}\zeta(2r+1)$, where $1\leqslant r\leqslant\frac{k-1}{2}$.
\end{thm}

\medskip

{\it Remark.}- It is reasonable to conjecture that $V$ is of dimension $k-1$ over $\bf Q$, $i.\,e.$ that each of the two
 sets of $k-1$ generators considered in Theorem \ref{th7} is free. At least for $k\leqslant 19$, we have found by running the PSLQ algorithm no $\bf Z$-linear relation between them, with coefficients of size at most $10^{80}$.

We can express Theorems \ref{TH7} and \ref{TH8}, for compositions of weight $k$, by a relation between block matrices
\begin{equation}
\left(\begin{array}{c}X\\\hdashline Y\end{array}\right)=\left(\begin{array}{c:c}
A&B\\ \hdashline C& D
\end{array}\right)\left(\begin{array}{c}U\\\hdashline V\end{array}\right)
\end{equation}
where $X$, $Y$, $U$, $V$ are the column vectors of size $\frac{k-1}{2}$ defined by
\begin{equation}
\begin{aligned}
X&=\Big(\zeta(\underbrace{2,\ldots,2}_{p-1},3,\underbrace{2,\ldots,2}_{\frac{k-1}{2}})\Big)_{1\leqslant p\leqslant \frac{k-1}{2}} \\Y&=\Big(\zeta(\underbrace{2,\ldots,2}_{\frac{k+1}{2}-p},1,\underbrace{2,\ldots,2}_{p-1})\Big)_{1\leqslant p\leqslant\frac{k-1}{2}}
\\U&=\Big((-1)^q\frac{(\pi/3)^{k-2q}}{(k-2q)!})\,2^{1-2q}\,\sqrt{3}\;{\rm L}(2q,\chi)\Big)_{1\leqslant q\leqslant\frac{k-1}{2}}
\\V&=\Big((-1)^q\frac{(\pi/3)^{k-2q-1}}{(k-2q-1)!})2^{-2q}3^{-3q}\zeta(2q+1)\Big)_{1\leqslant q\leqslant\frac{k-1}{2}}
\end{aligned}
\end{equation}
and the coefficients of the matrices $A$, $B$, $C$, $D$ are given, for $1\leqslant p \leqslant \frac{k-1}{2}$ and $1\leqslant q \leqslant \frac{k-1}{2}$, by
\begin{equation}
\begin{aligned}
a_{p,q}&=-\left({2q-1 \atop 2p-2}\right)(2^{2q-1}+1),
\\b_{p,q}&=-\left({2q \atop 2p-2}\right)(2^{2q}-1)(3^{2q}-1)+\left({2q \atop k+1-2p}\right)2^{2q+1}3^{2q},
\\c_{p,q}&=\left({2q-1 \atop k-2p}\right)\,2^{2q-1},
\\d_{p,q}&=-\left({2q \atop k-2p}\right)2^{2q}(3^{2q}-1)-2\left({2q \atop 2p-1}\right)(2^{2q}-1)3^{2q},
\end{aligned}
\end{equation}
where the binomial coefficient $\left({m \atop n}\right)$ is considered to be equal to $0$ whenever $m>n$.

To prove the theorem, it is sufficient to prove the following lemma:

\medskip

\begin{lem}
	The determinant of $\left(\begin{array}{c:c}
	A&B\\ \hdashline C& D
	\end{array}\right)$is different from $0$.
\end{lem}

To prove this lemma, we can replace the matrices $B$ and $D$ by the matrices $B'$ and $D'$ where
\begin{equation}
b_{p,q}'=b_{p,q}/2^{v_2(q)+2},\qquad d_{p,q}'=d_{p,q}/2^{v_2(q)+2},
\end{equation}
where $v_2$ is the $2$-adic valuation. The matrix \smash{$\left(\begin{array}{c:c}
	A&B'\\ \hdashline C& D'
	\end{array}\right)$} has a non zero determinant since:
\\- all its coefficients are  integers;
\\- all its diagonal coefficients  are odd;
\\- all its coefficients  below the diagonal are even.

Indeed, these assertions follow from the following observations:
\\- The coefficients $a_{p,q}$ are integers; the coefficient $a_{q,q}$ is equal to $-(2q-1)(2^{2q-1}+1)$, hence is odd; we have $a_{p,q}=0$ whenever $p>q$.
\\- The coefficients $b'_{p,q}$ are  integers since $v_2(3^{2q}-1)>v_2(q)+2$ and ${2q+1}>v_2(q)+2.$
\\- All coefficients $c_{p,q}$ are even integers.
\\- Since $\left({2q\atop 2p-1}\right)=\frac{2q}{2p-1}\left({2q-1 \atop 2p-2}\right)$, we have $v_2\left(2\Big({2q \atop 2p-1}\Big)\right)\geqslant v_2(q)+2$. On the other hand, we have  \mbox{$v_2(3^{2q}-1)>v_2(q)+2$}. Hence $d'_{p,q}$ is an integer, $d'_{q,q} $ is  odd, and  $d'_{p,q}$ is an even integer whenever $p>q$.

This proves the lemma.

\vspace{2em}

\subsection{A direct proof of Zagier's formula for $\zeta(2,\ldots,2,3,2,\ldots,2)$}\label{s2s8a}
\myindent The expressions that we get in  Theorems \ref{TH7} and \ref{TH8}  for the Ap{\'e}ry-like sums of the form $\sigma(2,\ldots,2,1,2,\ldots,2)$ and $\sigma(2,\ldots,2,3,2,\ldots,2)$ have striking similarities with those obtained by D. Zagier for $\zeta(2,\ldots,2,3,2,\ldots,2)$ in \cite{Zag}.  Zagier's method relies on a  delicate comparision of generating series. We present here a simpler way to get the same results  directly. It is based on the integral expressions obtained in Section \ref{s2s7a}. For the convenience of the reader, we restate Zagier's theorem in our notations, since his $\zeta(a_1,\ldots,a_r)$ is our $\zeta(a_r,\ldots,a_1)$.

\medskip

\begin{thm}\label{Th7}
	For all integers $a \geqslant 0$ and $b \geqslant 0$, we have\begin{equation}
	\begin{aligned}
	\zeta(\underbrace{2,\ldots,2}_{a},3,\underbrace{2,\ldots,2}_b)&= 2\sum_{r=b+1}^{a+b+1}(-1)^r\left({2r \atop 2b+2}\right)\frac{\pi^{k-1-2r}}{(k-2r)!}\zeta(2r+1)\\
	&-2\sum_{r=a+1}^{a+b+1}(-1)^r\left({2r \atop 2a+1}\right)\frac{\pi^{k-1-2r}}{(k-2r)!}(1-2^{-2r})\zeta(2r+1),		\end{aligned}
	\end{equation}
	where $k=2a+2b+3$.
\end{thm}

We indeed have, by Corollary \ref{cor2} of Theorem \ref{TH6}
\begin{equation}
\begin{aligned}
\zeta(\underbrace{2,\ldots,2}_{a},3,\underbrace{2,\ldots,2}_b)&= \frac{2^{k+1}}{\pi}{\rm It}\int_0^{\pi/2}(\underbrace{dt,\ldots,dt}_{2a+1},\frac{1}{{\rm tan}\, t}dt,\underbrace{dt,\ldots,dt}_{2b+2})\\
&= \frac{2^{k+1}}{\pi}\int_0^{\pi/2} \frac{(\frac{\pi}{2}-t)^{2a+1}}{(2a+1)!}\frac{1}{{\rm tan}\, t}\frac{t^{2b+2}}{(2b+2)!}dt.
\end{aligned}
\end{equation}
Therefore, by formula (\ref{EU76}) of Lemma \ref{LU1}, in which we replace all terms on the right-hand side by their real parts, we have
\begin{equation}
\begin{aligned}
\zeta(\underbrace{2,\ldots,2}_{a},3,\underbrace{2,\ldots,2}_b)&= 2\sum_{r=a+1}^{a+b+1}(-1)^r\left({2r \atop 2a+1}\right) \frac{\pi^{k-2r-1}}{(k-2r)!}{\rm Li}_{2r+1}(-1)\\
&+2\sum_{r=b+1}^{a+b+1}(-1)^r\left({2r \atop 2b+2}\right) \frac{\pi^{k-2r-1}}{(k-2r)!}\zeta(2r+1).
\end{aligned}
\end{equation}
The theorem follows, since ${\rm Li}_{2r+1}(-1)=-(1-2^{-2r})\zeta(2r+1)$ by Lemma \ref{LU2}.

\vspace{2em}

\subsection{Linear independence of tails of multiple Ap{\'e}ry-like  sums }\label{s2s7}
\myindent In this section, we state and prove a strong linear independence theorem for tails of multiple Ap{\' e}ry-like sums. It is the main tool needed to prove the unicity results in Section~\ref{s3}.

We recall that $\mathscr{C}$  denotes the set  of all compositions and  $\mathscr{C}^{\ast}$ the set $\mathscr{C}-\{\varnothing\}$.

\medskip
\begin{thm}\label{T2}
	Let $\left(f_{\bf a}\right)_{{\bf a}\in \mathscr{C}}$ be a family with finite support  indexed by $\mathscr{C}$ of rational ${\text{  functions}}$ in ${\bf C}(T)$. For each integer $n$ large enough, the numbers $f_{\bf a}(n)$ are then well defined. Assume that the real numbers
	\begin{equation}\label{E21}
	F(n)=\sum_{{\bf a} \in \mathscr{C}}f_{\bf a}(n)\sigma({\bf a})_n
	\end{equation}
	vanish for  $n$ large enough. Then $f_{\bf a}=0$ for all ${\bf a} \in \mathscr{C}$.
\end{thm} 

The proof is by double induction. Assume that not all $f_{\bf a}$ are  $0$, and then denote by $r$  the largest depth of the compositions ${\bf a}$ for which $f_{\bf a}\neq 0$. We have $r\geqslant1$. By multiplying the family $(f_{\bf a})_{{\bf a}\in \mathcal{C}}$ by a suitable non zero polynomial $P\in {\bf C}[T]$, we can reduce to the case where $f_{\bf a}$ is a polynomial for each $\bf a \in \mathcal{C}$ of depth $r$. The double induction will then be first on $r$, second  on the maximum degree of the polynomials $f_{\bf a}$ for ${\bf a}\in \mathscr{C}$ of depth~$r$.

By using Proposition \ref{P4}, we get, for $n$ large enough, 
\begin{eqnarray}\nonumber
F(n-1)-F(n)&=& \sum_{{\bf a} \in \mathscr{C}}f_{\bf a}(n-1)\sigma({\bf a})_{n-1}- \sum_{{\bf a} \in \mathscr{C}}f_{\bf a}(n)\sigma({\bf a})_{n}\\ \nonumber
&=& f_{\varnothing}(n-1)\left({2n-2 \atop n-1}\right)^{-1}-f_{\varnothing}(n)\left({2n \atop n}\right)^{-1}\\ \nonumber
&&+\sum_{{\bf a} \in \mathscr{C}^{\ast}} f_{\bf a}(n-1)\left(\sigma({\bf a})_{n}+n^{|{\bf a}^{\rm init}|-|{\bf a}|}\sigma({\bf a}^{\rm init})_{n}\right)-\sum_{{\bf a} \in \mathscr{C}^{\ast}}f_{\bf a}(n)\sigma({\bf a})_{n}
\\ \nonumber
&=&\left(f_{\varnothing}(n-1)\frac{4n-2}{n}-f_{\varnothing}(n)\right)\left({2n \atop n}\right)^{-1}\\ \nonumber
&&+\sum_{{\bf a} \in \mathscr{C}^{\ast}}\left(f_{\bf a}(n-1)-f_{\bf a}(n)\right) \sigma({\bf a})_{n}+\sum_{{\bf a} \in \mathscr{C}^{\ast}} f_{\bf a}(n-1)n^{|{\bf a}^{\rm init}|-|{\bf a}|}\sigma({\bf a}^{\rm init})_{n}.
\end{eqnarray}

Therefore we can write 
\begin{equation}
F(n-1)-F(n)=\sum_{{\bf a} \in \mathscr{C}}g_{\bf a}(n)\sigma({\bf a})_{n},
\end{equation}
where  $g_{\bf a}$ is the rational function defined by
\begin{eqnarray}\nonumber
g_{\bf a}(T)=\begin{cases}
f_{\varnothing}(T-1)\frac{4T-2}{T}-f_{\varnothing}(T)+\sum_{\substack{{\bf b} \in \mathscr{C},\\{\rm depth}(\bf{b})=1}}f_{\bf b}(T-1)T^{-|{\bf{b}}|}\quad {\rm when}\; {\bf a}=\varnothing,\\
f_{\bf a}(T-1)-f_{\bf a}(T)+ \sum_{\substack{{\bf b} \in \mathscr{C}\\{\bf b}^{\rm init}={\bf a}}}
f_{\bf b}(T-1)T^{|{\bf a}|-|{\bf b}|}\;\; \qquad\quad\;{\rm when } \;{\bf a} \neq \varnothing.
\end{cases}
\end{eqnarray}
We note that $g_{\bf a}=0$ for ${\bf a}$ of depth $>r$, and that, for ${\bf a}$ of depth $r$, $g_{\bf a}$ is the polynomial $f_{\bf a}(T-1)-f_{\bf a}(T)$, whose degree is smaller than the degree of $f_{\bf a}$ when $f_{\bf a} \neq 0$.

If $F(n)$ vanishes for $n$ large enough, the  same is true of $F(n-1)-F(n)$. By the induction hypotheses of our double induction, we have therefore
\begin{equation}\nonumber
g_{\bf a}=0 \quad \text{ for all } {\bf a} \in \mathscr{C}.
\end{equation}
The proof of the theorem will now be completed in three steps.
\\{\rm Step 1 : } {\it For each ${\bf a} \in \mathscr{C}^{\ast}$, $f_{\bf a}$ is a constant polynomial.}

Assume that this is false. Let then  $\bf a$ be  a composition of maximum depth for which $f_{\bf a}$ is not constant and $s$ be its depth. The relation $g_{\bf a}=0$ implies that
\begin{equation}\label{E24}
f_{\bf a}(T)-f_{\bf a}(T-1)=\sum_{\substack{{\bf b} \in \mathscr{C}\\{\bf b}^{\rm init}={\bf a}}}f_{\bf b}(T-1)T^{|{\bf a}|-|{\bf b}|}.
\end{equation}
In the right-hand side, all $f_{\bf b}$ occurring are constant rational functions since the relation ${\bf b}^{\rm init}={\bf a}$ implies depth(${\bf b}$)$>s$. The right-hand side is therefore a ${\bf C}$-linear combination of the rational functions $T^{-m}$, with $m\geqslant1$. These have no pole distinct from $0$.  
Let us assume that $f_{\bf a}$ has at least a pole in $\bf{C}$. Let $\alpha$ be a pole of $f_{\bf a}$ with largest real value and $\beta$ be a pole of $f_{\bf a}$ with smallest real value. Then $\alpha+1$ is a pole of $f_{\bf a}(T-1)$ and not of $f_{\bf a}(T)$, whereas $\beta$ is a pole of $f_{\bf a}(T)$ but not of $f_{\bf a}(T-1)$. Hence both $\alpha+1$ and $\beta$ are poles of $f_{\bf a}(T)-f_{\bf a}(T-1)$, and they are distinct.  This is incompatible with formula $(\ref{E24})$. Therefore $f_{\bf a}$ is a polynomial. Since $f_{\bf a}$ is not constant, $f_{\bf a}(T)-f_{\bf a}(T-1)$ is a non-zero polynomial. But this again is incompatible with formula $(\ref{E24})$. Hence step 1 follows.

\medskip
\noindent{\rm Step 2 :} {\it For each ${\bf b} \in \mathscr{C}$ of depth $\geqslant 2$, we have  $f_{\bf b}=0$}. 

For each ${\bf a} \in \mathscr{C}^{\ast}$, the relation $g_{\bf a}=0$ yields, by step 1,
\begin{equation}
\sum_{\substack{{\bf b} \in \mathscr{C}\\{\bf b}^{\rm init}={\bf a}}}f_{{\bf b}}(T-1)T^{|{\bf a}|-|{\bf b}|}=0.
\end{equation}
where all $f_{\bf b}$ occurring are constant rational functions. 
The rational functions \mbox{$T^{-m}$ ($m \geqslant 1$)} are linearly independent, and two non-empty compositions ${\bf b}$ and ${\bf c}$ such that ${\bf b}^{\rm init}={\bf a}$ and ${\bf c}^{\rm init}={\bf a}$ and $|{\bf b}|=|{\bf c}|$, are necessarily equal. Therefore  we have $f_{\bf b}=0$ for all ${\bf b} \in \mathscr{C}$ such that ${\bf b}^{\rm init}={\bf a}$. Since this is true for all ${\bf a}\in \mathscr{C}^{\ast}$, we have $f_{\bf b}=0$ whenever $\bf b$ has \mbox{depth $\geqslant 2$.}

\medskip
\noindent{\rm Step 3 : } {\it For each ${\bf b} \in \mathscr{C}$ of depth $\leqslant 1$, we have  $f_{\bf b}=0$.}

The relation $g_{\varnothing}=0$ can, due to step 2, be written as
\begin{equation}\label{E27}
f_{\varnothing}(T)-f_{\varnothing}(T-1)\frac{4T-2}{T}=\sum_{\substack{{\bf b} \in \mathscr{C},\\{\rm depth}({\bf b})=1}}f_{\bf b}(T-1)T^{-|{\bf b}|},
\end{equation}
where all $f_{\bf b}$ occuring are constant rational functions by step 1. 
Let us assume that $f_{\varnothing}$ has at least a pole in $\bf{C}$. Let $\alpha$ be a pole of $f_{\varnothing}$ with largest real value and $\beta$ be a pole of $f_{\varnothing}$ with smallest real value. If ${\rm Re}(\alpha) \geqslant 0$, then $\alpha+1$ is a pole of $f_{\varnothing}(T-1)$ and not of $f_{\varnothing}(T)$ and $\frac{4T-2}{T}$  does not vanish at $\alpha+1$, therefore $f_{\varnothing}(T)-f_{\varnothing}(T-1)\frac{4T-2}{T}$ has a pole at $\alpha+1 $. This is a contradiction since $\alpha+1 \neq 0$.  If ${\rm Re}(\alpha)<0$, then   ${\rm Re}(\beta) <0$ and $\beta$  is a pole of $f_{\varnothing}(T)$ but not of $f_{\varnothing}(T-1)$, and $\frac{4T-2}{T}$  has no pole at $\beta$, therefore $f_{\varnothing}(T)-f_{\varnothing}(T-1)\frac{4T-2}{T}$ has a pole at $\beta $. This is again a contradiction since $\beta \neq 0$. We thus have proved that $f_{\varnothing}$ has no pole in $\mathbb{C}$, hence is a polynomial. If this polynomial is not equal to $0$ and $aT^n$ denotes its  monomial of highest degree, then we have $f_{\varnothing}(t)-f_{\varnothing}(t-1)\frac{4t-2}{t}=-3at^n+o(t^n)$ when $t$ goes to $+\infty$. This  contradicts formula~$(\ref{E27})$. Hence $f_{\varnothing}=0$ and by the linear independence of the $T^{-m}$ for $m \geqslant1$, all $f_{\bf b}$  with depth($\bf b$)=$1$ are \mbox{equal to $0$}. This concludes the proof.

\vspace{2em}
\section{Expressing multiple zeta values in terms of  multiple   Ap{\'e}ry-like  sums }\label{s3}
\subsection{Definition of a graded map $\delta:{\bf Z}^{({\mathcal{B}})} \rightarrow {\bf Z}^{({\mathcal{A}})}$ }\label{s3as1}
\myindent Let $\mu $ denote the $\bf Z$-linear map from $
 {\bf Z}^{(\mathcal{A}^*)}$ to $ {\bf Z}^{(\mathcal{A})}$ such that $\mu({\bf a})={\bf a}^{\rm init}$  for each non-empty admissible composition  ${\bf a}$. For each integer $k \geqslant 1$, $\mu$ induces a $\bf Z$-linear map
\begin{equation}
\mu_k : {\bf Z}^{\mathcal{A}_k}\rightarrow \bigoplus_{0\leqslant k' <k}{\bf Z}^{\mathcal{A}_{k'}}.
\end{equation}
 We denote by $\mu_{k',k}:{\bf Z}^{\mathcal{A}_k}\rightarrow {\bf Z}^{\mathcal{A}_{k'}}$ its projection on the factor of index $k'$, for $0\leqslant k'<k$. The map $\mu_k$ is clearly bijective when $k \geqslant 2$.
 
 \medskip
\begin{thm}\label{Th6}
	There exists a unique $\bf Z$-linear map $\delta:{\bf Z}^{({\mathcal{B}})} \rightarrow {\bf Z}^{({\mathcal{A}})}$, graded of degree $0$ (for the weight), such that $\delta([\varnothing])=\varnothing$ and such that $\mu \circ \delta^*=\delta \circ \alpha$, where $\alpha:{\bf Z}^{(\mathcal{B}^*)}\rightarrow {\bf Z}^{(\mathcal{B})}$ is the $\bf Z$-linear map  defined in Section \ref{s1s11} and $\delta^*:{\bf Z}^{(\mathcal{B}^*)}\rightarrow{\bf Z}^{(\mathcal{A}^*)}$ is the map induced by $\delta$.
\end{thm}
A graded $\bf Z$-linear map $\delta:{\bf Z}^{(\mathcal{B})}\rightarrow{\bf Z}^{(\mathcal{A})}$ of degree $0$ satisfies the condition \mbox{ $\mu \circ \delta^*=\delta \circ \alpha$} if and only if it satisfies the following equivalent conditions:

a) We have 
\begin{equation}
\delta([{\bf a}])^{\rm init}=\delta([{\bf a}]^{\rm init}+[{\bf a}]^{\rm mid}+[{\bf a}]^{\rm fin})
\end{equation}
 for each non-empty admissible composition ${\bf a}$.

b) The homogeneous components $\delta_k$ of $\delta$ satisfy the relation
\begin{equation}\label{EU61}
\mu_k \circ \delta_k=\Big(\bigoplus_{0\leqslant k'<k}\delta_{k'}\Big)\circ \alpha_k\qquad \text{for $k \geqslant 1$.}
\end{equation}

c) We have 
\begin{equation}
\mu_{k',k}\circ \delta_k=\delta_{k'}\circ\alpha_{k',k} \qquad  {\rm for}\; 0\leqslant k'<k.
\end{equation}

The condition $\delta([\varnothing])= \varnothing$ determines $\delta_0$ uniquely, $\delta_1$ is the zero map since $\mathcal{B}_1$ is empty and then relation (\ref{EU61}) determines $\delta_k$ uniquely for $k \geqslant 2$ by induction on $k$, since the map $\mu_k $ is bijective when $k \geqslant 2$.

\medskip

{\it Remark.}- Note that the proof of Theorem \ref{Th6}  yields in fact an algorithm to compute $\delta_k$ combinatorially by induction on $k$.

\vspace{2em}
\subsection{Expressing multiple zeta values and their double tails in terms of multiple Ap{\'e}ry-like  sums and their tails}\label{s3s2}
\medskip
\myindent Let us recall that $\ell \rightarrow\zeta(\ell)_{n,n}$ denotes the ${\bf Z}$-linear map from ${\bf Z}^{(\mathcal{B})}$ to $\bf R$ that sends $[{\bf a}]$ to $\zeta({\bf a})_{n,n}$ for all ${\bf a}\in \mathcal{A}$ (see Section \ref{s1s8}). Similarly, let us denote by $\ell \rightarrow\sigma(\ell)_{n}$ the $\bf Z$-linear map from $\bf Z^{(\mathcal{A})}$ to $\bf R$ that sends $\bf a$ to $\sigma({\bf a})_{n}$ for ${\bf a}\in \mathcal{{A}}$.
\medskip

 \begin{thm}\label{T1}
 	For each $\ell \in {\bf Z}^{(\mathcal{B})}$ and each $n\geqslant0$, we have 
 	\begin{equation}\label{EF48}
 	\zeta(\ell)_{n,n}=\sigma(\delta(\ell))_{n}.
 	\end{equation}
 	In particular, we have
 	\begin{equation}\label{E16}
 	\zeta(\ell)=\sigma(\delta(\ell)).
 	\end{equation}
 \end{thm}
 
 By $\bf Z$-linearity, it suffices to prove the theorem when $\ell \in {\bf Z}^{\mathcal{B}_k}$ for some $k \geqslant 0$. We shall argue by induction on $k$.
 When $k=0$, our assertion follows from the fact that both $\zeta([\varnothing])_{n,n}$ and $\sigma(\varnothing)_n$ are equal to $\left({2n \atop n}\right)^{-1}$ by convention.

 Assume now $k \geqslant 1$. It follows from  formula (\ref{E12}) that $\zeta(\ell)_{n,n}$ tends to $0$ when $n$ goes to $+\infty$. Hence we deduce from the recurrence relation (\ref{EF20}) that we have
 \begin{equation}
 \zeta(\ell)_{n,n}=\sum_{m>n}\sum_{0\leqslant k'<k}m^{k'-k}\zeta(\alpha_{k',k}(\ell))_{m,m}
 \end{equation}
 and therefore, by the  induction hypothesis,
 \begin{equation}
 \begin{aligned}
  \zeta(\ell)_{n,n}=&\sum_{m>n}\sum_{0\leqslant k'<k}m^{k'-k}\sigma(\delta_{k'}(\alpha_{k',k}(\ell)))_m, \\ 
 =&\sum_{m>n}\sum_{0\leqslant k'<k} m^{k'-k}\sigma(\mu_{k',k}(\delta_k(\ell)))_m.
 \end{aligned}
 \end{equation}
 Hence, by Proposition \ref{P4}, we have
 \begin{eqnarray}
 \zeta(\ell)_{n,n}&=&\sum_{m>n}\big(\sigma(\delta_k(\ell))_{m-1}-\sigma(\delta_k(\ell))_{m}\big),
 \end{eqnarray}
 and since $\sigma(\delta_k(\ell))_m $ tends to $0$ when $m$ goes to $+\infty$,
 \begin{eqnarray}
 \zeta(\ell)_{n,n}&=& \sigma(\delta_k(\ell))_n.
 \end{eqnarray}
 
 {\it Remark.}- Let $\ell \in {\bf Z}^{(\mathcal{B})}$. Then $\delta(\ell)$ can be characterized as the unique element   $\sum_{ {\bf b}\in \mathcal{A}}\lambda_{\bf b}{\bf b}$ of ${\bf Z}^{(\mathcal{A})}$ such that $\zeta(\ell)_{n,n}=\sum_{ {\bf b}\in \mathcal{A}}\lambda_{\bf b}\sigma({\bf b})_{n}$ for all $n\geqslant 0$. Indeed, $\delta(\ell)$ satisfies this condition by  Theorem \ref{T1}, and the unicity follows from the linear independence of tails of multiple Ap{\' e}ry-like sums stated in Theorem \ref{T2} of Section \ref{s2s7}.
 
\vspace{2em}
 \subsection{Examples}\label{S3S3}
 \myindent A large number of identities found in the literature are consequences of Theorem \ref{T1} of Section \ref{s3s2}. Let us give here some examples.
 \subsubsection*{a) The case of weight 2}
 \myindent The set $\mathcal{B}_2$ has only one element, $[2]$. Since $\delta([2])$ is of weight $2$ and $\mu(\delta([2]))=\delta(\alpha([2]))=\delta(3[\varnothing])=3\varnothing$, we have  $\delta([2])=3(2)$. It implies that, for all $n \geqslant 0$,
 \begin{equation}\label{EU}
 \zeta(2)_{n,n}=3\sigma(2)_n,
 \end{equation}
 or more explicitely,
 \begin{equation}
 \sum_{k>n}\frac{1}{\left({k+n \atop n}\right)k^2}=3\sum_{k>n}\frac{1}{\left({2k \atop k}\right)k^2}\cdot
 \end{equation}
 In particular, when $n=0$, we get the famous formula, due to Euler,
 \begin{equation}\label{EU71}
 \zeta(2)=3\sigma(2),
 \end{equation}
 or more explicitely
  \begin{equation}
  \sum_{k=1}^{\infty}\frac{1}{k^2}=3\sum_{k=1}^{\infty}\frac{1}{\left({2k \atop k}\right)k^2}\cdot
  \end{equation}
  
  \medskip
\subsubsection*{b) The case of weight 3}
\myindent  The set $\mathcal{B}_3$ has only one element, $[3]=[2,1]$. Since $\delta([3])$ is of weight $3$ and $\mu(\delta([3]))=\delta(\alpha([3]))=\delta([2]+2[\varnothing])=3(2)+2\varnothing$, we have $\delta([3])=2(3)+3(2,1)$.
 
 It implies that, for all $n\geqslant 0$,
 \begin{equation}\label{EU1}
 \zeta(3)_{n,n}=2\sigma(3)_{n}+3\sigma(2,1)_n.
 \end{equation}
 Even the special case
 \begin{equation}\label{EU73}
 \zeta(3)=2\sigma(3)+3\sigma(2,1)
 \end{equation} 
  of the previous formula, {\it i. e.}
 \begin{equation}\label{E68}
 \sum_{k=1}^{\infty}\frac{1}{k^3}=2\sum_{k=1}^{\infty}\frac{1}{\left({2k \atop k}\right)k^3}+3 \sum_{k_1>k_2\geqslant1}\frac{1}{\left({2k_1 \atop k_1}\right)k_1^2k_2}
 \end{equation}
 does not seem to have been noticed previously.
 
 \medskip
  {\it  Remarks.-} 1)  It is very unlikely that $\sigma(3)$ and $\sigma(2,1)$ are rational multiples of $\zeta(3)$. Indeed we computed $\sigma(3)/\zeta(3)$ with 2500 exact digits and expanded it  as a continued fraction
 \begin{equation}
 \frac{\sigma(3)}{\zeta(3)}= \cfrac{1}{2+\cfrac{1}{3+\cfrac{1}{2+\cfrac{1}{1+\cfrac{1}{6+\cdots}}}}}\;\cdot \end{equation}
 We checked that the 2000th  convergent $\frac{p}{q}$ of this continued fraction  satisfies the inequalities $|\frac{\sigma(3)}{\zeta(3)}-\frac{p}{q}|<\frac{1}{q^2}$ and $q>10^{1000}$. Therefore, if $\frac{\sigma(3)}{\zeta(3)}$ were a rational number, its  denominator would be larger than $10^{1000}$.
 
 \medskip
 2) I. Zucker has shown in \cite{Zuk}, formula (2.9), that
 \begin{equation}\label{EU110}
 \begin{aligned}
 \sigma(3)&=\frac{\pi \sqrt{3}}{2}{\rm L}(2,\chi)-\frac{4}{3}\zeta(3),
 \end{aligned}
 \end{equation}
 where $\chi$ is the non-trivial Dirichlet character of conductor $3$.
Hence by (\ref{EU73}), we have
 \begin{equation}\label{EU111}
 \begin{aligned}
 \sigma(2,1)&=-\frac{\pi \sqrt{3}}{3}{\rm L}(2,\chi)+\frac{11}{9}\zeta(3).
 \end{aligned}
 \end{equation}
 Note that  (\ref{EU110}) is a particular case of Theorem~\ref{TH8}, and (\ref{EU111}) is a particular case of Theorem~\ref{TH7}.

 3) R. Ap{\' e}ry has proved in \cite{Apery} that $\zeta(3)$ is irational, starting from the fact that it is equal to $\frac{5}{2} \tilde{\sigma}(3)$ where $\tilde{\sigma}(3)=\sum_{n=1}^{\infty}\frac{(-1)^{n-1}}{\left({2n \atop n}\right)n^3}$. Can his method be extended to prove that $\sigma(3)$ and $\sigma(2,1)$ are irrational? If yes, can the linear independence over {\bf Q} of $1$, $\sigma(3)$, and $\sigma(2,1)$ (or equivalently of $1$, $\zeta(3)$ and $\pi \sqrt{3} {\rm L}(2,\chi)$) be proved by the techniques devoloped  by T. Rivoal in \cite{Rivoal}?
 
 \medskip

 \subsubsection*{c) The case of weight 4}
 \myindent The set $\mathcal{B}_4$ has three elements, $[4]=[2,1,1],\;[3,1]$ and $[2,2]$.
 
 \medskip
 \begin{lem}\label{L}
 	We have
 	\begin{equation}
 	\begin{aligned}
 	\delta([4])&= 2\;(4)+2\;(3,1)+3\;(2,1,1),\\ 
 	\delta([3,1])&= 4\;(3,1)+3\;(2,2)+6\;(2,1,1),\\ 
 	\delta([2,2])&= (4)+6\;(2,2).
 	\end{aligned}
 	\end{equation}
 \end{lem}
 
 The lemma indeed follows from the fact that $\delta([4]),\delta([3,1])$ and $\delta([2,2])$ are the elements of ${\bf Z}^{\mathcal{A}_4}$ satisfying the following identities:
 \begin{equation}
 \begin{aligned}
  \mu(\delta([4]))&=\delta(\alpha([4]))=\delta([3]+2[\varnothing])=2(3)+3(2,1)+2\varnothing, \\ 
 \mu(\delta([3,1]))&=\delta(\alpha([3,1]))=\delta(2[3]+[2])=4(3)+6(2,1)+3(2),\\ 
 \mu(\delta([2,2]))&=\delta(\alpha([2,2]))=\delta(2[2]+[\varnothing])=6(2)+\varnothing.
 \end{aligned}
 \end{equation}
 
 \medskip
 
 From Lemma \ref{L} and Theorem \ref{T1} we deduce that, for all $n \geqslant 0$,
  \begin{equation}
  \begin{aligned}
  \zeta(4)_{n,n}=\zeta(2,1,1)_{n,n}&=\; 2\; \sigma(4)_n + 2\; \sigma(3,1)_n +3\; \sigma(2,1,1)_n,\\
  \zeta(3,1)_{n,n}&=\;4\; \sigma(3,1)_n + 3\; \sigma(2,2)_n +6 \; \sigma(2,1,1)_n,\\
  \zeta(2,2)_{n,n}&=\; \sigma(4)_n+ 6\; \sigma(2,2)_n.
  \end{aligned}
  \end{equation}
  From these relations , we deduce that, for all $n\geqslant 0$,
  \begin{equation}
  \begin{aligned}
  \sigma(4)_n&=\frac{1}{9}\left(4\zeta(4)_{n,n}-2\zeta(3,1)_{n,n}+\zeta(2,2)_{n,n}\right),\\
  \sigma(2,2)_n&=\frac{1}{27}\left(-2\zeta(4)_{n,n}+\zeta(3,1)_{n,n}+4\zeta(2,2)_{n,n}\right),\\
  2\sigma(3,1)_n+3\sigma(2,1,1)_n&=\frac{1}{9}\left(\zeta(4)_{n,n}+4\zeta(3,1)_{n,n}-2\zeta(2,2)_{n,n}\right).
    \end{aligned}
  \end{equation}
  
  When $n=0$, these relations, together with the well known identities
  \begin{equation}
 \zeta(3,1)=\frac{1}{4}\zeta(4),\quad \zeta(2,2)=\frac{3}{4}\zeta(4),\quad \zeta(4)=\frac{\pi^4}{90},
  \end{equation} 
   yield
    \begin{equation}\label{E20}
    \begin{aligned}
    \sigma(4)&=\frac{17}{36}\zeta(4) &=\frac{17 \pi^4}{3240},\\
    \sigma(2,2)&= \frac{5}{108}\zeta(4)&=\frac{\pi^4}{1944},\\
    2\;\sigma(3,1)+3\; \sigma(2,1,1)&=\frac{1}{18}\zeta(4)&=\frac{\pi^4}{1620}.
    \end{aligned}
    \end{equation}
    The first of these three equalities appears in \cite{Comtet}, p.89, whereas the second one can be found in \cite{Wen}, p.431 and we were not able to locate the third one anywhere in the literature.
    
    \medskip
    
 {\it  Remark.-}  It is very unlikely that  each number $\sigma(3,1)$ and $\sigma(2,1,1)$ is a rational multiple of $\pi^4$. In the same way as in the remark of example b), we have checked that, were for example   $\frac{\sigma(3,1)}{\pi^4}$ a rational number, its  denominator would be larger than $10^{1000}$.
 \vspace{2em}
 \subsubsection*{d) The case of weight 5}
 \myindent The set $\mathcal{B}_5$ has four elements $[5]=[2,1,1,1]$, $[4,1]=[3,1,1]$, $[3,2]=[2,2,1]$ and $[2,3]=[2,1,2]$. Their images by $\delta$ are computed as in the example c). We get:
 \begin{equation}\label{EU85}
 \begin{aligned}
 \delta([5])&=2(5)+2(4,1)+2(3,1,1)+3(2,1,1,1),\\
\delta([4,1])&=2(4,1)+2(3,2)+6(3,1,1)+3(2,2,1)+3(2,1,2)+9(2,1,1,1),\\
 \delta([3,2])&=(4,1)+2(3,2)+3(2,3)+6(2,2,1)+3(2,1,2),\\
 \delta([2,3])&=(5)+2(3,2)+3(2,3)+3(2,1,2). 
 \end{aligned}
 \end{equation}
 We leave to the reader the task of expressing the corresponding multiple zeta values and their double tails in terms of Ap{\'e}ry-like sums and their tails.
 
 Any relation between multiple zeta values yields  a relation between Ap{\'e}ry-like sums: if $\ell \in {\bf Z}^{(\mathcal{B})}$ is such that $\zeta(\ell)=0$, then $\ell'=\delta(\ell)$ is an element of ${\bf Z}^{(\mathcal{A})}$ such that $\sigma(\ell')=0$. From the two relations
 \begin{equation}
 \begin{aligned}
 &\zeta(5)-\zeta(4,1)&=\zeta(3,2)+\zeta(2,3),\quad\\
 &5\zeta(4,1)&=-\zeta(3,2)+\zeta(2,3),\;
  \end{aligned}
 \end{equation}
 we deduce two independent $\bf Z$-linear relations between Ap{\' e}ry-like sums of weight $5$:
 \begin{equation}\label{EU87}
 \begin{aligned}
 \sigma(5)&=\sigma(4,1)+6\sigma(3,2)+4\sigma(3,1,1)+6\sigma(2,3)+9\sigma(2,2,1)\\&\quad+9\sigma(2,1,2)+6\sigma(2,1,1,1),\\
 \sigma(5)&=11\sigma(4,1)+10\sigma(3,2)+30\sigma(3,1,1)+21\sigma(2,2,1)+15 \sigma(2,1,2)\\&\quad+45\sigma(2,1,1,1).
 \end{aligned}
 \end{equation}
 But by PSLQ algorithm, we discovered a third apparent one:
 \addtocounter{equation}{+1}
 \begin{equation}\label{EU88}
 4\sigma(4,1)=6\sigma(2,2,1)+22\sigma(3,1,1)+33\sigma(2,1,1,1). 
 \end{equation}
 We didn't find any simple direct proof of it, but we can deduce it from the evaluations given in Section~\ref{s2s8}:
 
 \medskip
 \begin{prop}
 	Formula ${(\rm \ref{EU88})}$ holds.
 \end{prop}
 Indeed, we deduce from Theorem \ref{TH7} that
 \begin{equation}\label{EU127}
 \left(\begin{array}{c}
 \sigma(3,2)\\
 \sigma(2,3)\\
 \sigma(2,2,1)\\
 \sigma(2,1,2)
 \end{array}\right)=\left(\begin{array}{c c c c}
 \frac{1}{108}& -\frac{3}{8} & \frac{1}{27}& \frac{29}{27}\\
 0&-\frac{9}{8}&-\frac{2}{27}&  \frac{58}{9}\\
 0 &  \frac{1}{3}& \frac{1}{6}& -\frac{575}{162}\\
 -\frac{1}{162}&1 & -\frac{8}{81}& -\frac{575}{162}
 \end{array}\right)\;\left(\begin{array}{c}
 \pi^3\sqrt{3}\,{\rm L}(2,\chi)\\
  \pi\sqrt{3}\,{\rm L}(4,\chi)\\
  \pi^2 \zeta(3)\\
 \zeta(5)
 \end{array}\right)
 \end{equation}
 holds, where $\chi$ is the non principal Dirichlet character modulo $3$. On the other hand, one knows how to express the multiple zeta values of weight $5$ as $\bf Q$-linear combination of $\pi^2\zeta(3)$ and $\zeta(5)$: we have
 \begin{equation}\label{EU128}
 \left(\begin{array}{c}
 \zeta(5)\\
 \zeta(4,1)\\
 \zeta(3,2)\\
 \zeta(2,3)
 \end{array}\right)=\left(\begin{array}{c c}
 0 & 1\\
-\frac{1}{6}& 2\\
\frac{1}{2}& -\frac{11}{2}\\
-\frac{1}{3}& \frac{9}{2} 
 \end{array}\right)\; \left(\begin{array}{c}
 \pi^2\zeta(3)\\
 \zeta(5)
  \end{array}\right).
 \end{equation}
 By using (\ref{EU127}) and (\ref{EU128}) and the relations between multiple zeta values and multiple Ap{\'e}ry-like sums deduced from (\ref{EU85}), we get
 \begin{equation}\label{EU129}
 \left(\begin{array}{c}
 \sigma(5)\\
 \sigma(4,1)\\
 2\sigma(3,1,1)+2\sigma(2,1,1,1)
 \end{array}\right)=\left(\begin{array}{c c c}
\frac{9}{8} & \frac{1}{9 } & -\frac{19}{3}\\
 -\frac{7}{8} &-\frac{1}{18} &  \frac{134}{27}\\
-\frac{1}{2} & -\frac{1}{9} &  \frac{101}{27}
  \end{array}\right)\; \left(\begin{array}{c}
  \pi\sqrt{3}\,{\rm L}(4,\chi)\\
  \pi^2 \zeta(3)\\
    \zeta(5)
    \end{array}\right).
 \end{equation}
 Now relation (\ref{EU88}) is an immediate consequence of (\ref{EU127}) and (\ref{EU129}).
 \vspace{2em}
 \subsubsection*{e) The case of depth 1}
 \begin{thm}\label{Th8}
 	For each integer $a\geqslant 2$, we have
 	\begin{equation}\label{E75}
 	\delta([a])=2\;\sum_{b=3}^{a}(b,\underbrace{1,\ldots ,1}_{a-b})+3\;(2,\underbrace{1,\ldots ,1}_{a-2})
 	\end{equation}
 \end{thm}
 
 The proof is by induction on $a$. When $a=2$, relation (\ref{E75}) becomes $\delta([2])=3(2)$ and has been proved above  in Subsection a).

  Assume $a\geqslant 3$. We have $(a)^{\rm init}=\varnothing,(a)^{\rm mid}=\varnothing, (a)^{\rm fin}=(a-1)$, hence $\alpha([a])=2[\varnothing]+[a-1]$ and 
 \begin{equation}
 \begin{aligned}
 \delta(\alpha([a]))&=2\varnothing+\delta([a-1])\\
 &=2\varnothing+2\sum_{b=3}^{a-1}(b,\underbrace{1,\ldots ,1}_{a-b-1})+3(2,\underbrace{1,\ldots ,1}_{a-3}).
 \end{aligned}
 \end{equation}
 Since $\delta([a])$ is the unique element of ${\bf Z}^{\mathcal{A}_a}$ which is mapped by $\mu$ to $\delta(\alpha([a]))$, the theorem follows.
 
 \medskip
 \begin{cor}\label{CT2}
 	For each integer $a\geqslant 2$, we have
 	\begin{equation}
 	\zeta(a)=2\;\sum_{b=3}^{a}\sigma(b,\underbrace{1,\ldots ,1}_{a-b})+3\;\sigma(2,\underbrace{1,\ldots ,1}_{a-2}).
 	\end{equation}
 \end{cor}
 Note that this corollary generalizes  formulas (\ref{EU71}) and (\ref{EU73}). It does not seem to have been stated previously in the literature.

\medskip
\begin{cor}\label{C2}
	Let $x$ be an indeterminate. In the ring of formal power series ${\bf R}[[x]]$, endowed with the product topology of ${\bf R}^{\bf N}$, we have
	\begin{equation}\label{EN78}
	\sum_{k=1}^{\infty}\frac{1}{k(k-x)}=\sum_{k=1}^{\infty}\frac{3k-x}{\left(2k \atop k\right)k^2(k-x)}\prod_{m=1}^{k-1}\left(1+\frac{x}{m}\right).
	\end{equation}
\end{cor}

Indeed, the left-hand side of (\ref{EN78}) is equal to
\begin{equation}
\sum_{k=1}^{\infty}\sum_{a=2}^{\infty}\frac{x^{a-2}}{k^a}=\sum_{a=2}^{\infty}\zeta(a)x^{a-2},
\end{equation}
whereas the right-hand side is equal to 
\begin{equation}
\begin{aligned}
&\sum_{k=1}^{\infty}\frac{1}{\left({2k \atop k}\right)}\Big(\frac{3}{k^2}+2\sum_{b=3}^{\infty}\frac{x^{b-2}}{k^b}\Big)\sum_{\substack{r\geqslant 0\\k>m_1>\ldots >m_r>0}}\frac{x^r}{m_1\ldots m_r}
\\=\;&\sum_{a=2}^{\infty}\Big(3\; \sigma(2,\underbrace{1,\ldots ,1}_{a-2})+2\;\sum_{b=3}^{a}\sigma(b,\underbrace{1,\ldots ,1}_{a-b})\Big)\;x^{a-2},
\end{aligned}
\end{equation}
hence the Corollary \ref{C2} is a translation, in terms of generating series, of Corollary \ref{CT2}.

\medskip
{\it Remarks.}- 1) Alternatively, if one consider $x$ as a complex variable instead of an indeterminate, both sides of (\ref{EN78}) converge normally on each compact subset of the open set ${\bf C}-{\bf N}^*$ and their sums are equal in this open set. Note that if $f(x)$ denotes this sum, we have $f(-x)=\frac{1}{x}\left(\psi(x)+\frac{1}{x}+\gamma\right)$, where $\psi$ is the logarithmic derivative of the Gamma function and $\gamma$ is Euler's constant.

\medskip
 2) Theorem \ref{Th8} implies that we have more generally, for all integers $a\geqslant 2$ and $n \geqslant 0$,
\begin{equation}
 	\zeta(a)_{n,n}=2\;\sum_{b=3}^{a}\sigma(b,\underbrace{1,\ldots ,1}_{a-b})_n+3\;\sigma(2,\underbrace{1,\ldots ,1}_{a-2})_n.
\end{equation}
As above, we have for each $n \geqslant 0$, the following translation of these identities in terms of generating series
\begin{equation}
\sum_{k=n+1}^{\infty}\frac{1}{\left(k+n \atop n\right)k(k-x)}=\sum_{k=n+1}^{\infty}\frac{3k-x}{\left(2k \atop k\right)k^2(k-x)}\prod_{m=n+1}^{k-1}\left(1+\frac{x}{m}\right).
\end{equation}
\vspace{1.5 em}
  \subsubsection*{f) Some identities experimentally discovered by D. Bailey,  J. Borwein, D.Bradley}
   				\begin{thm}\label{Th9}For each even integer $k\geqslant 0$, we have
  					\begin{equation}\label{E78}
  					\delta \left(\sum_{ {\bf a}\in \mathscr{A}_{k}}(-1)^{{\rm depth}({\bf a})}[\bf{a}]\right)=\sum_{{\bf a}\in \mathscr{A}_{k}^{\rm even}}(-3)^{{\rm depth}({\bf a})}{\bf a},
  					\end{equation}
  					where $\mathcal{A}_k^{\rm even}$ is the set of admissible compositions ${\bf a}$ of weight $k$ with  even entries.
  				\end{thm}
For each  (not necessarily even) integer $p\geqslant0$, let $\ell_p$ denote the element of ${\bf Z}^{\mathcal{A}_p}$ defined by 
\begin{equation}
\ell_p=\sum_{ {\bf a}\in \mathcal{A}_p}(-1)^{{\rm depth}({\bf a})}{\bf a}.
\end{equation}  

We have in particular $\ell_0=\varnothing$ and $\ell_1=0$. Since ${\bf a}\rightarrow \overline{\bf a}$ is an involution of $\mathcal{A}_p$  and since   ${\rm depth}({\bf a})+{\rm depth}({\overline{\bf a}})=p$ for ${\bf a}\in \mathcal{A}_p$, we have 
\begin{equation}\label{E80}
\overline{\ell_p}=(-1)^p\ell_p.
\end{equation}

\smallskip
\begin{lem}
	For each even integer $k\geqslant 2$, we have
	\begin{equation}
	\alpha([\ell_k])=-3\sum_{\substack{0\leqslant q \leqslant k-2\\ q\; {\rm even}}}[\ell_q].
	\end{equation}
\end{lem}

We have $\ell_0^{\rm init}=\ell_0$ and for $p\geqslant 2$

\begin{equation}\label{E81}
\ell_p^{\rm init}=-\sum_{0\leqslant q<p}\ell_{q}=-\left(\ell_0+\sum_{q=2}^{p-1}\ell_{q}\right).
\end{equation}				

For $k\geqslant 2$ even, we have therefore
\begin{equation}\label{E82}
\begin{aligned}
\ell_k^{\rm fin}=\overline{\overline{\ell_k}^{\rm init}}=\overline{\ell_k^{\rm init}}=-\left(\ell_0+\sum_{p=2}^{k-1}(-1)^p\ell_{p}\right).
\end{aligned}
\end{equation}
and
\begin{equation}
\ell_k^{\rm mid}=\left(\ell_k^{\rm fin}\right)^{\rm init}=-\ell_0+\sum_{p=2}^{k-1}(-1)^p\left(\ell_0+\sum_{q=2}^{p-1}\ell_q\right)=-\sum_{\substack{0\leqslant q\leqslant k-2\\ q \;{\rm even}}}\ell_q.
\end{equation}	

We hence have
\begin{equation}
\ell_k^{\rm init}+\ell_k^{\rm mid}+\ell_k^{\rm fin}=-3\sum_{\substack{0\leqslant q\leqslant k-2\\ q \;{\rm even}}}\ell_q.
\end{equation}
This proves the lemma.

\medskip
Let us now prove Theorem \ref{Th9} by induction on $k$. When $k=0$, it follows from the relation $\delta([\varnothing])=\varnothing$. When $k\geqslant 2$ is even, we deduce from the lemma and the induction hypothesis that
\begin{equation}
\delta([\ell_k])^{\rm init}=\delta(\alpha([\ell_k]))=-3\sum_{\substack{0\leqslant q\leqslant k-2\\q\text{  even}}}\;\sum_{ {\bf b}\in \mathcal{A}_{q}^{\rm even }}(-3)^{{\rm depth}({\bf b})}{\bf b}
\end{equation}
and hence that
\begin{equation}
\delta([\ell_k])=\sum_{{\bf a}\in \mathscr{A}_{k}^{\rm even}}(-3)^{{\rm depth}({\bf a})}{\bf a}.
\end{equation}

\medskip	
\begin{customcor}{1}\label{CT1} For each even integer $k \geqslant 2$, we have
\begin{equation}\label{E88}
\zeta(k)=-\sum_{{\bf a}\in \mathscr{A}_{k}^{\rm even}}(-3)^{{\rm depth}({\bf a})}\sigma({\bf a}).
\end{equation}
\end{customcor}

When we apply Theorem \ref{T1} to  relation (\ref{E78}), we get, for each even integer $k\geqslant 2$,
\begin{equation}\label{E89}
\sum_{ {\bf a} \in \mathcal{A}_k}(-1)^{ {\rm depth}({\bf a})}\zeta({\bf a})=\sum_{ {\bf a} \in \mathcal{A}_k^{\rm even}} (-3)^{ {\rm depth}({\bf a})}\sigma({\bf a}).
\end{equation}

But it has been proved by A. Granville  in \cite{Grand}, as well as by Y. Ohno  and D. Zagier  in \cite{OZ}, that for each integer $r$ such that $1 \leqslant r \leqslant k-1$, we have
\begin{equation}
\sum_{ \substack{{\bf a} \in \mathcal{A}_k\\ {\rm depth}({\bf a})=r}}\zeta({\bf a})=\zeta(k).
\end{equation}
Hence the left-hand side of (\ref{E89}) is equal to $\sum_{r=1}^{k-1}(-1)^r
\zeta(k)=-\zeta(k)$, and the corollary follows.

\medskip

{\it Remarks.}- 1) The identities stated in Corollary \ref{CT1} were first discovered experimentally by D. Bailey, J. Borwein and D. Bradley in 2006 (see \cite{BBB}). Their result was stated in a slightly different form: for each non-empty sequence $(a_1,\ldots ,a_r)$ of positive even integers, they introduced the real number
\begin{equation}
\sigma(a_1;[a_2,\ldots ,a_r])=\sum_{n_1=1}^{\infty} \frac{1}{\left({2n_1 \atop n_1}\right)n_1^{a_1}}\prod_{i=2}^{r}\left(\sum_{n_i=1}^{n_1-1}\frac{1}{n_i^{a_i}}\right).
\end{equation}
They stated then their result, in weight $6$ for example, in the form
\begin{equation}\label{E91}
\zeta(6)=3\sigma(6;[])-9\sigma(4;[2])-\frac{45}{2}\sigma(2;[4])+\frac{27}{2}\sigma(2;[2,2]).
\end{equation}\label{EE103}
Since their $\sigma(6;[])$, $\sigma(4;[2])$, $\sigma(2;[4])$ and $\sigma(2;[2,2])$ are clearly equal to our $\sigma(6)$, $\sigma(4,2),$ $\sigma(2,4)$ and $2\sigma(2,2,2)+\sigma(2,4)$, their formula (\ref{E91}) is equivalent to our formula
\begin{equation}
\zeta(6)=3 \sigma(6)-9\sigma(4,2)-9\sigma(2,4)+27 \sigma(2,2,2).
\end{equation}
However, our formulas look simpler and more symmetric than theirs.

\medskip
\begin{customcor}{2}\label{C2T1}
		Let $x$ be an indeterminate. In the ring of formal power series ${\bf R}[[x]]$, endowed with the product topology of ${\bf R}^{\bf N}$, we have
	\begin{equation}\label{E93}
	\sum_{k=1}^{\infty}\frac{1}{k^2-x^2}=3\sum_{k=1}^{\infty}\frac{1}{\left({2k \atop k}\right)(k^2-x^2)	}\prod_{m=1}^{k-1}\frac{m^2-4x^2}{m^2-x^2}\cdot
	\end{equation}
\end{customcor}

\medskip
Indeed, the left-hand side of (\ref{E93}) can be written as
\begin{equation}
\sum_{k=1}^{\infty}\sum_{n=1}^{\infty}\frac{x^{2n-2}}{k^{2n}}=\sum_{n=1}^{\infty}\zeta(2n)x^{2n-2},
\end{equation}
whereas the right-hand side of (\ref{E93}) is equal to
\begin{equation*}
\begin{aligned}
&3\sum_{k=1}^{\infty}\frac{1}{\left({2k \atop k}\right)(k^2-x^2)}\prod_{m=1}^{k-1}\Big(1-3\sum_{a=1}^{\infty}\frac{x^{2a}}{m^{2a}}\Big)
\\=\;&3\sum_{k=1}^{\infty}\frac{1}{\left({2k \atop k}\right)}\sum_{a_1 =1}^{\infty}\frac{x^{2a_1-2}}{k^{2a_1}}\sum_{\substack{r\geqslant 1\\ a_2,\ldots ,a_r\geqslant 1}}\sum_{k>m_2>\ldots >m_r>0}\Bigg(\frac{-3x^{2a_2}}{m_2^{2a_2}}\Bigg)\ldots \Bigg(\frac{-3x^{2a_r}}{m_r^{2a_r}}\Bigg)
\\=\;&-\sum_{\substack{r\geqslant 1\\ a_1,\ldots ,a_r\geqslant 1}}(-3)^{r}\sigma(2a_1,\ldots ,2a_r)x^{2a_1+\ldots +2a_r-2}
=\sum_{n=1}^{\infty}\Bigg(-\sum_{ {\bf a}\in \mathcal{A}_{2n}^{\rm even}}(-3)^{{\rm depth}({\bf a})}\sigma({\bf a})\Bigg)x^{2n-2}.
\end{aligned}
\end{equation*}
Hence Corollary \ref{C2T1} follows from Corollary \ref{CT1}.

\medskip

{\it  Remarks.}- 2) In fact, D. Bailey, J. Borwein and D. Bradley first searched experimentally formulas similar to (\ref{E91}) for small even weights $k$. They found such formulas \mbox{for $k\leqslant 10$.} They then searched, with the help of Pad{\' e} approximants, generating series to guess how these formulas should  generalize to higher weights. In this way, they were led to conjecture that (\ref{E93}) holds, as an equality between meromorphic functions of $x$. Finally they succeeded to prove this equality, by analysing the growth of the two sides at infinity, and comparing their poles and residues. They reduced the proof of the equality of the residues to some identities about finite generalized hypergeometric series. These identities were in turn proved by applying Wilf-Zeilberger algorithm, as implemented in \mbox{MAPLE 9.5.}  We consider that the account given here provides  an elementary and natural proof of the same formulas.
\medskip

3) Theorem \ref{T1} applied to  relation (\ref{E78}) yields, for each even integer $k \geqslant 2$ and each integer $n \geqslant 0$, the identity
\begin{equation}\label{E97}
\sum_{ {\bf a} \in \mathcal{A}_k}(-1)^{{\rm depth}({\bf a})}\zeta({\bf a})_{n,n}=\sum_{ {\bf a} \in \mathcal{A}_k^{\rm even}}(-3)^{ {\rm depth}({\bf a})}\sigma({\bf a})_n,
\end{equation}
which generalizes (\ref{E89}). Moreover the left hand side of (\ref{E97}) vanishes for odd integers~$k$, since \mbox{${\rm depth}({\bf a})+{\rm depth}(\overline{\bf a})=k$} and $\zeta(\overline{\bf a })_{n,n}=\zeta({\bf a})_{n,n}$ for ${\bf a}\in \mathcal{A}_k$. In the same way as we derived Corollary \ref{C2T1} from Corollary \ref{CT1}, one deduces from all these facts that, for all $n\geqslant 0$,   the following equality of generating series holds: 
\begin{equation}
\sum_{k=n+1}^{\infty}\frac{1}{\left(k+n \atop n\right)\;k\;(k-x)}\prod_{m=n+1}^{k-1}\frac{m-2x}{m-x}=3\sum_{k=n+1}^{\infty}\frac{1}{\left(2k \atop k\right)(k^2-x^2)}\prod_{m=n+1}^{k-1}\frac{m^2-4x^2}{m^2-x^2}\cdot
\end{equation}

\medskip
\subsection*{g) Compositions of the form $(2,\ldots,2)$}
\begin{thm}\label{TH10}
	For each $m \geqslant 0$, we have
	\begin{equation}\label{EN103}
	\delta([\underbrace{2,\ldots,2}_{m}])=\sum_{{\bf a}\in \mathcal{A}_{2m}^{\{2,4\}}}c_{\bf a} {\bf a},
	\end{equation}
where $\mathcal{A}_{2m}^{\{2,4\}}$ is the set of compositions of weight $2m$ whose entries all belong to $\{2,4\}$ and where, for such a composition ${\bf a}=(a_1,\ldots,a_r),$
\begin{equation}
c_{\bf a}=\begin{cases}
3\cdot2^{s-1} & \text{ if } a_1=2,\\
2^s& \text{ otherwise},
\end{cases}
\end{equation}
where $s=2r-m$ is the number of indices $i \in \{1,\ldots,r\}$ such that $a_i=2$.
\end{thm}

The proof is by induction on $m$. For $m=0$ and $m=1$, (\ref{EN103}) follows from the equalities $\delta([\varnothing])=\varnothing$ and $\delta([2])=3(2)$. When $m \geqslant2$, we have $\alpha([\underbrace{2,\ldots,2}_{m}])=2[\underbrace{2,\ldots,2}_{m-1}]+[\underbrace{2,\ldots,2}_{m-2}]$, hence $\mu\Big(\delta([\underbrace{2,\ldots,2}_{m}])\Big)=2\delta([\underbrace{2,\ldots,2}_{m-1}])+\delta([\underbrace{2,\ldots,2}_{m-2}])$, and the lemma follows from the induction hypothesis.

\medskip
We get as corollary the following formulas, of the same flavour as that those obtained by D. Bailey, J. Borwein and D. Bradley, but which now hold also at the level of tails:

\medskip

\begin{customcor}{1} \label{CR1}
	For all natural numbers $m$ and $n$, we have
	\begin{equation}
	\zeta([\underbrace{2,\ldots,2}_{m}])_{n,n}=\sum_{{\bf a}\in \mathcal{A}_{2m}^{\{2,4\}}}c_{\bf a} \sigma({\bf a})_n,
	\end{equation}
	where the coefficients $c_{\bf a}$ are as in Theorem \ref{TH10}. We have in particular
	\begin{equation}
	\frac{\pi^{2m}}{(2m+1)!}=\sum_{{\bf a}\in \mathcal{A}_{2m}^{\{2,4\}}}c_{\bf a} \sigma({\bf a}).
	\end{equation}
\end{customcor}

Last formula follows from the classical equality $\zeta([\underbrace{2,\ldots,2}_{m}])=\frac{\pi^{2m}}{(2m+1)!}$.

\medskip

As in examples e)  and f), we can translate Corollary \ref{CR1} in a statement about generating series:

\medskip
\begin{customcor}{2}
		Let $x$ be an indeterminate. In the ring of formal power series ${\bf R}[[x]]$, endowed with the product topology of ${\bf R}^{\bf N}$, we have
\begin{equation}
\frac{{\rm sinh} (\pi x)}{\pi x}=1+x^2\sum_{k=1}^{\infty}\frac{3k^2+x^2}{\left(2k \atop k\right)k^4}\prod_{m=1}^{k-1}\left(1+\frac{x^2}{m^2}\right)^2,
\end{equation}
and more generally, for all $n \geqslant 0$,
\begin{equation}
\sum_{k=n+1}^{\infty}\frac{1}{\left(k+n \atop n\right)k^2}\prod_{m=n+1}^{k-1}\Big(1+\frac{x^2}{m^2}\Big)=\sum_{k=n+1}^{\infty}\frac{3k^2+x^2}{\left(2k \atop k\right)k^4}\prod_{m=n+1}^{k-1}\left(1+\frac{x^2}{m^2}\right)^2.
\end{equation}

\end{customcor}
\vspace{2em}

\subsection*{h) An identity of D. Leshchiner}

 \begin{thm}\label{Th11}
 	For each even weight $k=2m\geqslant2$
 	\begin{equation}\label{E109}
 	\delta\Big(\sum_{b=2}^{k}(-1)^b[b,\underbrace{1,\ldots,1}_{k-b}]\Big)=3(-1)^{m-1}(\underbrace{2,\ldots,2}_{m})+4\sum_{c=2}^{m}(-1)^{m-c}(2c,\underbrace{2,\ldots 
 		,2}_{m-c}).
 	\end{equation}
 \end{thm}
 The proof of the theorem is by induction on $k$. When $k=2$, we have $\delta([2])=3(2)$, hence (\ref{E109}) holds.
 
 We assume now that $k$ is even and $k\geqslant 4$. We have
 \begin{equation}
 \begin{aligned}
(b,\underbrace{1,\ldots,1}_{k-b})^{\rm init} &=\begin{cases}
(b,\underbrace{1,\ldots,1}_{k-1-b}) & \qquad\text{ if } 2 \leqslant b \leqslant k-1,\\
\varnothing & \qquad\text{ if } b=k,
\end{cases}
\\
(b,\underbrace{1,\ldots,1}_{k-b})^{\rm fin} &=\begin{cases}
(b-1,\underbrace{1,\ldots,1}_{k-b}) & \;\text{ if } 3 \leqslant b \leqslant k,\\
\varnothing & \;\text{ if } b=2,
\end{cases}
\\
(b,\underbrace{1,\ldots,1}_{k-b})^{\rm mid} &=\begin{cases}
(b-1,\underbrace{1,\ldots,1}_{k-1-b}) & \;\text{ if } 3\leqslant b \leqslant k-1,\\
\varnothing & \;\text{ if } b=2 \text{ or } b=k.
\end{cases}
 \end{aligned}
   \end{equation}
   Summing up, we get
   \begin{equation}
   \alpha(\sum_{b=2}^{k}(-1)^b[b,\underbrace{1,\ldots,1}_{k-b}])=4\varnothing+\sum_{b=3}^{k-1}(-1)^b[b-1,\underbrace{1,\ldots,1}_{k-1-b}]=4\varnothing-\sum_{b=2}^{k-2}(-1)^b[b,\underbrace{1,\ldots,1}_{k-2-b}].
   \end{equation}
   Hence
   \begin{equation}
   \begin{aligned}
  &\mu\Big(\delta\Big(\sum_{b=2}^{k}(-1)^b[b,\underbrace{1,\ldots,1}_{k-b}]\Big)\Big)=4\varnothing-\delta\Big(\sum_{b=2}^{k-2}(-1)^b[b,\underbrace{1,\ldots,1}_{k-2-b}]\Big)\\
  &\qquad\qquad\qquad=4\varnothing+3(-1)^{m-1}(\underbrace{2,\ldots,2}_{m-1})+4\sum_{c=1}^{m-1}(-1)^{m-c}(2c,\underbrace{2,\ldots,2}_{m-1-c})
     \end{aligned}
      \end{equation}
      by the induction hypothesis. Theorem \ref{Th11} follows.
      
      \medskip

 \begin{customcor}{1}\label{ccr1} For all integers $n \geqslant 0$ and even integers $k=2m \geqslant 2$, we have
 	 	\begin{equation}
 	 	\sum_{b=2}^{k}(-1)^b\zeta(b,\underbrace{1,\ldots,1}_{k-b})_{n,n}=3(-1)^{m-1}\sigma(\underbrace{2,\ldots,2}_{m})_n+4\sum_{c=2}^{m}(-1)^{m-c}\sigma(2c,\underbrace{2,\ldots 
 	 		,2}_{m-c})_n.
 	 	\end{equation}
  	 \end{customcor}
  	 As a consequence, we recover the following identity, due to Leshchiner (see \cite{Lesh}, formula (3a), in which however there is a minor flaw in the exponent of $2$ of the left hand side):
  	 
  	 \medskip
  	 \begin{customcor}{2}\label{ccr2}
  	 	For all even integers $k=2m\geqslant 2$, we have
  	 	\begin{equation}
  	 	\begin{aligned}
  	 	2(1-2^{1-k})\zeta(k)&=2\sum_{n=1}^{\infty}\frac{(-1)^{n-1}}{n^k}\\&=3(-1)^{m-1}\sigma(\underbrace{2,\ldots,2}_{m})+4\sum_{c=2}^{m}(-1)^{m-c}\sigma(2c,\underbrace{2,\ldots,2}_{m-c}).
  	 	\end{aligned}
  	 	\end{equation}
  	 \end{customcor}
  	 
  	 Indeed, Corollary \ref{ccr2} follows from Corollary \ref{ccr1} by the following lemma:
  	 
  	 \medskip
  	 \begin{lem}
  	 	For all even integers $k\geqslant 2$, we have
  	 	\begin{equation}
  	 	\sum_{b=2}^{k}(-1)^b\zeta(b,\underbrace{1,\ldots,1}_{k-b})=2\sum_{n=1}^{\infty}\frac{(-1)^{n-1}}{n^k}\cdot
  	 	\end{equation}
  	 \end{lem}
  	 We shall deduce this lemma from formula (6) of Y. Ohno and D. Zagier's paper \cite{OZ}, in which we take $y=-x$ and $z=0$. Their formula then becomes the following identity between formal power series in the indeterminate $x$
  	 \begin{equation}\label{EE116}
  	 1+\sum_{k=2}^{\infty}\Big(\sum_{b=2}^{k}(-1)^b\zeta(b,\underbrace{1,\ldots,1}_{k-b}\Big)x^k=\prod_{m=1}^{\infty}\left(1-\frac{x^2}{m^2}\right)^{-1}.
  	 \end{equation}
  	 Now the right hand side is the Taylor expansion at $0$ of $\frac{\pi x}{{\rm sin}(\pi x)}$. But we have
  	 \begin{equation}\label{EE117}
  	 \frac{\pi x}{{\rm sin}(\pi x)}=1+x\sum_{n=1}^{\infty}(-1)^n(\frac{1}{x-n}+\frac{1}{x+n})=1+2x^2\sum_{n=1}^{\infty}\frac{(-1)^{n-1}}{n^2-x^2}
  	 \end{equation}
  	 By comparing, when $k$ is even and $\geqslant 2$, the coefficient of $x^k$ in the left hand side of (\ref{EE116}) and in the right hand side of (\ref{EE117}), we get the lemma.
  	 
  	 \medskip
  	 
 As in the examples e)  and f), we can translate Corollary \ref{CR1} in a statement about generating series. We thus get the following identity, already  deduced from Leshchiner's identities by D. Baily, J. Borwein and D. Bradley (\cite{BBB}, formula (31)).
 
 \medskip
 \begin{customcor}{3}
 			Let $x$ be an indeterminate. In the ring of formal power series ${\bf R}[[x]]$, endowed with the product topology of ${\bf R}^{\bf N}$, we have
 			\begin{equation}\label{EE123}
 			\sum_{k=1}^{\infty}\frac{1}{k^2\left({2k \atop k}\right)}\frac{3k^2+x^2}{k^2-x^2}\prod_{m=1}^{k-1}\left(1-\frac{x^2}{m^2}\right)=2\sum_{n=1}^{\infty}\frac{(-1)^{n-1}}{n^2-x^2}=\frac{\pi}{x {\rm sin}(\pi x)}-\frac{1}{x^2}
 			\end{equation}
 
 \end{customcor}
 \vspace{2em}
 \subsection*{i) Compositions of the form $(2,1,\ldots,1,v)$}
 \myindent Let $\bf a$ be a composition of the form ${\bf a}=(2,\underbrace{1,\ldots,1}_{u-2},v) $, where $u$ and $v$ are integers~$\geqslant2$. Its associated binary word is $0\{1\}_{u-1}\{0\}_{v-1}1$. We have ${\bf a}^{\rm init}=(2,\underbrace{1,\ldots,1}_{u-2}) =\overline{(u)}$, ${\bf a}^{\rm fin}=(v) $, ${\bf a}^{\rm mid}=\varnothing $, hence
 \begin{equation}
 \begin{aligned}
 &\mu(\delta([{\bf a}]))=\delta([u])+\delta([v])+\delta(\varnothing)\\
&\qquad\;\;\; =2\sum_{b=3}^{u}(b,\underbrace{1,\ldots,1}_{u-b})+2\sum_{b=3}^{v}(b,\underbrace{1,\ldots,1}_{v-b})+3(2,\underbrace{1,\ldots,1}_{u-2})+3(2,\underbrace{1,\ldots,1}_{v-2})+\varnothing
 \end{aligned}
 \end{equation}
 by example e), and therefore
 \begin{equation}
 \begin{aligned}
 \delta([{\bf a}])&=2\sum_{b=3}^{u}(b,\underbrace{1,\ldots,1}_{u-b},v)+2\sum_{b=3}^{v}(b,\underbrace{1,\ldots,1}_{v-b},u)\\&\quad\quad+3(2,\underbrace{1,\ldots,1}_{u-2},v)+3(2,\underbrace{1,\ldots,1}_{v-2},u)+(u+v).
 \end{aligned}
 \end{equation}
 \vspace{1em}
 \section{An explicit formula for the map $\delta$}
 \myindent As we saw in the previous section, the map $\delta$ plays a crucial role in the expression of multiple zeta values in terms of multiple Ap{\'e}ry-like sums. It is therefore highly desirable to get an explicit combinatorial expression of $\delta({[{\bf a}]})$ for each admissible composition~$\bf a$. Such an expression will be given in Subsection \ref{s3as2} below.
 
 \medskip
 \subsection{A modified version of the stuffle product}\label{s3s1}

\myindent The stuffle product is a $\bf Z$-bilinear composition law denoted by $\ast$ on ${\bf Z}^{({\mathcal{C}})}$. Let us recall its definition. A {\it stuffling} of two natural numbers $r$ and $s$ in $\bf N$ is a triple $(t,I,J)$ where $t$ is a natural number and $I$, $J$ are sets such that $|I|=r$, $|J|=s$ and $I\cup J=\{1,\ldots,t\}$. We then denote by $\sigma_{I}$ (resp. $\sigma_J$) the unique increasing bijection from $I$ to $\{1,\ldots ,r\}$ (resp. from $J$ to $\{1,\ldots ,s\}$). If ${\bf a}=(a_1,\ldots ,a_r)$ and ${\bf b}=(b_1,\ldots ,b_s)$ are two compositions of depth $r$ and $s$ respectively, we define {\it the composition deduced from} $\bf a$ and $\bf b$ by the stuffling $(t,I,J)$ to be the composition ${\bf c}=(c_1,\ldots ,c_t)$, where
\begin{equation}
c_i=\begin{cases}
{ a}_{\sigma_{_{I}}(i)} &\text{ if  }\; i \in I, i\notin J,\\
{ b}_{\sigma_{_{J}}(i)} &\text{ if  }\; i \notin I, i\in J,
\\
{ a}_{\sigma_{_{I}}(i)}+{ b}_{\sigma_{_{J}}(i)} &\text{ if  }\; i \in I, i\in J.
\end{cases}
\end{equation}

The stuffle product ${\bf a}\ast {\bf b}$ is then defined  as the sum, over all stufflings $(t,I,J)$ of $r$ and $s$, of the compositions deduced from $\bf a$ and $\bf b$ by these stufflings. It is extended  to a $\bf Z$-bilinear composition law
\begin{equation}
\ast:{\bf Z}^{({\mathcal{C}})}\times{\bf Z}^{({\mathcal{C}})} \rightarrow {\bf Z}^{({\mathcal{C}})}
\end{equation}
called {\it the stuffle product}. The stuffle product is compatible with the grading by the weight. It is commutative, associative, has $\varnothing$ as unit element, and ${\bf Z}^{({\mathcal{A}})}$ is stable by this  composition law.

We introduce here a slightly modified version of the stuffle product. It will be defined only for non-empty compositions.  More precisely, it is the $\bf Z$-bilinear map
\begin{equation*}
\boxast:{\bf Z}^{({\mathcal{C^{\ast}}})}\times{\bf Z}^{({\mathcal{C^{\ast}}})} \rightarrow {\bf Z}^{({\mathcal{C^{\ast}}})}
\end{equation*}
defined as follows: if ${\bf a}=(a_1,\ldots ,a_r)$ and ${\bf b}=(b_1,\ldots ,b_s)$ are two non-empty compositions and
\begin{equation}
(a_2,..,a_r)\ast (b_2,\ldots ,b_s)=\sum_{{\bf c}\in \mathcal{C}}\lambda_{{\bf c}} \;{\bf c},
\end{equation} 
then
\begin{equation}
(a_1,\ldots ,a_s)\boxast (b_1,\ldots ,b_s)=\sum_{{\bf c}\in \mathcal{C}}\lambda_{{\bf c}} (a_1+b_1,{\bf c}).
\end{equation}
(In other words, we add the first components, and stuffle the remaining ones).

The composition law $\boxast $ thus defined on ${\bf Z}^{(\mathcal{C}^{\ast})}$ is commutative, associative and maps ${\bf Z}^{(\mathcal{C}^{\ast})} \times {\bf Z}^{(\mathcal{C}^{\ast})}$ to ${\bf Z}^{(\mathcal{A}^{\ast})}$.
\medskip

 {\it Example.-} We have $(3)\ast (4,1)=(3,4,1)+(4,3,1)+(4,1,3)+(7,1)+(4,4)$, hence
 $$(2,3)\boxast (1,4,1)=(3,3,4,1)+(3,4,3,1)+(3,4,1,3)+(3,7,1)+(3,4,4).$$
 
 The main property of the composition law $\boxast $ used in this paper is the following. For $p,q$ integers such that $0\leqslant p <q$ and ${\bf a}=(a_1,\ldots ,a_r)\in \mathcal{C}^*$, let $\phi_{p,q}({\bf a})$ denote the rational number defined by
 \begin{equation}
 \varphi_{p,q}({\bf a})=q^{-a_1}\sum_{q>n_2\ldots >n_r>p}n_2^{-a_2}\ldots n_r^{-a_r}.
 \end{equation}
 Let us extend this map $\varphi_{p,q}$ to ${\bf Z}^{^{(\mathcal{C}^{\ast})}}$ by $\bf Z$-linearity. Then we have, for ${\bf a}, {\bf b} \in \mathcal{C}^\ast $,
 \begin{equation}\label{EE61}
\varphi_{p,q}({\bf a}\boxast {\bf b})=\varphi_{p,q}({\bf a})\;\varphi_{p,q}({\bf b}).
 \end{equation}
 
\vspace{2em}
 \subsection{An explicit expression of $\delta([{\bf a}])$}\label{s3as2}
 \myindent Let $\bf a$ be a non empty admissible composition of weight $k$. Let $\varepsilon_1\ldots \varepsilon_k$ denote its associated binary word ${\bf w}({\bf a})$.  For $1 \leqslant i \leqslant k-1$, we denote by ${\bf a}_i$  the composition such that ${\bf w}({\bf a}_i)=\varepsilon_{i+1} \ldots  \varepsilon_k$ and hence by $\overline{\bf a}_{k-i}$  the composition such that  ${\bf w}(\overline{\bf a}_{k-i})=\overline{\varepsilon_{i}} \ldots  \overline{\varepsilon_{1}}.$   We note that even when the compositions ${\bf a}_i$ and $\overline{\bf a}_{k-i}$ are not  admissible, $\overline{\bf a}_{k-i}\boxast {\bf a}_i$ is an admissible composition.
 
 \medskip
  \begin{thm}\label{Th10} 	With the previous notations, we have
 \begin{equation}
 \delta([{\bf a}])=\sum_{i=1}^{k-1}(1+\overline{\varepsilon_i}+\varepsilon_{i+1})\;\overline{\bf a}_{k-i}\boxast {\bf a}_i.
 \end{equation}
  \end{thm}
  
  \medskip
  {\it Remark.}- Note that $\varepsilon_{i+1}$ is the first bit of the binary word ${\bf w}({\bf a}_i)$ and $\overline{\varepsilon_i}$ is the first bit of ${\bf w}(\overline{\bf a}_{k-i})$. Hence $1+\overline{\varepsilon_i}+\varepsilon_{i+1}$ is equal to $1$ when both ${\bf a}_i$ and $\overline{\bf a}_{k-i}$ are admissible, to $2$ when only one of them is admissible, to $3$ when none of them is admissible.
  \medskip
  
  \noindent{\it First proof, by reference to one of our previous papers \cite{akhi}}
  \\
   We deduce as a special case of Theorem 9 of \cite{akhi} that, for each integer $n \geqslant 0$, we have,
   \begin{equation}\label{E15}
   \zeta({\bf a})_{n,n}=\sum_{m>n}\left({2m \atop m}\right)^{-1}\sum_{i=1}^{k-1}(1+\overline{\varepsilon_i}+\varepsilon_{i+1})\varphi_{n,m}(\overline{\bf a}_{k-i}) \varphi_{n,m}({\bf a}_i),
   \end{equation}
   and hence by (\ref{EE61})
   \begin{equation}
   \begin{aligned}
      \zeta({\bf a})_{n,n}&=\sum_{m>n}\left({2m \atop m}\right)^{-1}\sum_{i=1}^{k-1}(1+\overline{\varepsilon_i}+\varepsilon_{i+1})\varphi_{n,m}(\overline{\bf a}_{k-i}\boxast {\bf a}_i) \\ \nonumber &= \sum_{i=1}^{k-1} (1+\overline{\varepsilon_i}+\varepsilon_{i+1})\;\sigma(\overline{\bf a}_{k-i}\boxast {\bf a}_i)_n.
   \end{aligned}
   \end{equation}
Theorem \ref{Th10} then follows from the unicity result stated in the remark of Section \ref{s3s2}.

\medskip
\noindent{\it Second proof, by a direct combinatorial argument.}

Let us define a $\bf Z$-linear map $d: {\bf Z}^{(\mathcal{A})}\rightarrow{\bf Z}^{(\mathcal{A})}$ by the formula
\begin{equation}
d({\bf a})=\begin{cases}
\varnothing & \text{   when   } {\bf a}=\varnothing,\\
\sum_{i=1}^{k-1}(1+\overline{\varepsilon_i}+\varepsilon_{i+1})\;\overline{\bf a}_{k-i}\boxast {\bf a}_i & \text{  when  } {\bf a}\in \mathcal{A}^*,\\
\end{cases}
\end{equation}
where in the second case the notations are those of the theorem. This map is graded of degree $0$ for the weight.

\medskip
\begin{lem}\label{l0}
We have $d({\overline{\bf a}})=d({\bf a})$ for each ${\bf a}\in \mathcal{A}$.
\end{lem}

 This is clear when ${\bf a}=\varnothing$. When ${\bf a}\neq \varnothing$ , it follows from the commutativity of the composition law $\boxast$ and the remark above.

\medskip

Let us extend  the definition of the $\bf Z$-bilinear composition law $\boxast$ to $ {\bf Z}^{(\mathcal{C})}\times {\bf Z}^{(\mathcal{C})}$ by the following conventions : ${\bf a}\boxast \varnothing=\varnothing \boxast {\bf a}=0$ when ${\bf a} \in \mathcal{A}^*$ and $\varnothing \boxast \varnothing=\varnothing$.

\medskip
\begin{lem}\label{l1}
	Let $\bf a$ , $\bf b$ be two non empty compositions. We  have
	\begin{equation}\label{EN113}
	\mu({\bf a}\boxast {\bf b})={\bf a}^{\rm init}\boxast {\bf b}+{\bf a}\boxast {\bf b}^{\rm init}+{\bf a}^{\rm init}\boxast {\bf b}^{\rm init}
	\end{equation}
\end{lem}

We first consider the case when $\bf a$ has depth $1$. We  then have 
  ${\bf a}=(a_1)$, ${\bf a}^{\rm init}=\varnothing$, \mbox{ ${\bf b}=(b_1,\ldots ,b_s)$} with $s\geqslant1$, ${\bf a}\boxast {\bf b}=(a_1+b_1,b_2,\ldots ,b_s)$. When $s\geqslant 2$, we have  $\mu({\bf a}\boxast {\bf b})=(a_1+b_1,b_2,\ldots,b_{s-1})={\bf a}\boxast {\bf b}^{\rm init}$  and ${\bf a}^{\rm init} \boxast {\bf b}=0$, ${\bf a}^{\rm init} \boxast {\bf b}^{\rm init}=0$. When $s=1$, we have $\mu({\bf a}\boxast {\bf b})=\varnothing$ and ${\bf a}^{\rm init} \boxast {\bf b}=0$, ${\bf a} \boxast {\bf b}^{\rm init}=0$, ${\bf a}^{\rm init} \boxast {\bf b}^{\rm init}=\varnothing$, hence (\ref{EN113}) holds. The lemma similarly holds when the depth of ${\bf b}$ is $1$, since $\boxast$ is commutative.

We now assume that ${\bf a}$ and ${\bf b}$ have depth $\geqslant 2$. In this case, we can write ${\bf a}=({\bf a}^{\rm init},a)$ and ${\bf b}=({\bf b}^{\rm init},b)$ with ${\bf a}^{\rm init},{\bf b}^{\rm init}$  non empty compositions, and we have, by definition of the composition law $\boxast$,
	\begin{equation}
	{\bf a}\boxast {\bf b}=({\bf a}^{\rm init}\boxast {\bf b},a)+({\bf a}\boxast {\bf b}^{\rm init},b)+({\bf a}^{\rm init}\boxast {\bf b}^{\rm init},a+b).
	\end{equation}
The lemma follows.

\medskip
\begin{lem}\label{l2}
	Let $\bf a$ be an admissible composition. With the notations of the beginning of the section, we have
	\begin{subequations}
	\begin{align}\label{E115a}
	\; d({\bf a}^{\rm init}) &=\sum_{i=1}^{k-1}(1+\overline{\varepsilon_i}+\varepsilon_{i+1})\;\overline{\bf a}_{k-i}\boxast ({\bf a}_i)^{\rm init}\qquad \qquad\text{ if }\quad {\bf a}^{\rm init}\neq \varnothing,
	\\ \label{E115b} d({\bf a}^{\rm fin}) &=\sum_{i=1}^{k-1}(1+\overline{\varepsilon_i}+\varepsilon_{i+1})\;(\overline{\bf a}_{k-i})^{\rm init}\boxast {\bf a}_i\qquad \qquad\text{ if }\quad {\bf a}^{\rm fin}\neq \varnothing,
	\\ \label{E115c} d({\bf a}^{\rm mid}) &=\sum_{i=1}^{k-1}(1+\overline{\varepsilon_i}+\varepsilon_{i+1})\;(\overline{\bf a}_{k-i})^{\rm init}\boxast ({\bf a}_i)^{\rm init}\qquad \text{ if }\quad {\bf a}^{\rm init}\neq \varnothing \text{ and } {\bf a}^{\rm fin}\neq \varnothing.
	\end{align}
	\end{subequations}

\end{lem}

Assume first that ${\bf a}^{\rm init} \neq \varnothing$. We then have ${\bf a}=(a_1,\ldots ,a_r)$ with $r \geqslant 2$,  ${\bf a}^{\rm init}=(a_1,\ldots ,a_{r-1})$ and therefore
\begin{equation}
 d({\bf a}^{\rm init}) =\sum_{i=1}^{k'-1}(1+\overline{\varepsilon_i}+\varepsilon_{i+1})\;(\overline{{\bf a}^{\rm init }})_{k'-i}\boxast ({\bf a}^{\rm init})_i\qquad
\end{equation}
where $k'=k-a_r$. Now, for $1 \leqslant i\leqslant k'-1$, $(\overline{{\bf a}^{\rm init }})_{k'-i}$ is equal to $\overline{\bf a}_{k-i}$ and  $({\bf a}^{\rm init})_i$ is equal to $({\bf a}_i)^{\rm init}$. On the other hand, for $k'\leqslant i \leqslant k-1$, we have $({\bf a}_i)^{\rm init}=\varnothing$ and $\overline{\bf a}_{k-i} \neq \varnothing$, hence $\overline{\bf a}_{k-i}\boxast {\bf a}_i^{\rm init}=0.$  Formula (\ref{E115a}) follows. 

Formula (\ref{E115b}) is proved similarly (or deduced from the previous case  by duality)  when ${\bf a}^{\rm fin}\neq\varnothing$.

Assume from now on that ${\bf a}^{\rm init}\neq\varnothing$ and ${\bf a}^{\rm fin}\neq\varnothing$ and let $k'$ be $k-a_r$ as above. We distinguish two cases:
\\a) When ${\bf a}^{\rm mid} \neq \varnothing$, we have
\begin{equation}
d({\bf a}^{\rm mid})=d(({\bf a}^{\rm init})^{\rm fin})=\sum_{i=1}^{k'-1}(1+\overline{\varepsilon_i}+\varepsilon_{i+1})\Big(\big(\overline{{\bf a}^{\rm init}}\big)_{k'-i}\Big)^{\rm init}\boxast ({\bf a}^{\rm init})_i,
\end{equation}
by formula (\ref{E115b}) applied to ${\bf a}^{\rm init}$, hence, by the same arguments as before,
\begin{equation}
d({\bf a}^{\rm mid})=\sum_{i=1}^{k'-1}(1+\overline{\varepsilon_i}+\varepsilon_{i+1})(\overline{\bf a}_{k-i})^{\rm init}\boxast ({\bf a}_i)^{\rm init}=\sum_{i=1}^{k-1}(1+\overline{\varepsilon_i}+\varepsilon_{i+1})(\overline{\bf a}_{k-i})^{\rm init}\boxast ({\bf a}_i)^{\rm init},
\end{equation}
since for $k' \leqslant i \leqslant k-1$, we have $({\bf a}_i)^{\rm init}=\varnothing$ and $(\overline{\bf a}_{k-i})^{\rm init} \neq \varnothing$. Hence (\ref{E115c}) holds. 
\\b) When ${\bf a}^{\rm mid}=\varnothing$, the binary word ${\bf w}({\bf a})$ associated to $\bf a$ has the form $\underbrace{01\ldots1}_{k'}\underbrace{0\ldots01}_{k-k'}$, with $2\leqslant k'\leqslant k-2$. We then have $({\bf a}_i)^{\rm init}\neq \varnothing$ for $1 \leqslant i \leqslant k'-1$, $({\bf a}_i)^{\rm init}=\varnothing$ for \mbox{$k' \leqslant i \leqslant k-1$},
$(\overline{{\bf a}}_{k-i})^{\rm init}=\varnothing$ for $1 \leqslant i \leqslant k'$, $(\overline{{\bf a}}_{k-i})^{\rm init}\neq \varnothing$ for $k'+1 \leqslant i \leqslant k-1$, hence $(\overline{{\bf a}}_{k-i})^{\rm init}\boxast ({\bf a}_i)^{\rm init}=0$ for $1\leqslant i \leqslant k-1$, $i\neq k'$. For $i=k'$ we have $\varepsilon_i=1$, $\varepsilon_{i+1}=0$, $1+\overline{\varepsilon_i}+\varepsilon_{i+1}=1$ and $(\overline{{\bf a}}_{k-i})^{\rm init}\boxast ({\bf a}_i)^{\rm init}=\varnothing$, hence (\ref{E115c}) holds.

\medskip

Theorem \ref{Th10} is now a consequence of the following:

\smallskip
\begin{lem}
We have $\delta([{\bf a}])=d({\bf a})$ for all ${\bf a}\in \mathcal{A}$.
\end{lem}

We shall prove the lemma by induction on the weight $k$ of $\bf a$. When ${\bf a}=\varnothing$, we have $\delta([{\bf a}])=\varnothing$ and $d({\bf a})= \varnothing$. When $\bf a$ has depth $1$, we have ${\bf a}=(k)$ with $k\geqslant2$. For $1\leqslant i \leqslant k-1$, we have $\overline{\bf a}_{k-i}=(\underbrace{1,\ldots,1}_{i})$, ${\bf a}_i=(k-i)$, $\overline{\bf a}_{k-i}\boxast {\bf a}_i=(k-i+1,\underbrace{1,\ldots,1}_{i-1})$, $\varepsilon_i=0$, and $\varepsilon_{i+1}$ is equal to $0$ if $1\leqslant i\leqslant k-2$ and to $1$ if $i=k-1$. Hence
\begin{equation}
d({\bf a})=2\sum_{i=1}^{k-2}(k-i+1,\underbrace{1,\ldots,1}_{i-1})+3(2,\underbrace{1,\ldots,1}_{k-2})=\delta([{\bf a}])
\end{equation}
by Theorem \ref{Th8}. When $\overline{\bf a}$ has depth $1$, we have $d({\bf a})=d(\overline{\bf a})=\delta([\overline{\bf a}])=\delta([{\bf a}])$ by lemma~\ref{l0}.

We now assume that $\bf a$ and $\overline{\bf a}$ have depth $\geqslant 2$, hence ${\bf a}^{\rm init}\neq \varnothing$ and ${\bf a}^{\rm fin}\neq\varnothing$. We have
\begin{equation}
\begin{aligned}
&\mu(d({\bf a}))=\sum_{i=1}^{k-1}(1+\overline{\varepsilon_i}+\varepsilon_{i+1})\;\mu(\;\overline{\bf a}_{k-i}\boxast {\bf a}_i)\\
&\qquad\quad=\sum_{i=1}^{k-1}(1+\overline{\varepsilon_i}+\varepsilon_{i+1})\Big({(\overline{\bf a}_{k-i})}^{\rm init}\boxast {{\bf a}_i}+{\overline{\bf a}_{k-i}}\boxast {({\bf a}_i)}^{\rm init}+{(\overline{\bf a}_{k-i})}^{\rm init}\boxast {({\bf a}_i)}^{\rm init}\Big)
\end{aligned}
\end{equation}
by Lemma~\ref{l1} applied to the non empty compositions $\overline{\bf a}_{k-i}$ and ${\bf a}_i$, hence
\begin{equation}
\mu(d({\bf a}))=d({\bf a}^{\rm init})+d({\bf a}^{\rm mid})+d({\bf a}^{\rm fin})
\end{equation} 
by Lemma \ref{l2}, and therefore $\mu(d({\bf a}))=\delta([{\bf a}^{\rm init}])+\delta([{\bf a}^{\rm mid}])+\delta([{\bf a}^{\rm fin}])=\delta(\alpha([{\bf a}]))=\mu(\delta([{\bf a}]))$ by the induction hypothesis. The equality  $d({\bf a})=\delta([{\bf a}])$ follows, since $d({\bf a})$ and $\delta([{\bf a}])$ both belong to ${\bf Z}^{\mathcal{A}_k}$, and $\mu$ is injective on  ${\bf Z}^{\mathcal{A}_k}$.
\vspace{ 2em}

\subsection{On the support of $\delta([{\bf a}])$}\label{s4s3}
\medskip
An admissible binary word $w$ can always be uniquely written as
\begin{equation}
\{0\}_{u_1}\{1\}_{v_1}\ldots \{0\}_{u_h}\{1\}_{v_h}=\underbrace{0\ldots 0}_{u_1}\underbrace{1\ldots1}_{v_1}\ldots\underbrace{0\ldots 0}_{u_h}\underbrace{1\ldots1}_{v_h}
\end{equation}
where $h\geqslant 0$ and $u_i\geqslant1, v_i \geqslant 1$ for $1 \leqslant i \leqslant h$. The integer $h$ is called  the {\it height}
of the word $w$. We define the {\it height of an admissible composition} $\bf a$ as the height of its associated binary word ${\bf w}({\bf a})$. It is the number of entries of ${\bf a}$ larger than or equal to $2$.

\medskip

\begin{prop}\label{TH13}
Let $\bf a$ be an admissible composition of weight $k$ and height $h$. Any admissible composition in the support of $\delta([{\bf a}])$ has at most $k-2h$ odd entries.	
\end{prop}
 Let $\varepsilon_1\ldots\varepsilon_k$ be the binary word ${\bf w}({\bf a})$ associated to $\bf a$. Note that $k-2h$ is the number of indices $i$ such that $1 \leqslant i \leqslant k-1$ and $\varepsilon_i=\varepsilon_{i+1}$.
 
 Fix an index $i$ such that $1\leqslant i \leqslant k-1$. As in Theorem \ref{Th10}, let ${\bf a}_i$ and $\overline{\bf a}_{k-i}$ be the compositions  with  associated binary words $\varepsilon_{i+1}\ldots\varepsilon_{k}$ and $\overline{\varepsilon_i}\ldots\overline{\varepsilon_1}$ respectively. Let $p_i$ (resp. $q_i$) denote the number of indices $j$ such that $i+1 \leqslant j \leqslant k-1$ (resp. $1\leqslant j\leqslant i-1$) and $\varepsilon_j=\varepsilon_{j+1}$.  From the previous alinea, we get
 \begin{equation}\label{EU149}
 k-2h=\begin{cases}
 p_i+q_i & \text{ if } \varepsilon_i\neq \varepsilon_{i+1},\\
 p_i+q_i+1   & \text{ if } \varepsilon_i= \varepsilon_{i+1}.
 \end{cases}
 \end{equation}
 Let $p_i'$ (resp. $q_i'$) denote the number of odd entries of ${\bf a}_i$ (resp. $\overline{\bf a}_{k-i}$). We have
 \begin{equation}\label{EU150}
  \begin{array}{c c}
     p_i'\leqslant\begin{cases}
     p_i & \text{if } \varepsilon_{i+1}=0,\\
     p_i+1 & \text{if } \varepsilon_{i+1}=1,
     \end{cases}&\qquad{q_i'\leqslant\begin{cases}
   q_i & \text{if } \varepsilon_{i}=1,\\
   q_i+1 & \text{if } \varepsilon_{i}=0.
   \end{cases}}
    \end{array} 
 \end{equation}
 Any admissible composition {\bf b} in the support of $\overline{{\bf a}}_{k-i}\boxast {\bf a}_i$ has at most $p_i'+q_i'$ odd entries; it has even at most $p_i'+q_i'-2$ odd entries when  $(\varepsilon_i,\varepsilon_{i+1})=(0,1)$, since in this  case both ${\bf a}_i $ and $\overline{\bf a}_{k-i}$ have $1$ as  first entry. Hence ${\bf b}$ has at most $k-2h$ entries, as is seen  by using relations (\ref{EU149}) and (\ref{EU150}).
 
 The multiple zeta values $\zeta({\bf a})$, where $\bf a$ is an admissible composition whose entries belong to $\{2,3\}$ are particularly interesting to study since, by a theorem of F. Brown~\cite{FB}, any other multiple zeta value is a $\bf Z$-linear combination of them and conjecturally, they are $\bf Z$-linearly independent. For these multiple zeta values, we have:
\medskip

\begin{customcor}{1}\label{pc1}
Let ${\bf a}$ be an admissible composition whose entries belong to $\{2,3\}$, and let~$s$ denote the number of  entries of $\bf a$ equal to $3$. Let $\bf b$ be an admissible composition in the support of $\delta([{\bf a}])$. The entries of $\bf b$ are at most $5$, and at most $s$ of them are odd.
\end{customcor} 

The first assertion follows from Theorem \ref{Th10}, and the second one from  Proposition~\ref{TH13}, since $k-2h$ is equal to $s$ for such a composition $\bf a$.
\medskip

Numerical experiment suggests that, for $\bf a$ as in Corollary \ref{pc1}, $\zeta({\bf a})$ can be expressed as a $\bf Q$-linear combination of a much smaller number of  Ap{\' e}ry-like sums than those occurring in $\sigma(\delta([{\bf a}]))$. We formulate in this direction the following conjecture:
\medskip
\begin{coj}\label{Cu1}
	Let $\bf a$ be an admissible composition  whose entries belong to $\{2,3\}$, and let $s$ denote the number of entries of $\bf a$ equal to $3$. Then $\zeta({\bf a})$ is a $\bf Q$-linear combination of Ap{\'e}ry-like sums of the form $\sigma({\bf b})$, where $\bf b$ runs over the admissible compositions of the same weight as $\bf a$, with entries in $\{1,2,3\}$ and at most $s$ of them odd.
\end{coj} 
\medskip

A partial result in this direction is:
\medskip

\begin{thm}
	Conjecture \ref{Cu1} holds for $s=0$ and $s=1$.
\end{thm}

For $s=0$, Conjecture \ref{Cu1} says that $\zeta(\underbrace{2,\ldots,2}_r)$ is a rational multiple of $\sigma(\underbrace{2,\ldots,2}_r)$. This follows from the fact that the first one is equal to $\frac{\pi^{2r}}{(2r+1)!}$ and the second one to $\frac{\pi^{2r}}{3^{2r}(2r)!}$ (see Section \ref{s2s7a}, formula (\ref{EU71a})).

By Zagier's Theorem (see Section \ref{s2s8a}, Theorem \ref{Th7}), $\zeta(\underbrace{2,\ldots,2}_a,3,\underbrace{2,\ldots,2}_b) $ belongs to the $\bf Q$-subspace $ W$ of $\bf R$ generated by the real numbers $\pi^{k-1-2r}\zeta(2r+1)$, where \mbox{$k=2a+2b+3$} and $1\leqslant r \leqslant\frac{k-1}{2}$. By our Theorem \ref{th7}, $ W$ is contained in the $\bf Q$-subspace $ V$ of $\bf R$ generated by the Ap{\'e}ry-like sums $\sigma({\bf a})$, where $\bf a$ is an admissible composition  of weight $k$, either of the form $(2,\ldots,2,3,2,\ldots,2)$ or of the form $(2,\ldots,2,1,2,\ldots,2)$. This implies Conjecture~\ref{Cu1} when $s=1$.

\medskip

{\it Examples.} 1) One deduces from formulas (\ref{EU127}) and (\ref{EU128}) that
\begin{equation}
\begin{aligned}
\zeta(3,2)&=\frac{1}{81}\left(140\sigma(3,2)+276\sigma(2,3)+459\sigma(2,2,1)+210\sigma(2,1,2)\right),\\
\zeta(2,3)&=\frac{1}{81}\left(610\sigma(3,2)+786\sigma(2,3)+594\sigma(2,2,1)+915\sigma(2,1,2)\right). 
\end{aligned}
\end{equation} 

2) One proves in a similar way that
\begin{equation*}
\begin{aligned}
6501255\zeta(3,2,2) =&-38573600\sigma(3,2,2)-24271152\sigma(2,3,2) +126621792\sigma(2,2,3)\\&+169412715\sigma(2,2,2,1) +48270762\sigma(2,2,1,2)- 57860400\sigma(2,1,2,2),\\
1300251\zeta(2,3,2)=& \;25510840\sigma(3,2,2)+ 33723696\sigma(2,3,2)+ 23152392\sigma(2,2,3)\\&+ 24948000\sigma(2,2,2,1)+ 27718389\sigma(2,2,1,2)+ 38266260\sigma(2,1,2,2),\\
2167085\zeta(2,2,3)=&\; 103897690\sigma(3,2,2)+ 138573066\sigma(2,3,2)+ 95423994\sigma(2,2,3)\\&+ 78682590\sigma(2,2,2,1)+ 117029394\sigma(2,2,1,2)+ 155846535\sigma(2,1,2,2).
\end{aligned}
\end{equation*}

For $s\geqslant 2$, we have no theoretical result at present, but only numerical evidence, based on PSLQ algorithm, that conjecture \ref{Cu1} should at least hold in weight $\leqslant 12$. In fact these experiments suggest that an even stronger pattern appears when $s\geqslant 2$. More precisely:

a) When $s=2$ (which implies that the weight $k$ is even and $\geqslant 6$), it seems that, in conjecture \ref{Cu1},  compositions $\bf b$  ending by $1$ can be omitted, and that the remaining Ap{\'e}ry-like sums $\sigma({\bf b})$ are then $\bf Z$-linearly independent.

b) When $s=3$  (which implies that the weight $k$ is odd and $\geqslant 9$) it seems that, in conjecture \ref{Cu1},  compositions $\bf b$  ending by $1$ or by $(1,2)$ can be omitted, but the remaining  Ap{\'e}ry-like sums $\sigma({\bf b})$ do not seem to be  $\bf Z$-linearly independent: there seems to exist one linear relation among them in weight $9$, and four independent ones in weight $11$.

c) Note that in the examples 1) and 2), the coefficients of two Ap{\'e}ry-like sums of the form $\sigma(3,{\bf c})$ and $ \sigma(2,1,{\bf c})$  are always proportional to $2$ and $3$. This can be proved for $s=1$ by using the results of \ref{s2s8} (under the assumption that the conjecture stated in the remark of section \ref{s2s9} holds). The numerical evidence we have gathered suggests that this might continue to hold for $s=2$ and $s=3$.

 \medskip
{\it Remark.}- It is worth here mentioning the following result which, according to a private communication of Francis Brown, follows from his theory of the unipotency degree, also sometimes referred to as the {\it motivic depth} (see \cite{Brown}, p. 29): any product of an even non-negative power of $\pi$ by a multiple zeta value of depth $\leqslant r$ is a $\bf Q$-linear combination of multiple zeta values of the form $\zeta({\bf a})$, where ${\bf a}$ runs over the compositions whose entries belongs to $\{2,3\}$, with at most $r$ of them equal to $3$. This is an indication that to get a proof of  Conjecture \ref{Cu1}, one perhaps needs to develope a motivic theory of Ap{\'e}ry-like sums and of their relations with motivic multiple zeta values.
\vspace{2em}
\subsection{Some further examples}\label{s4as4}
\myindent All formulas obtained in  Lemma \ref{L} and  in Theorems \ref{Th8}, \ref{Th9}, \ref{TH10} of the examples of Subsection \ref{S3S3} can of course be proved by a direct application of  Theorem \ref{Th10}. We add in this subsection some further examples.

\vspace{2em}
\subsection*{a) Compositions of height 1}
\myindent 
An admissible composition $\bf a$ of height $1$ (as defined in Section \ref{s3as2}) has the form $(u+1,\underbrace{1,\ldots,1}_{v-1})$ with $u\geqslant1$, $v\geqslant1$, its associated binary word being $\{0\}_u\{1\}_v$. Its weight is $u+v$ and its depth~is~$v$. Note that $\overline{\bf a}$ is then the composition $(v+1,\underbrace{1,\ldots,1}_{u-1})$, whose associated binary word is~$\{0\}_v\{1\}_u$.

\medskip
\begin{prop}
	Let $u,v$ be integers such that $1 \leqslant u \leqslant v$, and $\bf a$ be the composition \mbox{$(u+1,\underbrace{1,\ldots,1}_{v-1})$}. We have
	\begin{equation}
	\delta([{\bf a}])=\sum_{\bf b}c_{\bf b}{\bf b},
	\end{equation}
	where ${\bf b}=(b_1,\ldots,b_r)$ runs over all compositions of weight $u+v$ with the following properties: we have $2\leqslant b_1\leqslant v+1$, $b_i \in \{1,2\}$ for $2 \leqslant i \leqslant r$, and $r\geqslant{\rm max}(u,v+2-b_1).$ Furthermore, for such a composition ${\bf b}$, we have
	\begin{equation}\label{EU155}
	c_{\bf b}=\begin{cases}
	2\Big({s \atop r-u}\Big)& \text{ if $b_1 \geqslant 3$ and $r<v$},\\
	2\Big({s \atop r-u}\Big)+2\Big({s \atop r-v}\Big)	& \text{ if $b_1 \geqslant 3$ and $r\geqslant v$},\\
	3\Big({s \atop r-u}\Big)& \text{ if $b_1 =2$},
	\end{cases} 
	\end{equation}
	where $s=2r+b_1-u-v-2$ is the number of entries of $\bf b$ equal to $1$.
\end{prop}
The binary word ${\bf w}({\bf a})$ is $\{0\}_u\{1\}_v$ and its weight is $k=u+v$. We have therefore, with the notations of Theorem \ref{Th10},
\begin{equation}
{\bf a}_i=\begin{cases}
(u-i+1,\underbrace{1,\ldots,1}_{v-1}) & \text{ if $1 \leqslant i \leqslant u$,}\\
(\underbrace{1,\ldots,1}_{u+v-i}) & \;\text {if $\;u \leqslant i \leqslant u+v-1$,}
\end{cases}
\end{equation}

\begin{equation}
\overline{\bf a}_{k-i}=\begin{cases}
(\underbrace{1,\ldots,1}_{i})& \text{ if $1 \leqslant i \leqslant u$,}\\
 (i-u+1,\underbrace{1,\ldots,1}_{u-1}) & \;\text {if $\;u \leqslant i \leqslant u+v-1$,}
 	\end{cases} 
\end{equation}
hence

\begin{equation}
\overline{\bf a}_{k-i}\boxast {\bf a}_i =\begin{cases}
\big(u-i+2,(\underbrace{1,\ldots,1}_{i-1})\ast(\underbrace{1,\ldots,1}_{v-1} )\big)& \text{ if $1 \leqslant i \leqslant u$,}\\
\big(i-u+2,(\underbrace{1,\ldots,1}_{u-1})\ast(\underbrace{1,\ldots,1}_{u+v-i-1} ) \big)& \;\text {if $\;u \leqslant i \leqslant u+v-1$.}
\end{cases} 
\end{equation}
Applying Theorem \ref{Th10}, we get 
\begin{equation}\label{EU159}
\begin{aligned}
\delta([{\bf a}])=&2\sum_{i=1}^{u-1}\big(u-i+2,(\underbrace{1,\ldots,1}_{i-1})\ast(\underbrace{1,\ldots,1}_{v-1} )\big)+3\big(2,(\underbrace{1,\ldots,1}_{u-1})\ast(\underbrace{1,\ldots,1}_{v-1} )\big)\\
&+2\sum_{i=u+1}^{u+v-1}\big(i-u+2,(\underbrace{1,\ldots,1}_{u-1})\ast(\underbrace{1,\ldots,1}_{u+v-i-1}) \big)\\
=&2\sum_{b_1=3}^{u+1}\big(b_1,(\underbrace{1,\ldots,1}_{u+1-b_1})\ast(\underbrace{1,\ldots,1}_{v-1} )\big)+3\big(2,(\underbrace{1,\ldots,1}_{u-1})\ast(\underbrace{1,\ldots,1}_{v-1} )\big)\\
&2\sum_{b_1=3}^{v+1}\big(b_1,(\underbrace{1,\ldots,1}_{u-1})\ast(\underbrace{1,\ldots,1}_{v+1-b_1} ) \big).
\end{aligned}
\end{equation}
We see that all compositions occurring in this expression are necessarily of the form ${\bf b}=(b_1,\ldots,b_r)$, with the weight of $\bf b$ equal to $u+v$, $2\leqslant b_1\leqslant {\rm max}(u+1,v+1)=v+1$ and $b_i \in \{1,2\}$ for $2 \leqslant i\leqslant r$. 
\medskip

To continue the proof, we need the following lemma:

\medskip
\begin{lem} \label{L8}Let $m$ and $n$ be natural numbers.
	The stuffle product $(\underbrace{1,\ldots,1}_{m})\ast(\underbrace{1,\ldots,1}_{n} )$ is the sum $\sum_{\bf c}\lambda_{\bf c}{\bf c}$ where $\bf c$  runs over all  compositions of weight $m+n$ and depth $r\geqslant {\rm max}(m,n)$, whose entries are all equal to $1$ or $2$, the coefficients $\lambda_{\bf c}$ being $\left({s \atop r-m}\right)=\left({s \atop r-n}\right)$ where $s=2r-n-m$ is the number of entries of the composition ${\bf c}$ equal to $1$.
\end{lem}

A stuffling of $m$ and $n$ is a triple $(r,I,J)$ where $r$ is a natural number and $I,J$ sets such that $|I|=m$, $|J|=n$ and $I\cup J=\{1,\ldots,r\}$. Clearly this imposes that $r\geqslant {\rm max}(m,n).$ The composition deduced from $(\underbrace{1,\ldots,1}_{m})$ and $(\underbrace{1,\ldots,1}_{n})$ by this stuffling is the composition $\bf c$ whose $i$ th entry is $2$ if $i\in I \cap J$ and $1$ otherwise. Its weight is $m+n$. We have $|I-(I\cap J)|=r-n$ and $|J-(I\cap J)|=r-m$, hence the number $s$ of entries equal to $1$ in $\bf c$ is $2r-n-m$. Finally, when $\bf c$ is given, the number of stufflings $(r,I,J)$ leading to $\bf c$ is the number of subsets of cardinality $r-n$ in a set of $s$ elements. The lemma follows.
\medskip

Now let us consider a composition ${\bf b}=(b_1,\ldots,b_r)$ of weight $u+v$, with $2 \leqslant b_1 \leqslant v+1$ and $b_i \in \{1,2\}$ for $2 \leqslant i \leqslant r$. We want to compute the coefficient $c_{\bf b}$ with which it occurs in $\delta([{\bf a}])$.

By Lemma \ref{L8}, this composition $\bf b$ occurs in the first term of (\ref{EU159}) if and only if $3 \leqslant b_1 \leqslant u+1$ and $r\geqslant 1+{\rm max}(u+1-b_1,v-1)={\rm max}(u+2-b_1,v)=v$,  and its coefficient in this term is then $2\left({s \atop (r-1)-(v-1)}\right)=2\left(s \atop r-v\right)$, where $s$ is the number of entries of $\bf b$ equal to~$1$.

Similarly it occurs in the second term of (\ref{EU159}) if and only if $b_1=2$ and \mbox{$r\geqslant {\rm max}(u,v)=v$,} and its coefficient in this term is then $3\left({s \atop r-v}\right)=3\left({s \atop r-u}\right)$.

It occurs in the third term of (\ref{EU159}) if and only if  we have $3 \leqslant b_1 \leqslant v+1$ and $r\geqslant {\rm max}(u,v+2-b_1)$, and its coefficient in this term is then $2\left({s \atop r-u}\right).$

In all these cases, we have ${\rm max}(u,v+2-b_1)\leqslant v,$ and hence $r \geqslant {\rm max}(u,v+2-b_1).$ When this condition is satisfied, we see that the total of the previous contributions agree with the formula (\ref{EU155}), since $r\geqslant v$ implies $u+v=k\geqslant b_1+v-1$, hence $b_1 \leqslant u+1$, and since  $r \geqslant {\rm max}(u,v+2-b_1)$ implies $r \geqslant v$ when $b_1=2$.
\vspace{2em}
\subsection*{b) Compositions of the form $(a,\ldots,a)$}

\myindent In Theorem \ref{TH10} of Section \ref{S3S3}, example g), we have computed $\smash{\delta([\underbrace{2,\ldots,2}_{m}])}$. We now do the same for $\delta([\underbrace{a,\ldots,a}_{m}])$, when $a\geqslant 3$.

\medskip
\begin{thm}\label{Th13}
	For all integer $a\geqslant 3$ and $m \geqslant 1$ , we have
	\begin{equation}
	\delta(\underbrace{a,\ldots,a}_{m})=\sum c_{\bf b}{\bf b}
	\end{equation}
	where $\bf b$ runs over the compositions $(b_1,\ldots,b_r)$ with the following properties:
	\\ {\rm(i)} either $2 \leqslant b_1\leqslant a$ or $b_1=a+2$;
	\\ {\rm(ii)} for $2\leqslant i \leqslant r$, we have  $b_i \in \{1,2,a,a+1,a+2\};$
	\\ {\rm(iii)} the number $s$ of of indices $i \in \{2,\ldots,r\}$ such that $b_i\geqslant a$ is at most $ m-1$ when $2 \leqslant b_1 \leqslant a$ and at most $ m-2$ when $b_1=a+2$;
	\\ {\rm(iv)} the composition obtained by removing from $\bf b$ the first entry and those equal to $a$, and by then replacing those equal to $a+1$ by $1$ and those equal to $a+2$ by $2$, is
	\begin{equation}
	\underbrace{1,\ldots,1}_{u},\underbrace{2,\underbrace{1,\ldots,1,}_{a-2}\ldots ,2,\underbrace{1,\ldots,1}_{a-2}}_{v\; {\rm times}}\; ,
	\end{equation}
	where $(u,v)=(a-b_1,m-1-s)$ if $2 \leqslant b_1 \leqslant a$ and $(u,v)=(a-2,m-2-s)$ if $b_1=a+2$.
Furthermore, for such a composition $\bf b$, the coefficient $c_{\bf b}$ is equal to $3$ when $b_1=2$, to $2$ when $3 \leqslant b_1 \leqslant a$ and to $1$ when $b_1=a+2$.
\end{thm}

The binary word associated to the composition ${\bf a}=(\underbrace{a,\ldots,a}_{m})$ is
\begin{equation}
\underbrace{\{0\}_{a-1}1\{0\}_{a-1}1\ldots\{0\}_{a-1}1}_{m \;  {\rm times}}.
\end{equation}
With the notations of Theorem \ref{Th10}, we have
\begin{equation}
\delta([{\bf a}])=\sum_{i=0}^{m-1}\sum_{j=1}^{a-1}\big({\overline{\bf a}}_{ma-ia-j}\boxast {\bf a}_{ia+j}\big)\times\left \{ \substack{2 \quad \text{if } 1\leqslant j \leqslant a-2\;\\3 \quad \text{if } j=a-1\quad} \right. +\sum_{i=1}^{m-1}{\overline{\bf a}}_{ma-ia}\boxast {\bf a}_{ia}\cdot
\end{equation}
For $0 \leqslant i \leqslant m-1$ and $1 \leqslant j \leqslant a-1$, we have
\begin{equation}
\begin{aligned}
{\overline{\bf a}}_{ma-ia-j}&=\Big(\underbrace{1,\ldots,1}_{j},\underbrace{2,\underbrace{1,\ldots,1}_{a-2},\ldots,2,\underbrace{1,\ldots,1}_{a-2}}_{i \; {\rm times}}\Big),\\
{\bf a}_{ia+j}\quad\;&=(a-j,\underbrace{a,\ldots,a}_{m-1-i}),
\end{aligned}
\end{equation}
\begin{equation}\label{EE149}
{\overline{\bf a}}_{ma-ia-j}\boxast {\bf a}_{ia+j}=\Big(a-j+1,(\underbrace{1,\ldots,1}_{j-1},\underbrace{2,\underbrace{1,\ldots,1}_{a-2},\ldots,2,\underbrace{1,\ldots,1}_{a-2}}_{i \; {\rm times}})\ast (\underbrace{a,\ldots,a}_{m-1-i})\Big),
\end{equation}
For $1 \leqslant i \leqslant m-1$, we have
\begin{equation}
\begin{aligned}
{\overline{\bf a}}_{ma-ia}&=\Big(\underbrace{2,\underbrace{1,\ldots,1}_{a-2},\ldots,2,\underbrace{1,\ldots,1}_{a-2}}_{i \; {\rm times}}\Big),\\
{\bf a}_{ia}\;&=(\underbrace{a,\ldots,a}_{m-i}),
\end{aligned}
\end{equation}
\begin{equation}\label{EE151}
{\overline{\bf a}}_{ma-ia}\boxast {\bf a}_{ia}=\Big(a+2,(\underbrace{1,\ldots,1}_{a-2},\underbrace{2,\underbrace{1,\ldots,1}_{a-2},\ldots,2,\underbrace{1,\ldots,1}_{a-2}}_{i-1 \; {\rm times}})\ast (\underbrace{a,\ldots,a}_{m-1-i})\Big).
\end{equation}
We see that the compositions $\bf b$ occurring in the right hand sides of (\ref{EE149}) and (\ref{EE151}) satisfy the conditions $\rm (i)$, $\rm (ii)$, $\rm (iii)$, $\rm (iv)$ of Theorem \ref{Th13}, with $b_1=a-j+1$ in the first case, $b_1=a+2$ in the second case and $s=m-1-i$ in both cases. Moreover each such composition is obtained once and only once. The theorem follows.
\medskip

{\it Example.}- We have
\begin{equation}
\delta([3,3])=3(2,1,2,1)+3(3,3,1)+3(2,1,3)+3(2,4)+2(3,2,1)+2(3,3)+(5,1).
\end{equation}

\medskip
{\it Remark.}- The examples a) and b) given in this section let us think that there is no hope to get a simpler explicit expression of $\delta([{\bf a}])$ than the one stated in Theorem \ref{Th10}.
\vspace{2em}
\section{Linear relations between double tails of multiple zeta values}\label{s4}

\medskip
\subsection{The module $\rm M$ of $\bf Z$-linear relations between double tails of multiple zeta values}\label{s4s1}
\medskip

\myindent Let $\rm M$ denote the submodule of the $\bf Z$-module ${\bf Z}^{(\mathcal{B})}$ consisting of elements $\ell \in {\bf Z}^{(\mathcal{B})}$ such that
\begin{equation}
\forall n \geqslant 0, \quad\zeta(\ell)_{n,n}=0.
\end{equation} 
  As we have seen in Section \ref{s1s12}, $\rm M$ contains the  kernels of the maps $\alpha_k$ defined in Section \ref{s1s11}, and in particular is not equal to $\{0\}$. Our purpose in  this section  is to investigate the structure of $\rm M$.

\medskip

\begin{thm}\label{T3}
 The module $\rm M$  is the kernel of  the $\bf Z$-linear map $\delta: {\bf Z}^{(\mathcal{B})} \rightarrow {\bf Z}^{(\mathcal{A})}$ defined in Section~\ref{s3as1}.
 \end{thm}
 This immediately follows from the remark in Section \ref{s3s2}.
 
 \medskip
\begin{customcor}{1}\label{5c1}
 The $\bf Z$-module $\rm M$  is a graded submodule of ${\bf Z}^{(\mathcal{B})}$ (for the weight ).
 \end{customcor}
 This is a consequence of Theorem $\ref{T3}$, since the $\bf Z$-linear map $\delta: {\bf Z}^{(\mathcal{B})} \rightarrow {\bf Z}^{(\mathcal{A})}$ is graded of degree $0$ (when $ {\bf Z}^{(\mathcal{B})}$ and ${\bf Z}^{(\mathcal{A})}$ are graded by the weight).
 
 \medskip
 {\it Remark.-}  Corollary \ref{5c1} is a non trivial statement. Indeed one conjectures that the submodule of \mbox{${\bf Z}^{(\mathcal{B})}$} consisting of the elements $\ell \in{\bf Z}^{(\mathcal{B})}$ such that $\zeta(\ell)=0$ is graded, but that  is not known at present.
 
\vspace{2em}
\subsection{A description by induction of the homogeneous components of $\rm M$}\label{s4s4}
\medskip

\myindent Let us write ${\rm M}_k$ the homogeneous component of degree $k$ of the graded module $\rm M$, for each integer $k\geqslant 0$. We have ${\rm M}=\bigoplus_{k\geqslant 0}{\rm M}_k$, ${\rm M}_{0}=\{0\}$ and for $k\geqslant 1$,
\begin{equation}
{\rm Ker}(\alpha_k) \;\subset\; {\rm M}_k\;=\;{\rm Ker}(\delta_k)
\end{equation}
by Theorem \ref{T3}, where $\alpha_k:{\bf Z}^{\mathcal{B}_k}\rightarrow \bigoplus_{0 \leqslant k'< k} {\bf Z}^{\mathcal{B}_{k'}}$ has been defined in Section \ref{s1s11}.

\medskip
\begin{thm}\label{Th15}
For each $k\geqslant 1$, we have
\begin{equation}
{\rm M}_k=\alpha_k^{-1}\left(\bigoplus_{0\leqslant k'<k}{\rm M}_{k'}\right)=\bigcap_{0\leqslant k'<k}\alpha_{k',k}^{-1}({\rm M}_{k'}).
\end{equation}
\end{thm}
This indeed follows from the fact, that, for $k \geqslant 1$, we have
\begin{equation}
\mu_k\circ \delta_k=\left(\bigoplus_{0\leqslant k'<k}\delta_{k'}\right)\circ \alpha_k
\end{equation}
by Section \ref{s3as1}, and that $\mu_k$ is injective.

\medskip
\begin{customcor}{1}
For each element $\sum_{{\bf b}\in \mathcal{B}}\lambda_{\bf b}	{\bf b}$ of ${\rm M}$, we have $\sum_{{\bf b}\in \mathcal{B}}\lambda_{\bf b}=0$.
\end{customcor}

It is sufficient to prove this when the element $\ell=\sum_{{\bf b}\in \mathcal{B}}\lambda_{\bf b}	{\bf b}$ is homogeneous, let us say of weight $k$. We do it by induction on $k$. The result is clear when $k=0$. Assume now $k\geqslant 1$. For $0\leqslant k'<k$, we have $\alpha_{k',k}(\ell)\in {\rm M}_{k'}$, hence the sum of the coefficients of $\alpha_{k',k}(\ell)$ is $0$ by the induction hypothesis. Therefore the sum of the coefficients of $\alpha_k(\ell)$~is~$0$. But this sum is equal to $3 \sum_{{\bf b}\in \mathcal{B}}\lambda_{\bf b}$ by  definition of the map $\alpha$.

\medskip
Let $\mathcal{A}^{\{2,3\}}$ denote the set of admissible compositions whose entries  all belongs to $\{2,3\}$. It is expected (but far from being proved) that the real numbers $\zeta({\bf a})$, where ${\bf a}\in \mathcal{A}^{\{2,3\}}$, are $\bf Z$-linearly independent. Hence the elements $\delta([{\bf a}])$, where ${\bf a}\in \mathcal{A}^{\{2,3\}}$, should {\it a fortiori} be ${\bf Z}$-linearly independent. In fact, this latter statement is true. We even have the  following stronger result:

\medskip
\begin{prop}
	Let $\mathcal{A}^{\geqslant 2}$ (resp. $\mathcal{A}^{\leqslant 2}$) denote the set of admissible compositions whose entries are all $\geqslant 2$ (resp. $ \leqslant 2$). The map $\ell \rightarrow \delta([\ell])$ from ${\bf Z}^{(\mathcal{A}^{\geqslant 2})}$  to ${\bf Z}^{(\mathcal{A})}$ (resp. from~${\bf Z}^{(\mathcal{A}^{\leqslant 2})}$  to ${\bf Z}^{(\mathcal{A})}$) is injective.
\end{prop}

Since the involution ${\bf a}\rightarrow \overline{\bf a}$ of $\mathcal{A}$ maps $\mathcal{A}^{\geqslant 2}$ onto $\mathcal{A}^{\leqslant 2}$, it suffices to prove the assertion relative to $\mathcal{A}^{\geqslant 2}$. Since moreover $\delta$ is graded of degree $0$, it suffices to prove that, for each integer $k\geqslant 0$, the linear map $\ell \rightarrow \delta([\ell])$ from ${\bf Z}^{\mathcal{A}_k^{\geqslant 2}}$ to  ${\bf Z}^{\mathcal{A}_k}$ is injective. We argue by induction on $k$. The statement is clear when $k=0$ or $k=1$. We now assume $k \geqslant 2$.

Let $\ell=\sum_{{\bf a}\in \mathcal{A}^{\geqslant 2}}\lambda_{\bf a}{\bf a}$ be an element of ${\bf Z}^{\mathcal{A}_k^{\geqslant 2}}$ such that $\delta([\ell])=0$. For each ${\bf a} \in \mathcal{A}^{\geqslant 2}$, ${\bf a}^{\rm init}$ and ${\bf a}^{\rm mid}$ are of weight $\leqslant k-2$; moreover, if ${\bf a}=(a_1,\ldots,a_r)$, ${\bf a}^{\rm fin}$ is equal to $(a_1-1,a_2,\ldots,a_r)$ and  of weight $k-1$ if $a_1 \geqslant 3$, but it is equal to $(a_2,\ldots,a_r)$ and of weight $ k-2$ if $a_1=2$. Using these facts, one sees that
\begin{equation}
0=\mu_{k-1,k}(\delta([\ell]))=\delta(\alpha_{k-1,k}([\ell]))=\delta(\sum_{\substack{(a_1,\ldots,a_r)\in \mathcal{A}_k^{\geqslant 2}\\ a_1 \geqslant 3}}\lambda_{(a_1,\ldots,a_r)}[a_1-1,\ldots,a_r]).
\end{equation}

From the induction hypothesis, we deduce that $\lambda_{(a_1,\ldots,a_r)}=0$ whenever $a_1 \geqslant 3$. Writing in a similar way that $\mu_{k-2,k}(\delta([\ell]))=0$, we get $\delta([\ell'])=0$, where
\begin{equation}\label{EU174}
\ell'=\sum_{\substack{(a_1,\ldots,a_r)\in \mathcal{A}_k^{\geqslant 2}\\ a_1=2}}\lambda_{(a_1,\ldots,a_r)}(a_2,\ldots,a_r)+\sum_{\substack{(a_1,\ldots,a_r)\in \mathcal{A}_k^{\geqslant 2}\\ a_r=2}}\lambda_{(a_1,\ldots,a_r)}(a_1,\ldots,a_{r-1}).
\end{equation}
By the induction hypothesis we have $\ell'=0$.

For each ${\bf a}=(a_1,\ldots,a_r)\in \mathcal{A}_k^{\geqslant 2}$, let $s({\bf a})$ denote the largest integer $i \leqslant r$ such that $a_1=\ldots=a_i=2$. Assume that $\ell\neq 0$ and then let $\bf a$ be an element in the support of $\ell$ for which $s({\bf a})$ is minimal. If $s({\bf a})<k/2$, $(a_2,\ldots,a_r)$ occurs only once in the first sum (\ref{EU174}), and its coefficient there is $\lambda_{\bf a}$, whereas it does not occur in the second sum. Hence, we  get $\lambda_{\bf a}=0$, which is a contradiction. This implies that we have \smash{$\ell=\lambda_{(\underbrace{2,\dots,2}_{k/2})}(\underbrace{2,\dots,2}_{k/2})$}, $\ell'=2\lambda_{(\underbrace{2,\dots,2}_{k/2})}(\underbrace{2,\dots,2}_{k/2-1})$, and therefore $\lambda_{(\underbrace{2,\dots,2}_{k/2})}=0.$

\vspace{2em}

\subsection{Cases of equality of double tails of multiple zeta values  }\label{s4s5}
\begin{thm}\label{T7}
Let $\bf a$ and $\bf b$ be two admissible compositions such that $\zeta({\bf a})_{n,n}=\zeta({\bf b})_{n,n}$ for all $n \geqslant 0$. Then ${\bf b}$ is equal to $\bf a$ or $\overline{\bf a}$. 
\end{thm}

The hypothesis means that $\delta([{\bf a}])=\delta({[{\bf b}]})$ (Theorem \ref{T3}). Let $k$ be the weight of $\bf a$. By Theorem \ref{Th10}, $\delta([{\bf a}])$ is a non zero homogeneous element of ${\bf Z}^{(\mathcal{A})}$, of weight $k$. The same is then true of $\delta([{\bf b}])$, hence $\bf b$ has the same weight $k$. We shall prove Theorem \ref{T7} by induction on $k$. It is clear when $k=0$. From now on, we assume $k\geqslant 1$.

By Theorem  \ref{Th6}, we have
\begin{equation}
\mu(\delta({\bf a}))=\delta([{\bf a}^{\rm init}])+\delta([{\bf a}^{\rm mid }])+\delta([{\bf a}^{\rm fin }]).
\end{equation}
Each of the three terms in the right-hand side is a non zero homogeneous element of ${\bf Z}^{(\mathcal{A})}$, with non-negative coefficients (Theorem \ref{Th10}). Their respective weights are $|{\bf a}^{\rm init }|$, $|{\bf a}^{\rm mid }|$ and $|{\bf a}^{\rm fin}|$. Hence the set of integers $k'$ such that $0\leqslant k' <k$ and $\mu_{k',k}(\delta({[{\bf a}]}))\neq 0$ is $\{|{\bf a}^{\rm init}|, |{\bf a}^{\rm mid}|, |{\bf a}^{\rm fin}|\}$. Since $\delta([{\bf a}])=\delta([{\bf b}])$, this set is also $\{|{\bf b}^{\rm init}|, |{\bf b}^{\rm mid}|, |{\bf b}^{\rm fin}|\}$. We now distinguish two cases:

\smallskip
\noindent{\it First case:} {\it We have $\mu_{0,k}(\delta([{\bf a}]))=n\varnothing$ with $n \geqslant 2$.}

This happens if and only if ${\bf a}^{\rm init}$ or ${\bf a}^{\rm fin}$ is $\varnothing$, {\it i.e.} if and only if $\bf a$ is equal to either $(k)$ or $\overline{(k)}$. But the same is true then for $\bf b$, since $\delta([{\bf b}])=\delta([{\bf a}])$. Hence $\bf b$ is equal to $\bf a$ or $\overline{\bf a}$.

\smallskip
\noindent{\it Second case:} {\it We are not in the first case.}

Then $|{\bf a}^{\rm mid}|$ is smaller than $|{\bf a}^{\rm init}|$ and $|{\bf a}^{\rm fin}|$, and the same holds for $\bf b$. We therefore have $|{\bf a}^{\rm mid}|=|{\bf b}^{\rm mid}|$, and $\delta([{\bf a}^{\rm mid}])=\delta([{\bf b}^{\rm mid}])$. By the induction hypothesis, ${\bf b}^{\rm mid}$ is equal to ${\bf a}^{\rm mid}$ or $\overline{{\bf a}^{\rm mid}}$. By replacing if needed $\bf b$ by $\overline{\bf b}$, we can assume that ${\bf a}^{\rm mid}={\bf b}^{\rm mid}$ (see Section \ref{s1s9}). Let $r$ and $s$ be the integers $|{\bf a}^{\rm init}|-|{\bf a}^{\rm mid }|$ and $|{\bf a}^{\rm fin}|-|{\bf a}^{\rm mid }|$ respectively. We have $r\geqslant 1$, $s \geqslant 1$ and
\begin{equation}
{\bf w}({\bf a})=0\{1\}_{r-1}{\bf w}({\bf a}^{\rm mid})\{0\}_{s-1}1.
\end{equation}
We now distinguish two subcases:

\smallskip
 {\it First subcase:} {\it We have }$|{\bf a}^{\rm init}|=|{\bf b}^{\rm init}|$.

We then also have $|{\bf a}^{\rm fin}|=|{\bf b}^{\rm fin}|$. Hence
${\bf w}({\bf b})=0\{1\}_{r-1}{\bf w}({\bf b}^{\rm mid})\{0\}_{s-1}1$. This implies ${\bf w}({\bf b})={\bf w}({\bf a})$ since ${\bf w}({\bf b}^{\rm mid})={\bf w}({\bf a}^{\rm mid})$, and therefore ${\bf b}={\bf a}$.

\smallskip
 {\it Second subcase:} {\it We have }$|{\bf a}^{\rm init}|\neq|{\bf b}^{\rm init}|$.
 
 We then have  $|{\bf a}^{\rm init}|=|{\bf b}^{\rm fin}|$, $|{\bf a}^{\rm fin}|=|{\bf b}^{\rm init}|$, and since these integers differ, $\delta([{\bf a}^{\rm init}])=\delta([{\bf b}^{\rm fin}])$. From the induction hypothesis, we deduce that ${\bf b}^{\rm fin}$ is equal to either ${\bf a}^{\rm init}$ or to $\overline{{\bf a}^{\rm init}}$. But we have 
 \begin{equation}
 {\bf w}({\bf b})=0\{1\}_{s-1}{\bf w}({\bf b}^{\rm mid})\{0\}_{r-1}1=0\{1\}_{s-1}{\bf w}({\bf a}^{\rm mid})\{0\}_{r-1}1,
 \end{equation}
  hence ${\bf w}({\bf b}^{\bf fin})={\bf w}({\bf a}^{\rm mid})\{0\}_{r-1}1$, whereas ${\bf w}({\bf a}^{\rm init})=0\{1\}_{r-1}{\bf w}({\bf a}^{\rm mid})$. Equality ${\bf b}^{\rm fin}=\overline{{\bf a}^{\rm init}}$ implies that ${\bf a}^{\rm mid}=\overline{{\bf a}^{\rm mid}}$. Equality ${\bf b}^{\rm fin}={\bf a}^{\rm init}$ implies that $r=2$ and ${\bf w}({\bf a}^{\rm mid})$ is of the form $01\ldots01$, hence again that ${\bf a}^{\rm mid}=\overline{{\bf a}^{\rm mid}}$. Therefore in both cases, we have ${\bf b}=\overline{\bf a}$.

\medskip
{\it Remark.}- Theorem \ref{T7} is a non trivial statement. Indeed one conjectures that the relation $\zeta({\bf a})=\zeta({\bf b})$ implies that $\bf b$ is equal to $\bf a$ or $\overline{\bf a}$, but this is not known at present.

\vspace{2em}
\subsection{Numerical evaluation of the rank of ${\rm M}_k$ for $k\leqslant 16$}\label{s4s6}
\medskip

\myindent We have computed the rank of ${\rm M}_k$ for $k\leqslant 16$. We obtain the following table:
\begin{table}[H]
\centering
\begin{tabular}{|c|c|c|c|c|c|c|c|c|c|c|c|c|c|c|c|c|c|}
\hline
$k$ & 0 & 1 & 2 & 3 & 4 & 5 & 6 & 7 & 8 & 9 & 10 & 11 & 12 & 13 & 14 & 15 & 16 \\
\hline
 ${\rm{rk(M}}_k)$& 0 & 0 & 0 & 0 & 0 & 0 & 1 & 0 & 4 & 2 & 14 & 15 & 52 & 78 & 200 & 350& 789 \\
\hline
\end{tabular} 
\end{table}

Comparing this table with the one obtained in Section \ref{s1s12}, we see that in general, the inclusion ${\rm Ker}(\alpha_k)\subset {\rm Ker}(\delta_k)={\rm M}_k$ is strict. The first dimension for which this occurs is $k=8$.

We do not have at present any closed formula expressing the rank of ${\rm M}_k$ as an explicit function of $k$.

\vspace{2em}
\section{Some further results on the map $\delta$}\label{s5}
\subsection{On the  multiple Ap{\'e}ry-like sum $\sigma(2,\ldots,2)$}\label{s5s1}
We have seen in example 1 of section \ref{s2s7a}, that for each integer $r\geqslant 0$, we have
\begin{equation}\label{EU167}
\sigma(\underbrace{2,\ldots,2}_r)=\frac{\pi^{2r}}{3^{2r}(2r)!}\cdot
\end{equation}

This  formula  suggests that there might exist an element $\ell \in {\bf Z}^{\mathcal{B}_{2r}}$ such that  \mbox{$\delta(\ell)=(\underbrace{2,\ldots,2}_{r})$}, {\it i.e.} such that $\sigma(\underbrace{2,\ldots,2}_{r})_n=\zeta(\ell)_{n,n}$ for all $n\geqslant 0$. One of the difficulty here lies in the fact that this condition does not ensure the unicity of $\ell$, since $\delta$ is not injective. We are nevertheless able to  exhibit such an $\ell$, with a neat  expression:

 \medskip
 \begin{thm}\label{Th17}
 	For every integer $r \geqslant 1$, we have
 	\begin{equation}
 	\delta\Big(\sum_{{\bf a} \in \mathcal{A}_{2r}}(-1)^{{\rm depth}({\bf a})}\;4^{{\rm height}({\bf a})-1}\;[{\bf a}]\Big)=(-1)^{r}3^{2r-1}(\underbrace{2,\ldots,2}_{r}).
 	\end{equation}
 \end{thm}
 
 \medskip 
 Theorem \ref{Th17} can be deduced by the specialization  $t=4$ from the more general \mbox{Theorem~\ref{Th18}}, stated and proved in the next section.
\vspace{2em}
 \subsection{The elements $\ell_{k}(t)$ and their images by $\delta$}
 \begin{thm}\label{Th18}
 	 Let $t$ be an indeterminate. For every even integer $k \geqslant 0$, we have
 	\begin{equation}\label{E175}
 	\delta_{{\bf Z}[t]}\Big(\sum_{{\bf a} \in \mathcal{A}_{k}}(-1)^{{\rm depth}({\bf a})}\;t^{{\rm height}({\bf a})}\;[{\bf a}]\Big)=\sum_{{\bf b}\in \mathcal{A}_{k}^{\rm even}}c_{\bf b}(t){\bf b},
 	\end{equation}
 	where $\delta_{{\bf Z}[t]}:{\bf Z}[t]^{(\mathcal{B})}\rightarrow {\bf Z}[t]^{(\mathcal{A})}$ is the ${\bf Z}[t]$-linear extension of $\delta$, $\mathcal{A}_{k}^{\rm even}$ is the set of admissible compositions of weight $k$  with even entries, and for each  ${\bf b}\in \mathcal{A}_k^{\rm even}$,
 	\begin{equation}\label{E176}
 	c_{\bf b}(t)=(-1)^{s({\bf b})}(t^2-4t)^{{\rm depth}({\bf b})-s({\bf b})}\begin{cases}
 	3t(2t+1)^{s({\bf b})-1}& {\rm if} \quad b_1 = 2,\\
 	(2t+1)^{s({\bf b})} & \; {\rm otherwise},\\
 	\end{cases}
 	\end{equation}
 	where $s({\bf b})$ is the number of entries of $\bf b$ equal to $2$.
 \end{thm}
 
 \medskip
 {\it Remark.}- We recover  Theorem \ref{Th17} of Section \ref{s5s1} as a particular case of Theorem~\ref{Th18} by taking $t=4$ and Theorem \ref{Th9} of Section \ref{S3S3} by taking $t=1$.
 Some other interesting values of $t$ are $-1$, $-2$, $-\frac{1}{2}$. We get for them the following formulas:
 
 \begin{equation}
 \delta(\sum_{ {\bf a}\in \mathcal{A}_k}(-1)^{{\rm depth}({\bf a})-{\rm height}({\bf a})}[{\bf a}])=\sum_{ {\bf b}\in \mathcal{A}_{k}^{\rm even}}5^{{\rm depth}({\bf b})-s({\bf b})}{\bf b}\;\times\begin{cases}
 3 &{\rm if}\quad b_1=2,\\1 & {\rm otherwise},
 \end{cases}
 \end{equation}
 \begin{equation}
 \delta(\sum_{ {\bf a}\in \mathcal{A}_k}(-1)^{{\rm depth}({\bf a})-{\rm height}({\bf a})}2^{{\rm height}({\bf a})}[{\bf a}])=\sum_{ {\bf b}\in \mathcal{A}_{k}^{\rm even}}3^{{\rm depth}({\bf b})}4^{{\rm depth}({\bf b})-s({\bf b})}{\bf b}\;\times\begin{cases}
 2&{\rm if}\quad b_1=2,\\1 & {\rm otherwise},
 \end{cases}
 \end{equation}
 
 \begin{equation}\label{E179}
 \delta_{\bf Q}(\sum_{ {\bf a}\in \mathcal{A}_k}(-1)^{{\rm depth}({\bf a})-{\rm height}({\bf a})}2^{-{\rm height}({\bf a})}[{\bf a}])=\sum\nolimits' \left(\frac{3}{2}\right)^{2{\rm depth}({\bf b})}{\bf b}+\sum\nolimits''\left(\frac{3}{2}\right)^{2{\rm depth}({\bf b})-1}{\bf b},
  \end{equation}
  where in $\sum\nolimits'$ (resp. $\sum\nolimits''$ ), ${\bf b}$ runs over the admissible compositions of weight $k$ whose entries are all even and $\geqslant 4$ (resp. are all even and $\geqslant 4$, except the first one which is $2$).
 \medskip
 
 Theorem \ref{Th18} is clearly true for $k=0$ and $k=2$, since $\delta([\varnothing])=\varnothing$ and $\delta([2])=3(2)$. Its  proof, by induction on $k$, is  similar to the proof of Theorem \ref{Th9}. It follows from the following lemma, where
  \begin{equation}
 \ell_{p}(t)=\sum_{{\bf a}\in \mathcal{A}_p}(-1)^{{\rm depth}({\bf a})}t^{{\rm height}({\bf a})}{\bf a}
 \end{equation}
 for every integer $p \geqslant 0$.
 
 \medskip
 
 \begin{lem}
 	For every even integer $k\geqslant 4$, we have
 	\begin{equation}
 	\alpha_{{\bf Z}[t]}([\ell_k(t)])=(t^2-4t)\sum_{\substack{0 \leqslant q \leqslant k-4\\q \;{\rm even}}}[\ell_q(t)]-(2t+1)[\ell_{k-2}(t)].
 	\end{equation}
 \end{lem}
 
 We have $\ell_0(t)=\varnothing$, $\ell_1(t)=0$, $\ell_0(t)^{\rm init}=\ell_0(t)$ and $\ell_p(t)=(-1)^p\ell_p(t)$ for $p\geqslant 0$. Note that, if ${\bf a}=(a_1,\ldots,a_r)$ is a non-empty admissible                        composition of height $h$, ${\bf a}^{\rm init}=(a_1,\ldots,a_{r-1})$  is of height $h-1$ if $a_r \geqslant 2$, and of height $h$ if $a_r=1$. Therefore we have, for each integer $p \geqslant 2$,
 \begin{equation}
 \ell_p(t)^{\rm init}=-t\sum_{q=0}^{p-2}\ell_q(t)-\ell_{p-1}(t).
 \end{equation}
 
 Let now $k$ be an even integer $\geqslant 4$. We have
 \begin{equation}
 \begin{aligned}
 \ell_k(t)^{\rm init}&= -t\left(\ell_0(t)+\sum_{p=2}^{k-2}\ell_p(t)\right)-\ell_{k-1}(t),
 \end{aligned}
 \end{equation}
 \begin{equation}
 \begin{aligned}
 \ell_k(t)^{\rm fin}&=\overline{\overline{\ell_k(t)}^{\rm init}}=\overline{{\ell_k(t)}^{\rm init}}=-t\left(\ell_0(t)+\sum_{p=2}^{k-2}(-1)^p\ell_p(t)\right)+\ell_{k-1}(t),
\end{aligned}
\end{equation}
\begin{equation}
\begin{aligned}
 \ell_k(t)^{\rm mid}&=(\ell_k(t)^{\rm fin})^{\rm init}=-t\ell_0(t)+t^2\sum_{p=2}^{k-2}(-1)^p
\left(\ell_0(t)+\sum_{q=2}^{p-2}\ell_q(t)\right)\\& +t\sum_{p=2}^{k-2}(-1)^p\ell_{p-1}(t)-t\left(\ell_0(t)+\sum_{p=2}^{k-3}\ell_p(t)\right)-\ell_{k-2}(t) \\
&=(t^2-2t)\sum_{\substack{0 \leqslant q \leqslant k-4\\ q \; {\rm even}}}\ell_q(t)-\ell_{k-2}(t),
\end{aligned}
 \end{equation}
 and hence
 \begin{equation}
 \ell_k(t)^{\rm init}+\ell_k(t)^{\rm mid}+\ell_k(t)^{\rm fin}=(t^2-4t)\sum_{\substack{0 \leqslant q \leqslant k-4\\ q \; {\rm even}}}\ell_q(t)-(2t+1)\ell_{k-2}(t).
 \end{equation}
 The lemma follows.
 \medskip
 
 Let $k \geqslant 0$ be an even integer. One deduces from Theorem \ref{Th18} that we have
 \begin{equation}\label{EU199}
 \sum_{\bf a \in \mathcal{A}_k}(-1)^{{\rm depth}({\bf a})}\zeta({\bf a})t^{{\rm height}({\bf a})}=\sum_{{\bf b}\in \mathcal{A}_k^{\rm even}}c_{\bf b}(t)\sigma({\bf b })
 \end{equation}
 where the polynomials $c_{\bf b}(t)$ are given by the formula (\ref{E176}). Moreover, a generating series for the left-hand sides of (\ref{EU199}), when $k$ varies, has been gin by Y. Ohno and D.~Zagier~\cite{OZ}:
 
 \medskip
 \begin{lem}\label{L11}
 	For each even $k \geqslant 2$,
 	\begin{equation}
 	\sum_{ {\bf a}\in \mathcal{A}_k}(-1)^{{\rm depth}({\bf a})}\zeta({\bf a})t^{{\rm height}({\bf a})}
 	\end{equation}
 	is the coefficient of $x^k$ in the Taylor expansion of
 	\begin{equation}
 	\frac{t}{t-1}\prod_{m=1}^{\infty}\frac{m^2-tx^2}{m^2-x^2}=\frac{t^{1/2}{\rm sin}( \pi t^{1/2}x)}{(t-1){\rm sin}(\pi x)}\cdot
 	\end{equation} 
 	
 \end{lem}

 Indeed, Y. Ohno and D. Zagier defined  a formal power series $\Phi_0(x,y,z)$ in three indeterminates $x$, $y$, $z$  by
 \begin{equation}
 \Phi_0(x,y,z)=\sum_{\substack{{\bf a}\in \mathcal{A}\\ {\bf a} \neq \varnothing}}\zeta({\bf a})\;x^{{\rm weight}({\bf a})-{\rm depth}({\bf a})-{\rm height}({\bf a})}\;y^{{\rm depth}({\bf a})-{\rm height}({\bf a})}\;z^{{\rm height}({\bf a})-1}
 \end{equation}
 and  proved that it satisfies the identity
 \begin{equation}
 1-(xy-z)\Phi_0(x,y,z)=\prod_{m=1}^{\infty}\left(1-\frac{xy-z}{(m-x)(m-y)}\right)\cdot
 \end{equation}
 Therefore, for $k$ even $\geqslant 2$, $\sum_{ {\bf a}\in \mathcal{A}_k}(-1)^{{\rm depth}({\bf a})}\zeta({\bf a})t^{{\rm height}({\bf a})}$ is the coefficient of $x^k$ in $-tx^2\;\Phi_0(x,-x,-tx^2)$, and we have 
 \begin{equation}
 1+(1-t)\;x^2\;\Phi_0(x,-x,-tx^2)=\prod_{m=1}^{\infty}\frac{m^2-tx^2}{m^2-x^2}\cdot
 \end{equation}
 It follows that
 $\sum_{ {\bf a}\in \mathcal{A}_k}(-1)^{{\rm depth}({\bf a})}\zeta({\bf a})t^{{\rm height}({\bf a})}$ is the coefficient of $x^k$ in the Taylor expansion of
 \begin{equation}
 \frac{t}{t-1}\prod_{m=1}^{\infty}\frac{m^2-tx^2}{m^2-x^2}=\frac{t}{t-1}\cdot\frac{{\rm sin}(\pi t^{1/2}x)}{\pi t^{1/2}x}\cdot\frac{\pi x}{{\rm sin}(\pi x)}=\frac{t^{1/2}}{t-1}\frac{{\rm sin}(\pi t^{1/2}x)}{{\rm sin}(\pi x)}\cdot
 \end{equation}
\vspace{2em} 
\subsection{The elements $\ell_{k,h}$ and their images by $\delta$}
 \myindent Let $k$ and $h$ be integers such that $0 \leqslant h \leqslant k/2$. Let $\mathcal{A}_{k,h}$ denote the set of admissible compositions of weight $k$ and height $h$, and let $\ell_{k,h}$ denote the element of ${\bf Z}^{\mathcal{A}_{k,h}}$defined by
 \begin{equation}
 l_{k,h}=\sum_{ {\bf a}\in \mathcal{A}_{k,h}}(-1)^{{\rm depth}({\bf a})}{\bf a}.
 \end{equation}
 The set $\mathcal{A}_{k,h}$ is stable by ${\bf a}\rightarrow \overline{\bf a}$, and we have ${\rm depth}({\bf a})+{\rm depth}(\overline{\bf a})=k$ and $[{\bf a}]=[\overline{\bf a}]$ for ${\bf a}\in \mathcal{A}_{k,h}$, hence $[\ell_{k,h}]$ vanishes when $k$ is odd.
 
 Therefore, we assume from now on $k$ to be even. An explicit formula for $\delta([\ell_{k,h}])$ can then be deduced from Theorem \ref{Th18}: we have
 \begin{equation}
 \delta([\ell_{k,h}])=\sum_{{\bf b}\in \mathcal{A}_k^{\rm even}}c_{{\bf b},h}{\bf b},
 \end{equation}
 where $c_{{\bf b},h}$ is the coefficient of $t^h$ in the polynomial $c_{\bf b}(t)$ given by the formula (\ref{E176}).

 \medskip
 \begin{lem}\label{Le11}
 	Let $k$ be an even integer $\geqslant 2$. The elements $\delta([\ell_{k,h}])$ of ${\bf Z}^{\mathcal{A}_k^{\rm even}}$, where \mbox{$1\leqslant h\leqslant k/2$}, are $\bf Z$-linearly independent.
 \end{lem}
 
 The rank $r$ over $\bf Z$ of the family $(\delta([\ell_{k,h}]))_{1\leqslant h\leqslant k/2}$ of  elements ${\bf Z}^{\mathcal{A}_k^{\rm even}}$ is at most $k/2$. It is equal to the rank of the matrix with integer coefficients.
 \begin{equation}
 (c_{{\bf b},h})_{{\bf b}\in \mathcal{A}_k^{\rm even},1 \leqslant h \leqslant k/2}.
 \end{equation}
 Since the $c_{{\bf b},h}$, for given ${\bf b}\in \mathcal{A}_k^{\rm even}$ and $1 \leqslant h \leqslant k/2$, are the coefficients of the polynomial $c_{\bf b}(t)$, which has no constant term and is of of degree at most $k/2$, $r$ is also the rank over $\bf Z$ of the family of polynomials
 \begin{equation}
 (c_{\bf b}(t))_{{\bf b}\in \mathcal{A}_k^{\rm even}}.
 \end{equation}
 But, for $\bf b$ of the form $(k-2m,\underbrace{2,\ldots,2}_m)$, with $0 \leqslant m \leqslant \frac{k}{2}-1$, we have
 \begin{equation}
 c_{\bf b}(t)=\begin{cases}
 (-1)^m(t^2-4t)(2t+1)^m & \text{ if } m<\frac{k}{2}-1,\\
 3(-1)^{k/2}t(2t+1)^{\frac{k}{2}-1} & \text{ if }m=\frac{k}{2}-1.
 \end{cases}
 \end{equation}
 These $k/2$ polynomials being $\bf Z$-linearly independent, we have $r=k/2$.
 
 \medskip
 {\it Remark }1. Let $\rm N$ denote the submodule of ${\bf Z}^{\mathcal{A}_k^{\rm even}}$ consisting of the elements $\sum_{ {\bf b}\in \mathcal{A}_k^{\rm even}}\lambda_{\bf b}{\bf b}$ such that $\lambda_{\bf b}=0$ for all $\bf b$ of the form $(k-2m,\underbrace{2,\ldots,2}_m)$ with $0 \leqslant m\leqslant\frac{k}{2}-1$. It is a direct factor of corank $k/2$ of ${\bf Z}^{\mathcal{A}_k^{\rm even}}$. The proof of Lemma \ref{Le11} shows in fact that, if $\rm L$ is the submodule of ${\bf Z}^{\mathcal{A}_k^{\rm even}}$ generated by the elements $\delta([\ell_{k,h}])$, with $1 \leqslant h \leqslant k/2$, then ${\rm L} \cap {\rm N} =\{0\}$, and hence  ${\rm L}+{\rm N}$ is of finite index in ${\bf Z}^{\mathcal{A}_k^{\rm even}}$.
 
 \medskip
 
 By Lemma \ref{Le11}, the rank of  ${\rm Im}(\delta)\cap {\bf Z}^{\mathcal{A}_k^{\rm even}}$ is at least $k/2$ for each even integer $k\geqslant 2$. We conjecture the following:
 
 \medskip
 \begin{coj}\label{Cc2}
 	For each even integer $k \geqslant 2$, the rank of ${\rm Im}(\delta)\cap {\bf Z}^{\mathcal{A}_k^{\rm even}}$ is equal to $k/2$.
 \end{coj}
 
 We have checked that this conjecture indeed holds when $k \leqslant 16$.
 
 \medskip
 {\it Remark} 2. With the notations of Remark 1, Conjecture \ref{Cc2} can be restated as in any of the two following equivalent forms: 
 \\a) {\it We have} ${\rm N}\cap {\rm Im}(\delta)=\{0\}$.
 \\ b) { \it Any element of   ${\rm Im}(\delta)\cap {\bf Z}^{\mathcal{A}_k^{\rm even}}$ is a $\bf Q$-linear combination of the elements $\delta([\ell_{k,h}])$, with $1\leqslant h\leqslant k/2$.}

 Note that this $\bf Q$-linear combination is not always  $\bf Z$-linear: in weight $6$ for example, we have $\delta(2[4,2]+[3,3]-[3,2,1]-[1,2,3])=\frac{1}{3}\delta([\ell_{6,1}]+[\ell_{6,2}]+[\ell_{6,3}])$.
 
 \medskip
 
 Let $k$ be an even integer $\geqslant 2$. It follows from  Lemma \ref{L11} that the real number $\zeta(\ell_{k,h})$, which is equal to $\sigma(\delta([\ell_{k,h}]))$, belongs to $\pi^k {\bf Q}$ for $0 \leqslant h \leqslant k/2$. We conjecture the following:
 \medskip
 
 \begin{coj}\label{Co2a}
 	Any element $\ell \in {\bf Z}^{\mathcal{A}_k^{\rm even}}$ such that $\sigma(\ell)\in \pi^k{\bf Q}$ is a $\bf Q$-linear combination of the $\delta([\ell_{k,h}])$, where $1 \leqslant  h \leqslant k/2$.
 \end{coj}
 
 Conjecture \ref{Co2a} is supported by unsuccessful searches by PSLQ algorithm of elements $\ell$ such that $\sigma(\ell)\in \pi^k{\bf Q}$ and not of this form, for $k \leqslant 14$.
 
 \medskip
 {\it Remarks.} 3) With the notations of Remark 1, Conjecture \ref{Co2a} can be restated as follows: {\it If $\ell$ is an element of $\rm N$ such that $\sigma(\ell)\in \pi^k{\bf Q}$, then $\ell=0$.} Or equivalently: {\it The elements of $\mathcal{A}_k^{\rm even}$ not of the form $(2m,\underbrace{2,\ldots,2}_{k/2-m})$ with $2\leqslant m \leqslant k/2$, are $\bf Q$-linearly independent}.
 
 \medskip
  4) With the notations of Remark 1, Conjecture \ref{Cc2} is a statement of combinatorial nature, and might  be tractable. But Conjecture \ref{Co2a} is a much deeper statement about the $\bf Q$-linear independence of some Ap{\'e}ry-like sums, and looks out of reach at present.
 
 \medskip
  5) One could be tempted to generalise both  Conjecture \ref{Cc2} and \ref{Co2a} by asking whether any element $\ell \in \mathcal{A}_k^{\rm even}$ such that $\sigma(\ell)$ is a $\bf Q$-linear combination of multiple zeta values of weight $k$,  is itself a $\bf Q$-linear combination of the $\delta([\ell_{k,h}])$, with $1\leqslant h\leqslant k/2$. However, this is probably false: when $k=8$ for example, $\ell=18\;(2,6)+65\;(4,4)+12\;(2,2,4)$ is a non zero element of the module $\rm N$ introduced in Remark 1, hence is not a $\bf Q$-linear combination of the elements $\ell_{8,h}$ for $1\leqslant h\leqslant 4$, but PSLQ algorithm suggests that $\sigma(\ell)$ is equal to
 \begin{equation}
 \frac{16}{825}\left(\frac{1593337}{240}\zeta(2,2,2,2)-747\zeta(3,3,2)-818\zeta(3,2,3)-842\zeta(2,3,3)\right),
 \end{equation}
 an equality that we checked up to $1000$ digits, but for which we have no proof.
\vspace{2em}
\subsection{Finding a pattern in some identities of D. Bailey, J. Borwein, D. Bradley}
\myindent  D. Bailey, J. Borwein, D. Bradley give at the end of their paper \cite{BBB} some formulas that generalize the identity $\zeta(4)=\frac{36}{17}\sigma(4)$. They discovered them experimentally. With our notations, these formulas can be written as
\begin{equation}\label{E194}
\begin{aligned}
\zeta(6)&=\frac{36 \cdot 8}{163}\Big(\sigma(6)+\frac{3}{2}\sigma(2,4)\Big),\\
\zeta(8)&=\frac{36 \cdot 64}{1373}\Big(\sigma(8)+\frac{9}{4}\sigma(4,4)+\frac{3}{2}\sigma(2,6)\Big),\\
\zeta(10)&=\frac{36\cdot 512}{11143}\Big(\sigma(10)+\frac{9}{4}\sigma(6,4)+\frac{3}{2}\sigma(2,8)+\frac{9}{4}\sigma(4,6)+\frac{27}{8} \sigma(2,4,4)\Big).
\end{aligned}
\end{equation}
They conclude their paper by saying \textquotedblleft {\it  This pattern is not fruitful: the pattern stops at $n=10$}.\textquotedblright

Our purpose in this section is to show that actually such a pattern exists and can be derived from relation (\ref{E179}). Indeed  relation (\ref{E179}) implies that, for each even weight $k \geqslant 2$, we have
\begin{equation}\label{E195}
\sum_{{\bf a}\in \mathcal{A}_k}(-1)^{{\rm depth}({\bf a})}\left(\frac{-1}{2}\right)^{{\rm height}({\bf a})}\zeta({\bf a})=\sum\nolimits '\left(\frac{3}{2}\right)^{2\;{\rm depth}({\bf b})}\sigma({\bf b})+\sum\nolimits''\left(\frac{3}{2}\right)^{2\;{\rm depth}({\bf b})-1}\sigma({\bf b}),
\end{equation}
where in $\sum\nolimits'$ (resp. $\sum\nolimits''$ ), ${\bf b}$ runs over the admissible compositions of weight $k$ whose entries are all even and $\geqslant 4$ (resp. are all even and $\geqslant 4$, except the first one which is $2$). The right-hand side of (\ref{E195}) is equal to $\frac{9}{4}$ times
\begin{equation}
\begin{aligned}
&\sigma(4) \qquad\qquad &{\rm when} \; k=4,\\
&\sigma(6)+\frac{3}{2}\sigma(2,4)\qquad\qquad &{\rm when} \; k=6,\\
&\sigma(8)+\frac{9}{4}\sigma(4,4)+\frac{3}{2}\sigma(2,6)\qquad\qquad &{\rm when} \; k=8,\\
&\sigma(8)+\frac{9}{4}\sigma(6,4)+\frac{3}{2}\sigma(2,8)+\frac{9}{4}\sigma(4,6)+\frac{27}{8}\sigma(2,4,4)\qquad &{\rm when} \; k=10,\\
\end{aligned}
\end{equation}
thus providing the $\bf Q$-linear combinations of $\sigma$ values occuring in (\ref{E194}). On the other hand, we know by Lemma \ref{L11} that the left-hand side of  (\ref{E195}) is the coefficient if $x^k$ in the Taylor series expansion at $0$ of 
\begin{equation}
\frac{1}{3}\;\prod_{m=1}^{\infty}\frac{m^2+\frac{x^2}{2}}{m^2-x^2}\;=\;\frac{1}{3}\;\frac{{\rm sinh}\frac{\pi x}{\sqrt{2}}}{\frac{\pi x}{\sqrt{2}}}\;\;\frac{\pi x}{{\rm sin}\pi x}.
\end{equation}
Since by formula (\ref{EE117}), we have
\begin{equation}
\begin{aligned}
\frac{\pi x}{{\rm sin}\pi x}&=1+2\sum_{k \; {\rm even} \geqslant 2}\left(\sum_{n=1}^{\infty}\frac{(-1)^{n-1}}{n^k}\right)x^k=2\sum_{k \; {\rm even} \geqslant 2}(1-2^{1-k})\zeta(k)x^k\\
&=\sum_{k \; {\rm even} \geqslant 0}\frac{(2^k-2)\pi^k(-1)^{k/2-1}B_k}{k!}x^k,
\end{aligned}
\end{equation}
and on the other hand 
\begin{equation}
\frac{{\rm sinh}\frac{\pi x}{\sqrt{2}}}{\frac{\pi x}{\sqrt{2}}}=\sum_{k \; {\rm even} \geqslant 0}\frac{2^{-k/2}\pi^k x^k}{(k+1)!},
\end{equation}
we see that the left-hand side of (\ref{E195}) is equal to
\begin{equation}
\begin{aligned}
&\frac{1}{3}\Big(\sum_{ \substack{p+q=k\\p,q\; {\rm even}\geqslant 0}}\frac{(2^p-2)(-1)^{p/2-1}B_p}{p!}\frac{2^{-q/2}}{(q+1)!}\Big)\pi^k\\&=(-1)^{k/2-1}\frac{2^{1-k}k!}{3B_k}\sum_{\substack{p+q=k\\p,q\; {\rm even} \geqslant 0}}\frac{(2^p-2)(-1)^{p/2-1}B_p2^{-q/2}}{p!(q+1)!}\zeta(k).
\end{aligned}
\end{equation}
The numerical coefficient of $\zeta(k)$ in this last expression is equal to $\frac{17}{2^4}$ for $k=4$, $\frac{163}{2^7}$ for $k=6$, $\frac{1373}{2^{10}}$ for $k=8$, $\frac{11143}{2^{13}}$ for $k=10$. It is $\frac{61835987}{2^{16}\cdot 691}$ for $k=12$. Hence identities (\ref{E194}) follow from these considerations, and (\ref{E195}) extends their pattern to all even weight $\geqslant 2$.

\medskip
{\it Remark.}- It seems that D. Bailey, J. Borwein and D. Bradley found their identities (\ref{E194}) by searching experimentally identities of the form $\zeta(k)=\sum_{ {\bf b}\in \mathcal{A}_k^{\rm even}}\lambda_{\bf b}\sigma({\bf b})$ with rational coefficients $\lambda_{\bf b}$, involving $\sigma(k)$ ({$i.e.$} with $\lambda_{(k)\neq 0}$), and with the least possible number of non zero coefficients. Many such identities can be written, for example by taking $\lambda_{\bf b}$ proportional to $c_{\bf b}(t)$, where $t$ is a rational number distinct from $0$ and $4$ and the $c_{\bf b}(t)$ are given by the formula (\ref{E176}).  Now, among the identities of this form, the one with least number of non zero coefficients corresponds to the case when $t=-\frac{1}{2}$, as is seen on  formula (\ref{E176}). This indicates why  formula (\ref{E179}) yields exactly the identities~(\ref{E194}).
\vspace{2 em}

\subsection{A remark on the matrix of $\delta_k$}
The matrix of the $\bf Z$-linear map $\delta_k:{\bf Z}^{\mathcal{B}_k} \rightarrow {\bf Z}^{\mathcal{A}_k}$ in the canonical bases ${\mathcal{B}_k}$ and ${\mathcal{A}_k}$ of ${\bf Z}^{\mathcal{B}_k}$ and ${\bf Z}^{\mathcal{A}_k}$  is highly complicated. However, when $k$ is even  there is a big submatrix $\Delta_k$ which is nicely behaved: it is the one obtained by keeping only the columns indexed by the subset $\mathcal{B}_k^{\rm sd}$ of $\mathcal{B}_k$ consisting of the elements $[{\bf a}]$, where $\bf a$ is a self-dual admissible composition of weight $k$ ({\it i.e. }$\overline{\bf a}={\bf a}$), and by keeping only the lines indexed by the subset $\mathcal{A}_k^{\rm even}$ of $\mathcal{A}_k$ consisting of compositions with even entries.

Note that $\mathcal{B}_k^{\rm sd}$ and $\mathcal{A}_k^{\rm even}$ have both $2^{m(k)}$ elements, where $m(k)=\frac{k}{2}-1$ if $k \geqslant 2$ and $m(0)=0$. We want to consider  $\Delta_k$ as a square matrix  of size $2^{m(k)}\times2^{m(k)}$ with lines and columns  both indexed by $\{0,1,\ldots,2^{m(k)}-1\}$. For this purpose, we need to specify the bijections
\begin{equation}
\varphi: \{0,1,\ldots,2^{m(k)}-1\}\rightarrow \mathcal{A}_k^{\rm even}, \quad \psi: \{0,1,\ldots,2^{m(k)}-1\}\rightarrow \mathcal{B}_k^{\rm sd}
\end{equation}
indexing the lines and columns. For each integer $i \in \{0,1,\ldots,2^{m(k)}-1\}$, let $w(i)$ be the binary word of weight $k/2$ expressing of $i$ in base $2$: it starts by $0$ if $k\geqslant 2$. Then $\psi(i)$ is the class $[{\bf a}]$ of the self-dual composition $\bf a$ whose associated binary word is $w(i)\overline{w(i)}$, and $\varphi(i)$ is equal to $(2b_1,\ldots,2b_r)$, where $(b_1,\ldots,b_r)$ is the composition whose associated binary word is $\overline{w(i)}$.

With these conventions, the matrices $\Delta_k$, for even $\leqslant 8$, are given by:

\begin{equation}
\begin{aligned}
\Delta_0&=\left({\begin{array}{c} 1 \end{array}}\right)\qquad\qquad\qquad\qquad\quad &{\rm when}\; k=0,
\end{aligned}
\end{equation}
\begin{equation}
\begin{aligned}
\Delta_2&=\left({\begin{array}{c} 3 \end{array}}\right)\qquad\qquad\qquad\qquad\quad &{\rm when}\; k=2,
\end{aligned}
\end{equation}
\begin{equation}
\begin{aligned}
\Delta_4&=\left({\begin{array}{cc} 3 & 6\\ 0 & 1 \end{array}}\right)\qquad\qquad\qquad\quad &{\rm when}\; k=4,
\end{aligned}
\end{equation}
\begin{equation}
\begin{aligned}
\Delta_6&=\left({\begin{array}{cccc} 3 & 6& 12 & 0\\ 0 & 1 & 2 & 0 \\ 0 & 0 & 3 & 0\\ 0& 0 & 0 & 1 \end{array}}\right)\qquad\qquad &{\rm when}\; k=6,
\end{aligned}
\end{equation}
\begin{equation}
\begin{aligned}
\Delta_8&=\left({\begin{array}{cccc:cc:c:c} 3 & 6& 12 & 0 & 6 & 24 & 0 &0\\ 0 & 1 & 2 & 0 & 0 & 4 & 0 & 0\\ 0 & 0 & 3 & 0 & 0 & 6 & 0 & 0\\ 0& 0 & 0 & 1 & 0 & 0 & 0 & 0\\ \hdashline 0 & 0 & 0 & 0& 3 & 6 & 0 & 0\\ 0 & 0 & 0 & 0 &0 & 1 & 0 & 4\\ \hdashline 0 & 0 & 0 & 0& 0 & 0 & 3 & 0\\ \hdashline 0 & 0 & 0 & 0 & 0 &0 & 0 & 1 \end{array}}\right)\quad &{\rm when}\; k=8.
\end{aligned}
\end{equation}

In the case $k=8$ for example, the columns of the matrix are indexed by the elements $[5,1,1,1]$, $[4,2,1,1]$, $[3,2,2,1]$, $[3,1,3,1]$, $[2,3,1,2]$, $[2,2,2,2]$, $[2,1,2,3]$, $[2,1,1,4]$ of $\mathcal{B}_8^{\rm sd}$, whereas  the lines are indexed  by the elements $(2,2,2,2)$, $(4,2,2)$, $(2,4,2)$, $(6,2)$, $(2,2,4)$, $(4,4)$, $(2,6)$, $(8)$ of $\mathcal{A}_8^{\rm even}$. The dotted lines in the matrix $\Delta_8$ materialise a block decomposition from which some pattern can be observed. More details about this pattern are given in Proposition~\ref{P7} below.
\medskip

{\it Remarks }1) The matrix $\Delta_k$ is the matrix of the composed linear map
\begin{equation}
{\bf Z}^{\mathcal{B}_k^{\rm sd}}\hookrightarrow {\bf Z}^{\mathcal{B}_k} \xrightarrow{\delta} {\bf Z}^{\mathcal{A}_k} \twoheadrightarrow {\bf Z}^{\mathcal{A}_k^{\rm even}}
\end{equation}
where the first map is the canonical injection and the last the canonical surjection.

2) The ordering of $\mathcal{B}_k^{\rm sd}$ corresponding by $\psi$ to the natural ordering of $\{0,\ldots, 2^{m(k)}-1\} $ is the descending lexicographical ordering.
\medskip

\begin{prop}\label{P7}
	Let $k$ be an even integer $\geqslant 4$. Let us write $\Delta_k$ as a block matrix $\big(U_{i,j}\big)_{\substack{1\leqslant i\leqslant k/2\\ 1\leqslant j \leqslant k/2}}$ where $U_{i,j}$ is a matrix of size $2^{m(k-2i)}\times 2^{m(k-2j)}$. Then $U_{i,i}=\Delta_{k-2i}$ and
	$U_{i,j}=0$ if j is neither $i$ or $2i$ . In particular, $\Delta_k$ is an upper triangular matrix  whose diagonal elements are alternatively equal to $3$ and $1$. 
\end{prop}

The columns of $U_{i,j}$ correspond to the columns of $\Delta_k$ indexed by the elements of $\mathcal{B}_k$ of the form $[{\bf a}]$, where $\bf a$ is a self-dual composition of weight $k$ with last entry $j$. The lines of $U_{i,j}$ correspond to the lines of $\Delta_k$ indexed by the even compositions of weight $k$ with  last entry $2i$. When $\bf a$ is a self-dual composition of weight $k$ with last entry $j$, ${\bf a}^{\rm init }$ and ${\bf a}^{\rm fin }$ are of  weight $k-j$ and ${\bf a}^{\rm mid }$ is of weight $k-2j$. Hence $\mu_{k-2j,k}(\delta({\bf a}))=\delta({\bf a}^{\rm mid })$ and $\mu_{k',k}(\delta({\bf a}))=0$ if $k'$ is different from $k-j$ and $k-2j$. The proposition follows.

\medskip

\begin{customcor}{1}
	The restriction of the map $\delta$ to ${\bf Z}^{\mathcal{B}^{\rm sd}}$ is injective.
\end{customcor}

Since $\delta$ is graded, it suffices to prove that the restriction of $\delta$ to ${\bf Z}^{\mathcal{B}_k^{\rm sd}}$ is injective for all even integers $k\geqslant 0$. But this follows from Proposition \ref{P7}.

\bigskip
{\bf Acknowledgement :}  I am deeply grateful to Prof.~J.~Oesterl\'{e} of the University of Paris~\rom{6} for guiding me through this work. The bulk of the work described here was carried out at  the University of Paris \rom{6}  with leave from the Institute of Mathematical Sciences, Chennai and support from the IRSES Moduli project. I wish to extend my heartfelt thanks these institutions for their facilities and to Prof. P. Philippon of the University of Paris~\rom{6} for his support throughout this work. The final draft of this paper was prepared  during a visit to the Harish-Chandra Research Institute, Allahabad. The author is an INSPIRE fellow  at the Chennai Mathematical Institute.

\Addresses

\end{document}